\documentclass{birkmult}

\usepackage[all]{xy}

\theoremstyle{plain}
  \newtheorem{theorem}{Theorem}[section]
  \newtheorem{corollary}[theorem]{Corollary}
  \newtheorem{lemma}[theorem]{Lemma}
  \newtheorem{proposition}[theorem]{Proposition}
  \newtheorem{conjecture}[theorem]{Conjecture}
\theoremstyle{definition}
  \newtheorem{definition}[theorem]{Definition}
  \newtheorem{example}[theorem]{Example}
  \newtheorem{remark}[theorem]{Remark}

\DeclareMathOperator{\rank}{rank}
\DeclareMathOperator{\GL}{GL}

\DeclareMathOperator{\ind}{ind}
\DeclareMathOperator{\Det}{Det}

\DeclareMathOperator{\IND}{IND}

\DeclareMathOperator{\Int}{Int}
\DeclareMathOperator{\grad}{grad}

\DeclareMathOperator{\id}{id}

\DeclareMathOperator{\cs}{cs}

\DeclareMathOperator*\LIM{LIM}
\newcommand\ii{\mathbf i}
\newcommand\img{\operatorname{img}}
\newcommand\inc{\operatorname{inc}}
\newcommand\ad{\operatorname{ad}}
\newcommand\Cl{\operatorname{Cl}}
\newcommand\Aut{\operatorname{Aut}}
\newcommand\ec{\operatorname{e}}
\newcommand\tr{\operatorname{tr}}
\newcommand\hol{\operatorname{hol}}
\newcommand\eend{\operatorname{end}}
\newcommand\str{\operatorname{str}}
\newcommand\sdet{\operatorname{sdet}}
\newcommand\const{\operatorname{const}}
\newcommand\rk{\operatorname{rk}}
\newcommand\gl{\mathfrak{gl}}
\newcommand\Sym{\operatorname{Sym}}
\newcommand{\Z}{\mathbb Z}
\newcommand{\Q}{\mathbb Q}
\newcommand{\N}{\mathbb N}

\newcommand{\R}{\mathbb R}
\newcommand{\C}{\mathbb C}

\newcommand{\Eul}{\mathfrak{Eul}}

\newcommand{\e}{\mathfrak e}
\newcommand\comb{\text{\rm comb}}

\newcommand{\can}{\text{\rm can}}

\newcommand{\an}{\text{\rm an}}

\newcommand{\sing}{\text{\rm sing}}
\newcommand\even{\text{\rm even}}
\newcommand\odd{\text{\rm odd}}

\newcommand\itemref[1]{(\ref{#1})}
\renewcommand{\theenumi}{\roman{enumi}}
\renewcommand{\labelenumi}{(\theenumi)}

\begin{document}

\title{Complex valued Ray--Singer torsion}

\author{Dan Burghelea}

\address{Dept. of Mathematics,
         The Ohio State University,\\ 
         231 West 18th Avenue,
         Columbus, OH 43210, USA.}

\email{burghele@mps.ohio-state.edu}

\author{Stefan Haller}

\address{Department of Mathematics,
         University of Vienna,\\
         Nordbergstra{\ss}e 15,
         A-1090 Vienna, Austria.} 

\email{stefan.haller@univie.ac.at}

\thanks{Part of this work was done while the second author enjoyed the
        hospitality of the Max Planck Institute for Mathematics in Bonn.
        During the final stage of preparation of this manuscript the second author
        was partially supported by the Austrian Science Fund (FWF), grant P19392-N13}



\date{\today}

\begin{abstract}
In the spirit of Ray and Singer we define a complex valued analytic torsion using
non-selfadjoint Laplacians. We establish an anomaly formula which permits
to turn this into a topological invariant. Conjecturally this analytically
defined invariant computes the complex valued Reidemeister torsion, including its phase.
We establish this conjecture in some non-trivial situations.
\end{abstract}

\maketitle

\setcounter{tocdepth}{1}
\tableofcontents

\section{Introduction}\label{S:intro}

Let $M$ be a closed connected smooth manifold with Riemannian metric $g$.
Suppose $E$ is a flat complex vector bundle over $M$. Let $h$ be a 
Hermitian metric on $E$. Recall the deRham differential
$d_E:\Omega^*(M;E)\to\Omega^{*+1}(M;E)$ on the space of $E$-valued
differential forms. Let $d^*_{E,g,h}:\Omega^{*+1}(M;E)\to\Omega^*(M;E)$
denote its formal adjoint with respect to the Hermitian scalar product on
$\Omega^*(M;E)$ induced by $g$ and $h$. Consider the Laplacian $\Delta_{E,g,h}
=d_Ed^*_{E,g,h}+d^*_{E,g,h}d_E:\Omega^*(M;E)\to\Omega^*(M;E)$. Recall the
(inverse square of the) Ray--Singer torsion \cite{RS71}
$$
\prod_q\bigl({\det}'(\Delta_{E,g,h,q})\bigr)^{(-1)^qq}\in\R^+.
$$
Here ${\det}'(\Delta_{E,g,h,q})$ denotes the zeta regularized product
of all non-zero eigen values of the Laplacian acting in degree $q$.
This is a positive real number which coincides, up to a computable
correction term, with the absolute value of the Reidemeister torsion, see \cite{BZ92}.

The aim of this paper is to introduce a complex valued Ray--Singer torsion
which, conjecturally, computes the Reidemeister torsion, including its phase.
This is accomplished by replacing the Hermitian fiber metric $h$ with a fiber wise
non-degenerate symmetric bilinear form $b$ on $E$. The bilinear form $b$ permits
to define a formal transposed
$d^\sharp_{E,g,b}$ of $d_E$, and an in general not selfadjoint Laplacian
$\Delta_{E,g,b}:=d_Ed^\sharp_{E,g,b}+d^\sharp_{E,g,b}d_E$.
The (inverse square of the) complex valued Ray--Singer torsion is then defined by
\begin{equation}\label{E:XX}
\prod_q\bigl({\det}'(\Delta_{E,g,b,q})\bigr)^{(-1)^qq}\in\C^\times:=\C\setminus\{0\}.
\end{equation}
The main result proved here, see Theorem~\ref{T:anom} below, is an anomaly formula for the
complex valued Ray--Singer torsion, i.e.\ we compute the variation of the quantity \eqref{E:XX}
through a variation of $g$ and $b$.
This ultimately permits to define a smooth invariant, the analytic torsion.

The paper is roughly organized as follows.
In Section~\ref{S:eul} we recall Euler and coEuler structures.
These are used to turn the Reidemeister torsion and the complex valued
Ray--Singer torsion into topological invariants referred to as
combinatorial and analytic torsion, respectively.
In Section~\ref{S:comb} we discuss some finite dimensional linear algebra 
and recall the combinatorial torsion which was also called Milnor--Turaev
torsion in \cite{BH03}.
Section~\ref{S:anator} contains the definition 
of the proposed complex valued analytic torsion.
In Section~\ref{S:BZ} we formulate a conjecture, see Conjecture~\ref{C:main}, 
relating the complex valued analytic torsion with the combinatorial torsion. 
We establish this conjecture in some non-trivial cases via analytic continuation 
from a result of Cheeger \cite{C77, C79}, M\"uller \cite{M78} and Bismut--Zhang \cite{BZ92}. 
Section~\ref{S:anom} contains the derivation of the 
anomaly formula. This proof is based on the computation of leading and subleading
terms in the asymptotic expansion of the heat kernel associated with a certain class
of Dirac operators.
This asymptotic expansion is formulated and proved in Section~\ref{S:ass_ex}, see
Theorem~\ref{T:ass_exp}. In Section~\ref{S:anomprop} we apply this result
to the Laplacians $\Delta_{E,g,b}$ and therewith complete the proof of the anomaly formula.

We restrict the presentation to the case of vanishing Euler--Poincar\'e characteristics
to avoid geometric regularization, see \cite{BH03} and \cite{BH06}.
With minor modifications everything can easily be extended to the general situation.
This is sketched in Section~\ref{S:gen}. The analytic core of the results, 
Theorem~\ref{T:ass_exp} and its corollaries Propositions~\ref{P:strbD} and \ref{P:strgD},
are formulated and proved without any restriction on the Euler--Poincar\'e 
characteristics.

Let us also mention the series of recent preprints 
\cite{BK05b, BK05c, BK05a, BK06a, BK06}. In these papers Braverman 
and Kappeler construct a \lq\lq refined analytic torsion\rq\rq\ based on the odd signature 
operator on odd dimensional manifolds. 
Their torsion is closely related to the analytic torsion proposed in this paper.
For a comparison result see Theorem~1.4 in \cite{BK06}.
Some of the results below which partially establish Conjecture~\ref{C:main},
have first appeared in \cite{BK06}, and were not contained in the first version of this paper.
The proofs we will provide have been inspired by \cite{BK06}
but do not rely on the results therein.

Recently, in October 2006, two preprints \cite{BH06b} and \cite{SZ06}
have been posted on the internet providing the proof of Conjecture~\ref{C:main}.
In \cite{BH06b} Witten--Helffer--Sj\"ostrand theory has been extended to
the non-self\-adjoint Laplacians discussed here, and used along the lines of 
\cite{BFK96}, to establish Conjecture~\ref{C:main} for odd dimensional manifolds, up to
sign. Comments were made how to derive the conjecture in full generality on these lines.
A few days earlier, by adapting the methods in \cite{BZ92} to the non-selfadjoint 
situation, Su and Zhang in \cite{SZ06} provided a proof of the conjecture.


The definition of the complex valued analytic torsion was sketched in \cite{BH05b}.

We thank the referees for useful remarks and for pointing out several sign mistakes.

\section{Preliminaries}\label{S:eul}

Throughout this section $M$ denotes a closed connected smooth manifold
of dimension $n$.
For simplicity we will also assume vanishing Euler--Poincar\'e characteristics,
$\chi(M)=0$. At the expense of a base point everything can easily be extended 
to the general situation, see \cite{B99}, \cite{BH03}, \cite{BH06} and Section~\ref{S:gen}.

\subsection*{Euler structures}

Let $M$ be a closed connected smooth manifold of dimension $n$ with $\chi(M)=0$.
The set of \emph{Euler structures with integral coefficients} $\Eul(M;\Z)$ 
is an \emph{affine
version} of $H_1(M;\Z)$. That is, the homology group $H_1(M;\Z)$ acts free
and transitively on $\Eul(M;\Z)$ but in general there is no distinguished
origin. Euler structures have been introduced by Turaev \cite{T90} in order 
to remove the ambiguities in the definition of the Reidemeister torsion.
Below we will briefly recall a possible definition. For more details we refer
to \cite{BH03} and \cite{BH06}.

Recall that a vector field $X$ is called non-degenerate if $X:M\to TM$ is transverse to 
the zero section.
Denote its set of zeros by $\mathcal X$. Recall that every $x\in\mathcal X$
has a \emph{Hopf index} $\IND_X(x)\in\{\pm1\}$.
Consider pairs $(X,c)$ where $X$ is a non-degenerate vector field and
$c\in C_1^\sing(M;\Z)$ is a singular $1$-chain satisfying
$$
\partial c=\ec(X):=\sum_{x\in\mathcal X}\IND_X(x)x.
$$
Every non-degenerate vector field admits such $c$ since 
we assumed $\chi(M)=0$.

We call two such pairs $(X_1,c_1)$ and $(X_2,c_2)$ equivalent if
$$
c_2-c_1=\cs(X_1,X_2)\in C_1^\sing(M;\Z)/\partial C_2^\sing(M;\Z).
$$
Here $\cs(X_1,X_2)\in C_1^\sing(M;\Z)/\partial C_2^\sing(M;\Z)$
denotes the \emph{Chern--Simons class} which is represented by the 
zero set of a generic homotopy connecting $X_1$ with $X_2$. 
It follows from $\cs(X_1,X_2)+\cs(X_2,X_3)=\cs(X_1,X_3)$ that this indeed
is an equivalence relation.

Define $\Eul(M;\Z)$ as the
set of equivalence classes $[X,c]$ of pairs considered above. The action of $[\sigma]\in H_1(M;\Z)$
on $[X,c]\in\Eul(M;\Z)$ is simply given by $[X,c]+[\sigma]:=[X,c+\sigma]$.
Since $\cs(X,X)=0$ this action is well defined and free. 
Because of $\partial\cs(X_1,X_2)=\ec(X_2)-\ec(X_1)$ it is transitive.

Replacing singular chains with integral coefficients by singular chains with
real or complex coefficients we obtain in exactly the same way 
\emph{Euler structures with real coefficients} $\Eul(M;\R)$ and
\emph{Euler structures with complex coefficients} $\Eul(M;\C)$. These are affine version 
of $H_1(M;\R)$ and $H_1(M;\C)$, respectively.
There are obvious maps $\Eul(M;\Z)\to\Eul(M;\R)\to\Eul(M;\C)$
which are affine over the homomorphisms $H_1(M;\Z)\to H_1(M;\R)\to H_1(M;\C)$.
We refer to the image of $\Eul(M;\Z)$ in $\Eul(M;\R)$ or $\Eul(M;\C)$ as the \emph{lattice of integral
Euler structures}.

Since we have $\ec(-X)=(-1)^n\ec(X)$ and
$\cs(-X_1,-X_2)=(-1)^n\cs(X_1,X_2)$, the assignment $\nu([X,c]):=[-X,(-1)^nc]$
defines \emph{affine involutions} on $\Eul(M;\Z)$, $\Eul(M;\R)$ and $\Eul(M;\C)$.
If $n$ is even, then the involutions on $\Eul(M;\R)$ and $\Eul(M;\C)$ are affine
over the identity and so we must have $\nu=\id$.
If $n$ is odd the involutions on $\Eul(M;\R)$ and $\Eul(M;\C)$ are affine over $-\id$
and thus must have a unique fixed 
point $\e_\can\in\Eul(M;\R)\subseteq\Eul(M;\C)$.
This \emph{canonic Euler structure} permits to naturally identify 
$\Eul(M;\R)$ resp.\ $\Eul(M;\C)$ 
with $H_1(M;\R)$ resp.\ $H_1(M;\C)$, provided $n$ is odd.
Note that in general none of these statements is true for the involution on
$\Eul(M;\Z)$. This is due to the fact that in general 
$H_1(M;\Z)$ contains non-trivial elements of order $2$, and elements which are not
divisible by $2$.

Finally, observe that the assignment $[X,c]\mapsto[X,\bar c]$
defines a \emph{conjugation} $\e\mapsto\bar\e$ on $\Eul(M;\C)$ which is
affine over the complex conjugation $H_1(M;\C)\to H_1(M;\C)$, $[\sigma]\mapsto[\bar\sigma]$.
Clearly, the set of fixed points of this conjugation coincides with 
$\Eul(M;\R)\subseteq\Eul(M;\C)$.

\begin{lemma}\label{L:spray}
Let $M$ be a closed connected smooth manifold with $\chi(M)=0$, let $\e\in\Eul(M;\Z)$ be an
Euler structure, and let $x_0\in M$ be a base point. Suppose $X$ is a
non-degenerate vector field on $M$ with zero set $\mathcal X\neq\emptyset$. 
Then there exists a collection of paths
$\sigma_x$, $\sigma_x(0)=x_0$, $\sigma_x(1)=x$,
$x\in\mathcal X$, so that 
$\e=\bigl[X,\sum_{x\in\mathcal X}\IND_X(x)\sigma_x\bigr]$.
\end{lemma}

\begin{proof}
For every zero $x\in\mathcal X$ choose a path $\tilde\sigma_x$ with
$\tilde\sigma_x(0)=x_0$ and $\tilde\sigma_x(1)=x$. Set $\tilde
c:=\sum_{x\in\mathcal X}\IND_X(x)\tilde\sigma_x$.
Since $\chi(M)=0$ we clearly have $\partial\tilde c=\ec(X)$.
So the pair $(X,\tilde c)$ represents an Euler structure
$\tilde\e:=[X,\tilde c]\in\Eul(M;\Z)$. Because $H_1(M;\Z)$ acts transitively
on $\Eul(M;\Z)$ we find $a\in H_1(M;\Z)$ with $\tilde\e+a=\e$.
Since the Hur\'ewicz homomorphism is onto we can represent $a$ by a 
closed path $\pi$ with $\pi(0)=\pi(1)=x_0$. Choose $y\in\mathcal X$.
Define $\sigma_y$ as the concatenation of $\tilde\sigma_y$
with $\pi^{\IND_X(y)}$, and set $\sigma_x:=\tilde\sigma_x$ for $x\neq y$.
Then the pair $(x,\sum_{x\in\mathcal X}\IND_X(x)\sigma_x)$ represents
$\tilde\e+a=\e$.
\end{proof}

\subsection*{CoEuler structures}

Let $M$ be a closed connected smooth manifold of dimension $n$ with $\chi(M)=0$.
The set of \emph{coEuler structures} $\Eul^*(M;\C)$ is an affine version of
$H^{n-1}(M;\mathcal O_M^\C)$. That is the cohomology group
$H^{n-1}(M;\mathcal O_M^\C)$ with values in the complexified orientation
bundle $\mathcal O_M^\C$ acts free and transitively on $\Eul^*(M;\C)$.
CoEuler structures are well suited to remove the metric dependence from
the Ray--Singer torsion. Below we will briefly recall their definition, and
discuss an affine version of Poincar\'e duality relating Euler with coEuler
structures. For more details and the general situation we refer to \cite{BH03} 
or \cite{BH06}.

Consider pairs $(g,\alpha)$, $g$ a Riemannian metric on $M$,
$\alpha\in\Omega^{n-1}(M;\mathcal O_M^\C)$, which satisfy
$$
d\alpha=\ec(g).
$$
Here $\ec(g)\in\Omega^n(M;\mathcal O_M^\C)$ denotes the \emph{Euler form}
associated with $g$. In view of the Gauss--Bonnet theorem every $g$ admits such 
$\alpha$ for we assumed $\chi(M)=0$.

Two pairs $(g_1,\alpha_1)$ and $(g_2,\alpha_2)$ as above are called equivalent
if
$$
\alpha_2-\alpha_1=\cs(g_1,g_2)\in
\Omega^{n-1}(M;\mathcal O_M^\C)/d\Omega^{n-2}(M;\mathcal O_M^\C).
$$
Here $\cs(g_1,g_2)\in\Omega^{n-1}(M;\mathcal O_M^\C)/d\Omega^{n-2}(M;\mathcal O_M^\C)$
denotes the \emph{Chern--Simons class} \cite{CS74} associated with $g_1$ and $g_2$.
Since $\cs(g_1,g_2)+\cs(g_2,g_3)=\cs(g_1,g_3)$ this is indeed an equivalence
relation.

Define the set of \emph{coEuler structures with complex coefficients}
$\Eul^*(M;\C)$ as the set of equivalence classes $[g,\alpha]$ of pairs considered above.
The action of $[\beta]\in H^{n-1}(M;\mathcal O_M^\C)$ on
$[g,\alpha]\in\Eul^*(M;\C)$ is
defined by $[g,\alpha]+[\beta]:=[g,\alpha-\beta]$.
Since $\cs(g,g)=0$ this action is well defined and free. Because of
$d\cs(g_1,g_2)=\ec(g_2)-\ec(g_1)$ it is transitive too.

Replacing forms with values in $\mathcal O_M^\C$ by forms with values in the
real orientation bundle $\mathcal O_M^\R$
we obtain in exactly the same way \emph{coEuler structures with real coefficients}
$\Eul^*(M;\R)$, an affine version of $H^{n-1}(M;\mathcal O_M^\R)$. 
There is an obvious map $\Eul^*(M;\R)\to\Eul^*(M;\C)$
which is affine over the homomorphism 
$H^{n-1}(M;\mathcal O_M^\R)\to H^{n-1}(M;\mathcal O_M^\C)$.

In view of $(-1)^n\ec(g)=\ec(g)$ and $(-1)^n\cs(g_1,g_2)=\cs(g_1,g_2)$ the
assignment $\nu([g,\alpha]):=[g,(-1)^n\alpha]$ defines \emph{affine involutions}
on $\Eul^*(M;\R)$ and $\Eul^*(M;\C)$. For even $n$ these involutions are 
affine over the identity and so we must have $\nu=\id$.
For odd $n$ they are affine over $-\id$ and thus must have a unique fixed point 
$\e^*_\can\in\Eul^*(M;\R)\subseteq\Eul^*(M;\C)$.
Since $\ec(g)=0$ in this case, we have $\e^*_\can=[g,0]$ where
$g$ is any Riemannian metric. This \emph{canonic coEuler structure} provides a natural
identification of $\Eul^*(M;\R)$ resp.\ $\Eul^*(M;\C)$ with $H^{n-1}(M;\mathcal O_M^\R)$
resp.\ $H^{n-1}(M;\mathcal O_M^\C)$, provided the dimension is odd.

Finally, observe that the assignment $[g,\alpha]\mapsto[g,\bar\alpha]$
defines a \emph{complex conjugation} $\e^*\mapsto\bar\e^*$ on $\Eul^*(M;\C)$ which is
affine over the complex conjugation $H^{n-1}(M;\mathcal O_M^\C)\to H^{n-1}(M;\mathcal O_M^\C)$,
$[\beta]\mapsto[\bar\beta]$.
Clearly, the set of fixed points of this conjugation coincides with the image of
$\Eul^*(M;\R)\subseteq\Eul^*(M;\C)$.

\subsection*{Poincar\'e duality for Euler structures}

Let $M$ be a closed connected smooth manifold of dimension $n$ with $\chi(M)=0$.
There is a canonic isomorphism
\begin{equation}\label{E:664}
P:\Eul(M;\C)\to\Eul^*(M;\C)
\end{equation}
which is affine over the Poincar\'e duality $H_1(M;\C)\to
H^{n-1}(M;\mathcal O_M^\C)$.
If $[X,c]\in\Eul(M;\C)$ and $[g,\alpha]\in\Eul^*(M;\C)$ then
$P([X,c])=[g,\alpha]$ iff we have
\begin{equation}\label{E:665}
\int_{M\setminus\mathcal X}\omega\wedge(X^*\Psi(g)-\alpha)=\int_c\omega
\end{equation}
for all closed one forms $\omega$ which vanish in a neighborhood 
of $\mathcal X$, the zero set of $X$.
Here $\Psi(g)\in\Omega^{n-1}(TM\setminus M;\pi^*\mathcal O_M^\C)$ denotes the 
\emph{Mathai--Quillen form} \cite{MQ86} associated with $g$, and $\pi:TM\to M$ 
denotes the projection.
With a little work one can show that \eqref{E:665} does indeed
define an assignment as in \eqref{E:664}. Once this is established
\eqref{E:664} is obviously affine over the Poincar\'e duality and hence an isomorphism.
It follows immediately from $(-X)^*\Psi(g)=(-1)^nX^*\Psi(g)$ that $P$ intertwines
the involution on $\Eul(M;\C)$ with the involution on $\Eul^*(M;\C)$.
Moreover, $P$ obviously intertwines the complex conjugations on 
$\Eul(M;\C)$ and $\Eul^*(M;\C)$. Particularly, \eqref{E:664} restricts
to an isomorphism
$$
P:\Eul(M;\R)\to\Eul^*(M;\R)
$$
affine over the Poincar\'e duality $H_1(M;\R)\to H^{n-1}(M;\mathcal O_M^\R)$.

\subsection*{Kamber--Tondeur form}

Suppose $E$ is a flat complex vector bundle over a smooth manifold $M$.
Let $\nabla^E$ denote the flat connection on $E$.
Suppose $b$ is a fiber wise non-degenerate symmetric bilinear form on $E$.
The \emph{Kamber--Tondeur form} is the one form 
\begin{equation}\label{E:KTKT}
\omega_{E,b}:=-\tfrac12\tr(b^{-1}\nabla^Eb)\in\Omega^1(M;\C).
\end{equation}
More precisely, for a vector field $Y$ on $M$ we have
$\omega_{E,b}(Y):=\tr(b^{-1}\nabla^E_Yb)$.
Here the derivative of $b$ with respect to the induced flat connection
on $(E\otimes E)'$ is considered as $\nabla_Y^Eb:E\to E'$.
Then $b^{-1}\nabla^E_Yb:E\to E$ and
$\omega_{E,b}(Y)$ is obtained by taking the fiber wise trace.

The bilinear form $b$ induces a non-degenerate bilinear form $\det b$ on 
$\det E:=\Lambda^{\rk(E)}E$. From $\det b^{-1}\nabla^{\det E}(\det b)=
\tr(b^{-1}\nabla^Eb)$ we obtain
\begin{equation}\label{E:612}
\omega_{\det E,\det b}=\omega_{E,b}.
\end{equation}
Particularly, $\omega_{E,b}$ depends on the flat line bundle $\det E$ and the
induced bilinear form $\det b$ only. Since $\nabla^E$ is flat,
$\omega_{E,b}$ is a closed $1$-form, cf.~\eqref{E:612}.

Suppose $b_1$ and $b_2$ are two fiber wise non-degenerate symmetric bilinear forms on $E$.
Set $A:=b_1^{-1}b_2\in\Aut(E)$, i.e.\ $b_2(v,w)=b_1(Av,w)$ for all $v,w$ in the same 
fiber of $E$. Then $\det b_2=\det b_1\det A$, hence
$$
\nabla^{\det E}(\det b_2)
=\nabla^{\det E}(\det b_1)\det A+(\det b_1)d\det A
$$ 
and therefore
\begin{equation}\label{E:wb12}
\omega_{E,b_2}=\omega_{E,b_1}+\det A^{-1}d\det A.
\end{equation}

If $\det b_1$ and $\det b_2$ are homotopic as fiber wise non-degenerate bilinear forms on 
$\det E$, then the function $\det A:M\to\C^\times:=\C\setminus\{0\}$ is homotopic to the constant function $1$. 
So we find a function $\log\det A:M\to\C$ with $d\log\det A=\det A^{-1}d\det A$, and
in view of \eqref{E:wb12} the cohomology classes of $\omega_{E,b_1}$ and $\omega_{E,b_2}$
coincide. We conclude that the cohomology class 
$[\omega_{E,b}]\in H^1(M;\C)$ depends on the flat line bundle $\det E$ and the homotopy 
class $[\det b]$ of the induced non-degenerate bilinear form $\det b$ on
$\det E$ only.

If $E_1$ and $E_2$ are two flat vector bundles with fiber wise non-degenerate symmetric
bilinear forms $b_1$ and $b_2$ then
\begin{equation}\label{E:KBE12}
\omega_{E_1\oplus E_2,b_1\oplus b_2}=\omega_{E_1,b_1}+\omega_{E_2,b_2}.
\end{equation}
If $E'$ denotes the dual of a flat vector bundle $E$, and if $b'$ denotes the bilinear form
on $E'$ induced from a fiber wise non-degenerate symmetric bilinear form $b$ on $E$ then
clearly
\begin{equation}\label{E:KTdual}
\omega_{E',b'}=-\omega_{E,b}.
\end{equation}
If $\bar E$ denotes the complex conjugate of a flat complex vector bundle $E$, and if
$\bar b$ denotes the complex conjugate bilinear form of a fiber wise non-degenerate
symmetric bilinear form $b$ on $E$, then obviously
\begin{equation}\label{E:KTconj}
\omega_{\bar E,\bar b}=\overline{\omega_{E,b}}.
\end{equation}
Finally, if $F$ is a real flat vector bundle and $h$ is a fiber wise non-degenerate
symmetric bilinear form on $F$ one defines in exactly the same way a real Kamber--Tondeur
form $\omega_{F,h}:=-\frac12\tr(h^{-1}\nabla^Fh)$ which is closed too. If $F^\C:=F\otimes\C$ denotes the 
complexification of $F$ and $h^\C$ denotes the complexification of $h$ then clearly
\begin{equation}\label{E:KTcx}
\omega_{F^\C,h^\C}=\omega_{F,h}
\end{equation}
in $\Omega^1(M;\R)\subseteq\Omega^1(M;\C)$.
Note that all such $h$ give rise to the same cohomology class $[\omega_{F,h}]\in H^1(M;\R)$, see
\eqref{E:612} and \eqref{E:wb12}.
To see this also note that the induced fiber wise non-degenerate bilinear form $\det h$ on $\det F$
has to be positive definite or negative definite, but $\omega_{\det F,-\det h}=\omega_{\det F, \det h}$.

\subsection*{Holonomy}

Suppose $E$ is a flat complex vector bundle over a connected smooth manifold
$M$. Let $x_0\in M$
be a base point. Parallel transport along closed loops provides an anti
homomorphism $\pi_1(M,x_0)\to\GL(E_{x_0})$, where $E_{x_0}$ denotes the
fiber of $E$ over $x_0$. Composing with the inversion in $\GL(E_{x_0})$ we
obtain the \emph{holonomy representation} of $E$ at $x_0$
$$
\hol^E_{x_0}:\pi_1(M,x_0)\to\GL(E_{x_0}).
$$

Applying this to the flat line bundle $\det E:=\Lambda^{\rk(E)}E$ we
obtain a homomorphism $\hol_{x_0}^{\det E}:\pi_1(M,x_0)\to\GL(\det
E_{x_0})=\C^\times$ which factors to a homomorphism
\begin{equation}\label{E:theta}
\theta_E:H_1(M;\Z)\to\C^\times.
\end{equation}

\begin{lemma}\label{L:A}
Suppose $b$ is a non-degenerate symmetric bilinear form on $E$. Then
$$
\theta_E(\sigma)=\pm e^{\langle[\omega_{E,b}],\sigma\rangle},\qquad
\sigma\in H_1(M;\Z).
$$
Here $\langle[\omega_{E,b}],\sigma\rangle\in\C$ denotes the natural pairing of
the cohomology class $[\omega_{E,b}]\in H^1(M;\C)$ and $\sigma\in H_1(M;\Z)$.
\end{lemma}

\begin{proof}
Let $\tau:[0,1]\to M$ be a smooth path with $\tau(0)=\tau(1)=x_0$.
Consider the flat vector bundle $(\det E)^{-2}:=(\det E\otimes\det E)'$.
Let $\beta:[0,1]\to(\det E)^{-2}$ be a section over $\tau$ which is parallel. 
Since $\det b$ defines a global nowhere vanishing section of 
$(\det E)^{-2}$ we find $\lambda:[0,1]\to\mathbb C$ so that
$\beta=\lambda\det b$. Clearly,
\begin{equation}\label{E:511}
\lambda(1)\hol^{(\det E)^{-2}}_{x_0}([\tau])=\lambda(0).
\end{equation}
Differentiating $\beta=\lambda\det b$ we obtain
$0=\lambda'\det b+\lambda\nabla^{(\det E)^{-2}}_{\tau'}(\det b)$. Using
\eqref{E:612} this yields $0=\lambda'-2\lambda\omega_{E,b}(\tau')$.
Integrating we get
$$
\lambda(1)=\lambda(0)\exp\Bigl(\int_0^12\omega_{E,b}(\tau'(t))dt\Bigr)
=\lambda(0)e^{2\langle[\omega_{E,b}],[\tau]\rangle}.
$$
Taking \eqref{E:511} into account we obtain
$\hol^{(\det E)^{-2}}_{x_0}([\tau])=e^{-2\langle[\omega_{E,b}],[\tau]\rangle}$,
and this gives
$\hol^{\det E}_{x_0}([\tau])=\pm e^{\langle[\omega_{E,b}],[\tau]\rangle}$.
\end{proof}

\section{Reidemeister torsion}\label{S:comb}

The combinatorial torsion is an invariant
associated to a closed connected smooth manifold $M$, an Euler 
structure with integral coefficients $\e$, and a flat complex vector bundle $E$ over $M$.
In the way we consider it here this invariant is a non-degenerate
bilinear form $\tau_{E,\e}^\comb$ on the complex line $\det H^*(M;E)$
--- the graded determinant line of the cohomology with values in
(the local system of coefficients provided by) $E$.
If $H^*(M;E)$ vanishes, then $\tau_{E,\e}^\comb$ becomes a non-vanishing complex number.
The aim of this section is to recall these definitions, and to provide
some linear algebra which will be used in the analytic approach to this
invariant in Section~\ref{S:anator}.

Throughout this section $M$ denotes a closed connected smooth manifold of dimension $n$.
For simplicity we will also assume vanishing Euler--Poincar\'e characteristics,
$\chi(M)=0$.
At the expense of a base point everything can easily be extended 
to the general situation, see \cite{B99}, \cite{BH03}, \cite{BH06} and Section~\ref{S:gen}.

\subsection*{Finite dimensional Hodge theory}

Suppose $C^*$ is a finite dimensional graded complex over $\C$ with differential
$d:C^*\to C^{*+1}$. Its cohomology is a finite dimensional graded vector space 
and will be denoted by $H(C^*)$. Recall that there is a canonic isomorphism of complex
lines 
\begin{equation}\label{E:fintor}
\det C^*=\det H(C^*).
\end{equation}
Let us explain the terms appearing in \eqref{E:fintor} in more details.
If $V$ is a finite dimensional vector space its \emph{determinant line}
is defined to be the top exterior product $\det V:=\Lambda^{\dim(V)}V$.
If $V^*$ is a finite dimensional graded vector space its \emph{graded determinant line}
is defined by $\det V^*:=\det V^\even\otimes(\det V^\odd)'$. Here
$V^\even:=\bigoplus_qV^{2q}$ and $V^\odd:=\bigoplus_qV^{2q+1}$ are considered as 
ungraded vector spaces and $V':=L(V;\C)$ denotes the dual space. For more details
on determinant lines consult for instance \cite{KM76}.
Let us only mention that every short exact sequence of graded vector spaces
$0\to U^*\to V^*\to W^*\to0$ provides a canonic isomorphism of determinant lines
$\det U^*\otimes\det W^*=\det V^*$. 
The complex $C^*$ gives rise to two short exact sequences 
\begin{equation}\label{E:001}
0\to B^*\to Z^*\to H(C^*)\to0\quad\text{and}\quad
0\to Z^*\to C^*\xrightarrow{d}B^{*+1}\to0
\end{equation} 
where $B^*$ and $Z^*$ denote the boundaries and cycles in $C^*$, respectively.
The isomorphism \eqref{E:fintor} is then obtained
from the isomorphisms of determinant lines induced by \eqref{E:001} together with the 
canonic isomorphism $\det B^*\otimes\det B^{*+1}=\det B^*\otimes(\det B^*)'=\C$.

Suppose our complex $C^*$ is equipped with a \emph{graded non-degenerate symmetric 
bilinear form} $b$. That is, we have a non-degenerate symmetric bilinear form 
on every homogeneous component $C^q$, and different homogeneous components are
$b$-orthogonal.
The bilinear form $b$ will induce a non-degenerate bilinear form on $\det C^*$.
Using \eqref{E:fintor} we obtain a non-degenerate bilinear form on $\det H(C^*)$
which is called the \emph{torsion associated with $C^*$ and $b$}. It will 
be denoted by $\tau_{C^*,b}$.

\begin{remark}\label{R:biline}
Note that a non-degenerate bilinear form on a complex line essentially is a non-vanishing
complex number. If $C^*$ happens to be acyclic, i.e.\ $H(C^*)=0$, then canonically
$\det H(C^*)=\C$ and $\tau_{C^*,b}\in\C^\times$ is a genuine non-vanishing
complex number --- the entry in the $1\times1$-matrix representing this
bilinear form. 
\end{remark}

\begin{example}\label{Ex:simtor}
Suppose $q\in\Z$, $n\in\N$ and $A\in\GL_n(\C)$. Let $C^*$ denote 
the acyclic complex $\C^n\xrightarrow{d=A}\C^n$
concentrated in degrees $q$ and $q+1$.
Let $b$ denote the standard non-degenerate symmetric bilinear form on $C^*$.
In this situation we have
$\tau_{C^*,b}=(\det A)^{(-1)^{q+1}2}=(\det AA^t)^{(-1)^{q+1}}$.
\end{example}

The bilinear form $b$ permits to define the \emph{transposed} $d^\sharp_b$ of $d$
$$
d_b^\sharp:C^{*+1}\to C^*,\quad
b(dv,w)=b(v,d_b^\sharp w),\quad v,w\in C^*.
$$
Define the \emph{Laplacian} $\Delta_b:=dd_b^\sharp+d_b^\sharp d:C^*\to C^*$.
Let us write $C^*_b(\lambda)$ for the generalized $\lambda$-eigen space of $\Delta_b$. 
Clearly,
\begin{equation}\label{E:findeco}
C^*=\bigoplus_\lambda C^*_b(\lambda).
\end{equation}

Since $\Delta_b$ is symmetric with respect to $b$, different generalized eigen spaces 
of $\Delta$ are $b$-orthogonal. It follows that the restriction of $b$ to
$C^*_b(\lambda)$ is non-degenerate.

Since $\Delta_b$ commutes with $d$ and $d_b^\sharp$ the latter two will preserve the
decomposition \eqref{E:findeco}. Hence every eigen space $C^*_b(\lambda)$ is a subcomplex
of $C^*$. The inclusion $C^*_b(0)\to C^*$ induces an isomorphism in cohomology.
Indeed, the Laplacian factors to an invertible map on $C^*/C^*_b(0)$ and thus
induces an isomorphism on $H(C^*/C^*_b(0))$. On the other hand, the equation
$\Delta_b=dd_b^\sharp+d_b^\sharp d$ tells that the Laplacian will induce the zero map on 
cohomology. Therefore $H(C^*/C^*_b(0))$ must vanish and $C^*_b(0)\to C^*$ is indeed 
a quasi isomorphism. Particularly, we obtain a canonic isomorphism of complex lines
\begin{equation}\label{E:finhodge}
\det H(C^*_b(0))=\det H(C^*).
\end{equation}

\begin{lemma}\label{L:fintor}
Suppose $C^*$ is a finite dimensional graded complex over $\C$ which is equipped
with a graded non-degenerate symmetric bilinear form $b$. Then via \eqref{E:finhodge}
we have
$$
\tau_{C^*,b}=\tau_{C^*_b(0),b|_{C^*_b(0)}}\cdot\prod_q({\det}'(\Delta_{b,q}))^{(-1)^qq}
$$
where ${\det}'(\Delta_{b,q})$ denotes the product over all non-vanishing
eigen values of the Laplacian acting in degree $q$, $\Delta_{b,q}:=\Delta_b|_{C^q}:C^q\to C^q$.
\end{lemma}

\begin{proof}
Suppose $(C_1^*,b_1)$ and $(C_2^*,b_2)$ are finite dimensional complexes
equipped with graded non-degenerate symmetric bilinear forms. Clearly,
$H(C_1^*\oplus C_2^*)=H(C_1^*)\oplus H(C_2^*)$ and we obtain a canonic
isomorphism of determinant lines
$$
\det H(C_1^*\oplus C_2^*)=\det H(C_1^*)\otimes\det H(C_2^*).
$$
It is not hard to see that via this identification we have
\begin{equation}\label{E:sums}
\tau_{C_1^*\oplus C_2^*,b_1\oplus b_2}=
\tau_{C_1^*,b_1}\otimes\tau_{C_2,b_2}.
\end{equation}
In view of the $b$-orthogonal decomposition \eqref{E:findeco} we may therefore
w.l.o.g.\ assume $\ker\Delta_b=0$. Particularly, $C^*$ is acyclic.

Then $\img d\cap\ker d_b^\sharp\subseteq\ker d\cap\ker d_b^\sharp\subseteq\ker\Delta_b=0$. 
Since $\img d$ and $\ker d_b^\sharp$
are of complementary dimension we conclude $\img d\oplus\ker d_b^\sharp=C^*$.
The acyclicity of $C^*$ implies $\ker d_b^\sharp=\img d_b^\sharp$ and hence
$\img d\oplus\img d_b^\sharp=C^*$.
This decomposition is $b$-orthogonal and invariant under $\Delta_b$.
We obtain
$$
{\det}'(\Delta_{b,q})=\det(\Delta_{b,q})=\det(\Delta_b|_{C^q\cap\img d})
\cdot\det(\Delta_b|_{C^q\cap\img d_b^\sharp}).
$$
Since $d:C^q\cap\img d_b^\sharp\to C^{q+1}\cap\img d$ is an isomorphism commuting with $\Delta$
$$
\det(\Delta|_{C^q\cap\img d_b^\sharp})=
\det(\Delta|_{C^{q+1}\cap\img d}).
$$
A telescoping argument then shows
\begin{equation}\label{E:002}
\prod_q({\det}'(\Delta_{b,q}))^{(-1)^qq}
=\prod_q\det(\Delta_b|_{C^q\cap\img d})^{(-1)^q}.
\end{equation}
On the other hand, the $b$-orthogonal decomposition of complexes
$$
C^*=\bigoplus_q\Bigl(C^q\cap\img d_b^\sharp\xrightarrow{d}C^{q+1}\cap\img d\Bigr)
$$
together with \eqref{E:sums} and the computation in Example~\ref{Ex:simtor}
imply
$$
\tau_{C^*,b}=\prod_q\det\bigl(dd_b^\sharp|_{C^{q+1}\cap\img d}\bigr)^{(-1)^{q+1}}
$$
which clearly coincides with \eqref{E:002} since 
$\Delta_b|_{\img d}=dd_b^\sharp|_{\img d}$.
\end{proof}

\begin{example}\label{Ex:1}
Suppose $0\neq v\in\C^2$ satisfies $v^tv=0$. Moreover, suppose
$0\neq z\in\C$ and set $w:=zv^t$. Let $C^*$ denote the acyclic complex
$\C\xrightarrow{v}\C^2\xrightarrow{w}\C$ concentrated in degrees $0$, $1$ and $2$.
Equip this complex with the standard symmetric bilinear form $b$.
Then $\Delta_{b,0}=v^tv=0$, $\Delta_{b,2}=ww^t=0$, $\Delta_{b,1}=(1+z^2)vv^t$,
$(\Delta_{b,1})^2=0$. Thus all of this complex is contained in the generalized
$0$-eigen space of $\Delta_b$. The torsion of the complex computes to 
$\tau_{C^*,b}=-z^2$. 
Observe that the kernel of $\Delta_b$ does not compute the cohomology; 
that the bilinear form becomes degenerate when restricted 
to the kernel of $\Delta_b$; and that the torsion cannot be computed from the
spectrum of $\Delta_b$. 
\end{example}

\subsection*{Morse complex}

Let $E$ be a flat complex vector bundle over a closed connected
smooth manifold $M$ of dimension $n$.
Suppose $X=-\grad_g(f)$ is a \emph{Morse--Smale vector field} on $M$, see \cite{S93}.
Let $\mathcal X$ denote the zero set of $X$. Elements in $\mathcal X$ are
called \emph{critical points} of $f$.
Every $x\in\mathcal X$ has a \emph{Morse index} $\ind(x)\in\N$ which
coincides with the dimension of the unstable manifold of $x$ with respect to $X$. 
We will write $\mathcal X_q:=\{x\in\mathcal X\mid \ind(x)=q\}$ for the
set of critical points of index $q$.

Recall that the Morse--Smale vector field provides a \emph{Morse complex} $C^*(X;E)$ 
with underlying finite dimensional graded vector space 
$$
C^q(X;E)=\bigoplus_{x\in\mathcal X_q}E_x\otimes_{\{\pm1\}}\mathcal O_x. 
$$
Here $E_x$ denotes the fiber of $E$ over $x$, and $\mathcal O_x$ denotes the set
of orientations of the unstable manifold of $x$. The Smale condition
tells that stable and unstable manifolds intersect transversally. It follows
that for two critical points of index difference one there is only a finite
number of unparametrized trajectories connecting them. The differential
in $C^*(X;E)$ is defined with the help of these isolated trajectories and parallel
transport in $E$ along them.

Integration over unstable manifolds provides a homomorphism of complexes
\begin{equation}\label{E:int}
\Int:\Omega^*(M;E)\to C^*(X;E)
\end{equation}
where $\Omega^*(M;E)$ denotes the deRham complex with values in $E$.
It is a folklore fact that \eqref{E:int} induces an isomorphism on
cohomology, see \cite{S93}.
Particularly, we obtain a canonic isomorphism of complex lines
\begin{equation}\label{E:rm}
\det H^*(M;E)=\det H(C^*(X;E)).
\end{equation}

Suppose $\chi(M)=0$ and let $\e\in\Eul(M;\Z)$ be an Euler structure.
Choose a base point $x_0\in M$.
For every critical point $x\in\mathcal X$ choose a path $\sigma_x$ with $\sigma(0)=x_0$
and $\sigma_x(1)=x$ so that $\e=\bigl[-X,\sum_{x\in\mathcal
X}(-1)^{\ind(x)}\sigma_x\bigr]$.
This is possible in view of Lemma~\ref{L:spray}. Also
note that $\IND_{-X}(x)=(-1)^{\ind(x)}$.
Choose a non-degenerate symmetric bilinear form $b_{x_0}$ on the fiber $E_{x_0}$ over
$x_0$. For $x\in\mathcal X$ define a bilinear form
$b_x$ on $E_x$ by parallel transport of $b_{x_0}$ along $\sigma_x$. The
collection of bilinear forms $\{b_x\}_{x\in\mathcal X}$ defines a
non-degenerate symmetric bilinear form on the Morse complex $C^*(X;E)$.
It is elementary to check that the induced bilinear form on $\det C^*(X;E)$
does not depend on the choice of $\{\sigma_x\}_{x\in\mathcal X}$, and
because $\chi(M)=0$ it does not depend on $x_0$ or $b_{x_0}$ either. 
Hence the corresponding torsion is a non-degenerate bilinear form on
$\det H(C^*(X;E))$ depending on $E$, $\e$ and $X$ only.
Using \eqref{E:rm} we obtain a non-degenerate bilinear form on $\det H^*(M;E)$
which we will denote by $\tau_{E,\e,X}^\comb$.
For the following non-trivial statement we refer to \cite{M66}, \cite{T90}
or \cite{L92}.

\begin{theorem}[Milnor, Turaev]\label{T:MT}
The bilinear form $\tau_{E,\e,X}^\comb$ does not depend on $X$.
\end{theorem}

In view of Theorem~\ref{T:MT} we will denote $\tau^\comb_{E,\e,X}$ by
$\tau^\comb_{E,\e}$ from now on.

\begin{definition}[Combinatorial torsion]\label{D:comb}
The non-degenerate bilinear form $\tau^\comb_{E,\e}$ on $\det H^*(M;E)$
is called the \emph{combinatorial torsion} associated with the flat complex
vector bundle $E$ and the Euler structure $\e\in\Eul(M;\Z)$.
\end{definition}

\begin{remark}\label{R:combeul}
The combinatorial torsion's dependence on the Euler structure is very
simple. For $\e\in\Eul(M;\Z)$ and $\sigma\in H_1(M;\Z)$ we obviously 
have, see \eqref{E:theta}
$$
\tau^\comb_{E,\e+\sigma}=\tau^\comb_{E,\e}\cdot\theta_E(\sigma)^2.
$$
\end{remark}

The dependence on $E$, i.e.\ the dependence on the flat connection, is 
subtle and interesting. Let us only mention the following

\begin{example}[Torsion of mapping tori]\label{Ex:MT}
Consider a mapping torus
$$
M=N\times[0,1]/_{(x,1)\sim(\varphi(x),0)}
$$ 
where $\varphi:N\to N$ is a
diffeomorphism. Let $\pi:M\to S^1=[0,1]/_{0\sim1}$ denote the canonic
projection. The set of vector fields which project to the vector 
field $-\frac\partial{\partial\theta}$ on $S^1$ is contractible and thus
defines an Euler structure $\e\in\Eul(M;\Z)$ represented by $[X,0]$
where $X$ is any of these vector fields.
Let $\tilde E^z$ denote the flat line bundle over $S^1$ with
holonomy $z\in\C^\times$, i.e.\ $\theta_{\tilde E^z}:H_1(S^1;\Z)=\Z\to\C^\times$,
$\theta_{\tilde E^z}(k)=z^k$.
Consider the flat line bundle
$E^z:=\pi^*\tilde E^z$ over $M$. It follows from the Wang sequence of
the fibration $\pi:M\to S^1$ that for generic $z$ we will have
$H^*(M;E^z)=0$. In this case
$$
\tau^\comb_{E^z,\e}=(\zeta_\varphi(z))^2
$$
where 
\begin{align*}
\zeta_\varphi(z)&
=\exp\left(\sum_{k\geq1}\str\Bigl(H^*(N;\Q)\xrightarrow{(\varphi^k)^*}H^*(N;\Q)\Bigr)\frac{z^k}k\right)
\\
&=\sdet\Bigl(H^*(N;\C)\xrightarrow{1-z\varphi^*} H^*(N;\C)\Bigr)^{-1}
\end{align*}
denotes the Lefschetz zeta function of $\varphi$.
Here we wrote $\str$ and $\sdet$ for the super trace and the super determinant, respectively.
For more details and proofs we refer to \cite{F87} and \cite{BH03}.
\end{example}

\begin{remark}
Often the combinatorial 
torsion is considered as an element in (rather than
a bilinear form on) $\det H^*(M;E)$. This element is one of the two unit
vectors of $\tau^\comb_{E,\e}$. It is a non-trivial task (and requires the
choice of a homology orientation) to fix the sign, i.e.\ to describe
which of the two unite vectors actually is the torsion \cite{FT00}.
Considering bilinear forms this sign issue disappears.
\end{remark}

\subsection*{Basic properties of the combinatorial torsion}

If $E_1$ and $E_2$ are two flat vector bundles over $M$ 
then we have a canonic isomorphism
$H^*(M;E_1\oplus E_2)=H^*(M;E_1)\oplus H^*(M;E_2)$ which induces a canonic
isomorphism of complex lines
$\det H^*(M;E_1\oplus E_2)=\det H^*(M;E_1)\otimes\det H^*(M;E_2)$.
Via this identification we have
\begin{equation}\label{E:combsum}
\tau^\comb_{E_1\oplus E_2,\e}=\tau^\comb_{E_1,\e}\otimes\tau^\comb_{E_2,\e}.
\end{equation}
This follows from $C^*(X;E_1\oplus E_2)=C^*(X;E_1)\oplus C^*(X;E_2)$ and \eqref{E:sums}.

If $E'$ denotes the dual of a flat vector bundle $E$ then Poincar\'e duality
induces an isomorphism $H^*(M;E'\otimes\mathcal O_M)=H^{n-*}(M;E)'$ which induces a canonic
isomorphism $\det H^*(M;E'\otimes\mathcal O_M)=(\det H^*(M;E))^{(-1)^{n+1}}$. Via this identification we
have
\begin{equation}\label{E:dualcomb}
\tau^\comb_{E'\otimes\mathcal O_M,\nu(\e)}=(\tau^\comb_{E,\e})^{(-1)^{n+1}}
\end{equation}
where $\nu$ denotes the involution on $\Eul(M;\Z)$ discussed in Section~\ref{S:eul}. 
To see that use a
Morse--Smale vector field $X$ to compute $\tau^\comb_{E,\e}$ and use the
Morse--Smale vector field $-X$ to compute $\tau^\comb_{E'\otimes\mathcal O_M,\nu(\e)}$. Then there
is an obvious isomorphism of complexes $C^*(-X;E'\otimes\mathcal O_M)=C^{n-*}(X;E)'$ which
induces Poincar\'e duality on cohomology.

If $V$ is a complex vector space let $\bar V$ denote the complex conjugate
vector space. If $b$ is a bilinear form on $V$ let $\bar b$ denote the complex 
conjugate bilinear form on $\bar V$, that is $\bar b(v,w)=\overline{b(v,w)}$.
Let $\bar E$ denote the complex conjugate of a flat vector bundle $E$.
Then we have a canonic isomorphism $H^*(M;\bar E)=\overline{H^*(M;E)}$ which induces
a canonic isomorphism of complex lines $\det H^*(M;\bar E)=\overline{\det H^*(M;E)}$.
Via this identification we have
\begin{equation}\label{E:combconj}
\tau^\comb_{\bar E,\e}=\overline{\tau^\comb_{E,\e}}.
\end{equation}
This follows from $C^*(X;\bar E)=\overline{C^*(X;E)}$.

If $V$ is a real vector space we let $V^\C:=V\otimes\C$ denote its complexification.
If $h$ is a real bilinear form on $V$ we let $h^\C$ denote its complexification, 
more explicitly $h^\C(v_1\otimes z_1,v_2\otimes z_2)=h(v_1,v_2)z_1z_2$.
If $F$ is real flat vector bundle its torsion, defined analogously to the
complex case, is a real non-degenerate bilinear form on $\det H^*(M;F)$. 
Let $F^\C=F\otimes\C$ denote the complexification of the flat vector bundle $F$.
We have
a canonic isomorphism $H^*(M;F^\C)=H^*(M;F)^\C$ which induces a canonic isomorphism
of complex lines $\det H^*(M;F^\C)=(\det H^*(M;F))^\C$. Via this identification we have
\begin{equation}\label{E:cxcomb}
\tau^\comb_{F^\C,\e}=(\tau^\comb_{F,\e})^\C.
\end{equation}
This follows from $C^*(X;F^\C)=C^*(X;F)^\C$. Note that $\tau^\comb_{F,\e}$ 
is positive definite.

\section{Ray--Singer torsion}\label{S:anator}

The analytic torsion defined below is an invariant 
associated to a closed connected smooth manifold $M$, a complex
flat vector bundle $E$ over $M$, a coEuler structure $\e^*$
and a homotopy class $[b]$ of fiber wise non-degenerate symmetric bilinear forms on $E$.
In the way considered below, this
invariant is a non-degenerate symmetric bilinear form $\tau^\an_{E,\e^*,[b]}$
on the complex line $\det H^*(M;E)$. If $H^*(M;E)$ vanishes, then
$\tau^\an_{E,\e^*,[b]}$ becomes a non-vanishing complex number.

Throughout this section $M$ denotes a closed connected smooth manifold of dimension $n$.
For simplicity we will also assume vanishing Euler--Poincar\'e characteristics,
$\chi(M)=0$.
At the expense of a base point everything can easily be extended 
to the general situation, see \cite{B99}, \cite{BH03}, \cite{BH06} and Section~\ref{S:gen}.

\subsection*{Laplacians and spectral theory}

Suppose $M$ is a closed connected smooth manifold of dimension $n$. Let $E$ be a flat vector bundle over
$M$. We will denote the flat connection of $E$ by $\nabla^E$. Suppose
there exists a \emph{fiber wise non-degenerate symmetric bilinear} form $b$ on $E$.
Moreover, let $g$ be a Riemannian metric on $M$. This permits to define a
symmetric bilinear form $\beta_{g,b}$ on the space of $E$-valued differential
forms $\Omega^*(M;E)$,
$$
\beta_{g,b}(v,w):=\int_Mv\wedge(\star_g\otimes b)w,
\qquad v,w\in\Omega^*(M;E).
$$
Here $\star_g\otimes b:
\Omega^*(M;E)\to\Omega^{n-*}(M;E'\otimes\mathcal O_M)$ denotes the isomorphism
induced by the Hodge star operator\footnote{The normalization of the Hodge star operator we are using is
$\alpha_1\wedge\star_g\alpha_2=\langle\alpha_1,\alpha_2\rangle_g\Omega_g$, where $\alpha_1,\alpha_2\in\Omega(M;\R)$,
$\Omega_g\in\Omega^n(M;\mathcal O_M)$ denotes the volume density associated with $g$, and $\langle\alpha_1,\alpha_2\rangle_g$
denotes the inner product on $\Lambda^*T^*M$ induced by $g$, see \cite[Section~2.1]{J05}. Although we will frequently 
refer to \cite{BGV92} in the subsequent sections, the convention for the Hodge star operator we are using
differs from the one in \cite{BGV92}.}
$\star_g:\Omega^*(M;\R)\to\Omega^{n-*}(M;\mathcal O_M)$
and the isomorphism of vector bundles $b:E\to E'$.
The wedge product is computed with respect to the canonic pairing of
$E\otimes E'\to\C$.

Let $d_E:\Omega^*(M;E)\to\Omega^{*+1}(M;E)$ denote the deRham differential.
Let 
$$
d^\sharp_{E,g,b}:\Omega^{*+1}(M;E)\to\Omega^*(M;E)
$$
denote its formal transposed with respect to $\beta_{g,b}$. A straight forward
computation shows that $d^\sharp_{E,g,b}:\Omega^q(M;E)\to\Omega^{q-1}(M;E)$
is given by
\begin{equation}\label{E:YY}
d^\sharp_{E,g,b}
=(-1)^q(\star_g\otimes b)^{-1}\circ d_{E'\otimes\mathcal O_M}\circ(\star_g\otimes b).
\end{equation}
Define the \emph{Laplacian} by
\begin{equation}\label{E:Laplace}
\Delta_{E,g,b}:=
d_E\circ d^\sharp_{E,g,b}+d^\sharp_{E,g,b}\circ d_E.
\end{equation}
These are generalized Laplacians in the sense that their principal symbol coincides with
the symbol of the Laplace--Beltrami operator.

In the next proposition we collect some well known facts concerning the
spectral theory of $\Delta_{E,g,b}$. For details we refer to \cite{S01},
particularly Theorems~8.4 and 9.3 therein.

\begin{proposition}\label{P:spec}
For the Laplacian $\Delta_{E,g,b}$ constructed above the following hold:
\begin{enumerate}
\item\label{P:spec:i}
The spectrum of $\Delta_{E,g,b}$ is discrete.
For every $\theta>0$ all but finitely many points of the spectrum
are contained in the angle $\{z\in\C\mid-\theta<\arg(z)<\theta\}$.
\item\label{P:spec:ii}
If $\lambda$ is in the spectrum of $\Delta_{E,g,b}$ then the image of
the associated spectral projection 
is finite dimensional and contains smooth forms only. 
We will refer to this image as the (generalized) $\lambda$-eigen space
of $\Delta_{E,g,b}$ and denote it by $\Omega^*_{g,b}(M;E)(\lambda)$.
There exists $N_\lambda\in\N$ such that
$$
(\Delta_{E,g,b}-\lambda)^{N_\lambda}|_{\Omega^*_{g,b}(M;E)(\lambda)}=0.
$$
We have a $\Delta_{E,g,b}$-invariant $\beta_{g,b}$-orthogonal decomposition
\begin{equation}\label{E:deco}
\Omega^*_{g,b}(M;E)=\Omega^*_{g,b}(M;E)(\lambda)
\oplus\Omega^*_{g,b}(M;E)(\lambda)^{\perp_{\beta_{g,b}}}.
\end{equation}
The restriction of $\Delta_{E,g,b}-\lambda$ to
$\Omega^*_{g,b}(M;E)(\lambda)^{\perp_{\beta_{g,b}}}$ is invertible.
\item\label{P:spec:iii}
The decomposition \eqref{E:deco} 
is invariant under $d_E$ and $d_{E,g,b}^\sharp$.
\item\label{P:spec:iv}
For $\lambda\neq\mu$ the eigen spaces $\Omega^*_{g,b}(M;E)(\lambda)$
and $\Omega^*_{g,b}(M;E)(\mu)$ are orthogonal with respect to $\beta_{g,b}$.
\end{enumerate}
\end{proposition}

In view of Proposition~\ref{P:spec} the generalized $0$-eigen space
$\Omega^*_{g,b}(M;E)(0)$ is a finite dimensional subcomplex of
$\Omega^*(M;E)$.
The inclusion 
\begin{equation}\label{E:incl0}
\Omega^*_{g,b}(M;E)(0)\to\Omega^*(M;E)
\end{equation} 
induces an isomorphism in
cohomology. Indeed, in view of Proposition~\ref{P:spec}\itemref{P:spec:ii}
the Laplacian $\Delta_{E,g,b}$ induces an isomorphism on
$\Omega^*_{g,b}(M;E)/\Omega^*_{g,b}(M;E)(0)$ and thus an isomorphism on
$H(\Omega^*_{g,b}(M;E)/\Omega^*_{g,b}(M;E)(0))$. On the other hand
\eqref{E:Laplace} tells that $\Delta_{E,g,b}$ induces $0$ on cohomology,
hence $H(\Omega^*_{g,b}(M;E)/\Omega^*_{g,b}(M;E)(0))$ must vanish and \eqref{E:incl0}
is indeed a quasi isomorphism. We obtain a canonic isomorphism of complex lines
\begin{equation}\label{E:detdr}
\det H(\Omega^*_{g,b}(M;E)(0))=\det H^*(M;E).
\end{equation}

In view of Proposition~\ref{P:spec}\itemref{P:spec:ii} the bilinear form
$\beta_{g,b}$ restricts to a non-de\-ge\-ner\-ate bilinear form on
$\Omega^*_{g,b}(M;E)(0)$. Using the linear algebra discussed in
Section~\ref{S:comb} we obtain a non-degenerate bilinear form
on $\det H(\Omega^*_{g,b}(M;E)(0))$. Via \eqref{E:detdr} this
gives rise to a non-degenerate bilinear form on $\det H^*(M;E)$
which will be denoted by $\tau^\an_{E,g,b}(0)$.

Let $\Delta_{E,g,b,q}$ denote the Laplacian acting in degree $q$.
Define the zeta regularized product of its non-vanishing eigen values, as
$$
{\det}'(\Delta_{E,g,b,q})
:=\exp\left(-\frac\partial{\partial
s}\Big|_{s=0}\tr\Bigl(\bigl(\Delta_{E,g,b,q}|{\Omega^q_{g,b}(M;E)(0)^{\perp_{\beta_{g,b}}}}\bigr)^{-s}\Bigr)\right).
$$
Here the complex powers are defined with respect to any non-zero Agmon angle which
avoids the spectrum of 
$\Delta_{E,g,b,q}|{\Omega^q_{g,b}(M;E)(0)^{\perp_{\beta_{g,b}}}}$, see
Proposition~\ref{P:spec}\itemref{P:spec:i}. Recall that
for $\Re(s)>n/2$ the operator
$(\Delta_{E,g,b,q}|{\Omega^q_{g,b}(M;E)(0)^{\perp_{\beta_{g,b}}}})^{-s}$
is trace class. As a function in $s$ this trace extends to a meromorphic function on the
complex plane which is holomorphic at $0$, see \cite{Se67} or \cite[Theorem~13.1]{S01}. 
It is clear from
Proposition~\ref{P:spec}\itemref{P:spec:i} that ${\det}'(\Delta_{E,g,b,q})$
does not depend on the Agmon angle used to define the complex powers.

Assume $\chi(M)=0$ and suppose $\alpha\in\Omega^{n-1}(M;\mathcal O_M^\C)$ 
such that $d\alpha=\ec(g)$.
Consider the non-degenerate bilinear form on $\det H^*(M;E)$ defined by, cf.\ \eqref{E:KTKT},
$$
\tau^\an_{E,g,b,\alpha}
:=\tau^\an_{E,g,b}(0)\cdot\prod_q\bigl({\det}'(\Delta_{E,g,b,q})\bigr)^{(-1)^qq}\cdot
\exp\Bigl(-2\int_M\omega_{E,b}\wedge\alpha\Bigr).
$$

In Section~\ref{S:anom} we will provide a proof of the following result
which can be interpreted as an anomaly formula for the complex valued
Ray--Singer torsion \eqref{E:XX}.

\begin{theorem}[Anomaly formula]\label{T:anom}
Let $M$ be a closed connected smooth manifold with vanishing Euler--Poincar\'e characteristics.
Let $E$ be a flat complex vector bundle over $M$.
Suppose $g_u$ is a smooth one-parameter family of Riemannian metrics on $M$,
and $\alpha_u\in\Omega^{n-1}(M;\mathcal O_M^\C)$ is a smooth one-parameter
family so that $[g_u,\alpha_u]$ represent the same coEuler structure in $\Eul^*(M;\C)$.
Moreover, suppose $b_u$ is a smooth one-parameter family of fiber wise 
non-degenerate symmetric bilinear forms on $E$.
Then, as bilinear forms on $\det H^*(M;E)$, we have
$\frac\partial{\partial u}\tau^\an_{E,g_u,b_u,\alpha_u}=0$.
\end{theorem}

In view of Theorem~\ref{T:anom} the bilinear form $\tau^\an_{E,g,b,\alpha}$
does only depend on the flat vector bundle $E$, the coEuler structure
$\e^*\in\Eul^*(M;\C)$ represented by $(g,\alpha)$, and the homotopy class $[b]$ of $b$.
We will denote it by $\tau^\an_{E,\e^*,[b]}$ from now on.

\begin{definition}[Analytic torsion]
The non-degenerate bilinear form $\tau^\an_{E,\e^*,[b]}$ on $\det H^*(M;E)$
is called the \emph{analytic torsion} associated to the flat complex vector
bundle $E$, the coEuler structure $\e^*\in\Eul^*(M;\C)$ and the homotopy class $[b]$ of
fiber wise non-degenerate symmetric bilinear forms on $E$.
\end{definition}

\begin{remark}\label{R:aneul}
The analytic torsion's dependence on the coEuler structure is very simple. For $\e^*\in\Eul^*(M;\C)$ and 
$\beta\in H^{n-1}(M;\mathcal O_M^\C)$ we obviously have:
$$
\tau^\an_{E,\e^*+\beta,[b]}
=\tau^\an_{E,\e^*,[b]}\cdot\Bigl(e^{\langle[\omega_{E,b}]\cup\beta,[M]\rangle}\Bigr)^2
$$
Here $\langle[\omega_{E,b}]\cup\beta,[M]\rangle\in\C$ denotes the evaluation
of $[\omega_{E,b}]\cup\beta\in H^n(M;\mathcal O_M^\C)$ on the fundamental class
$[M]\in H_n(M;\mathcal O_M)$.
\end{remark}

\begin{remark}\label{R:cancoeul}
Recall from Section~\ref{S:eul} that for odd $n$ there is a canonic 
coEuler structure $\e^*_\can\in\Eul^*(M;\C)$ given by $\e^*_\can=[g,0]$.
The corresponding analytic torsion is:
$$
\tau^\an_{E,\e^*_\can,[b]}
=\tau^\an_{E,g,b}(0)\cdot\prod_q\bigl({\det}'(\Delta_{E,g,b,q})\bigr)^{(-1)^qq}
$$
Note however that in general this does depend on the homotopy class $[b]$, see for instance
the computation for the circle in Section~\ref{S:BZ} below.
This is related to the fact that $\e^*_\can$ in general is not integral, 
cf.\ Remark~\ref{R:bindep} below.
\end{remark}

\subsection*{Basic properties of the analytic torsion}

Suppose $E_1$ and $E_2$ are two flat vector bundles with fiber wise non-degenerate
symmetric bilinear forms $b_1$ and $b_2$. Via the canonic isomorphism of
complex lines $\det H^*(M;E_1\oplus E_2)=\det H^*(M;E_1)\otimes\det
H^*(M;E_2)$ we have:
\begin{equation}\label{E:578}
\tau^\an_{E_1\oplus E_2,\e^*,[b_1\oplus b_2]}
=\tau^\an_{E_1,\e^*,[b_1]}\otimes\tau^\an_{E_2,\e^*,[b_2]}
\end{equation}
For this note that via the identification $\Omega^*(M;E_1\oplus E_2)=\Omega^*(M;E_1)\oplus\Omega^*(M;E_2)$
we have $\Delta_{E_1\oplus E_2,g,b_1\oplus b_2}=\Delta_{E_1,g,b_1}\oplus\Delta_{E_2,g,b_2}$,
hence ${\det}'(\Delta_{E_1\oplus E_2,g,b_1\oplus b_2,q})=
{\det}'(\Delta_{E_1,g,b_1,q}){\det}'(\Delta_{E_2,g,b_2,q})$.
Moreover, recall \eqref{E:KBE12} for the correction terms.

Suppose $E'$ is the dual of a flat vector bundle $E$. Let $b'$ denote the
bilinear form on $E'$ dual to the non-degenerate symmetric bilinear form $b$
on $E$. The bilinear form $b'$ induces a fiber wise non-degenerate symmetric bilinear form
on the flat vector bundle $E'\otimes\mathcal O_M$ which will be denoted by
$b'$ too. Via the canonic isomorphism of complex lines
$\det H^*(M;E'\otimes\mathcal O_M)=(\det H^*(M;E))^{(-1)^{n+1}}$ induced by Poincar\'e duality we have
\begin{equation}\label{E:dualan}
\tau^\an_{E'\otimes\mathcal O_M,\nu(\e^*),[b']}=\bigl(\tau^\an_{E,\e^*,[b]}\bigr)^{(-1)^{n+1}}
\end{equation}
where $\nu$ denotes the involution introduced in Section~\ref{S:eul}.
This follows from the fact that $\star_g\otimes
b:\Omega^q(M;E)\to\Omega^{n-q}(M;E'\otimes\mathcal O_M)$ intertwines
the Laplacians $\Delta_{E,g,b,q}$ and $\Delta_{E'\otimes\mathcal
O_M,g,b',n-q}$, see
\eqref{E:YY}. Therefore $\Delta_{E,g,b,q}$ and $\Delta_{E'\otimes\mathcal
O_M,g,b',n-q}$ are isospectral and thus
${\det}'(\Delta_{E,g,b,q})={\det}'(\Delta_{E'\otimes\mathcal O_M,g,b',n-q})$.
Here one also has to use $\prod_q\bigl({\det}'(\Delta_{E,g,b,q})\bigr)^{(-1)^q}=1$,
and \eqref{E:KTdual}.

Let $\bar E$ denote the complex conjugate of a flat vector bundle $E$.
Let $\bar b$ denote the complex conjugate of a fiber wise non-degenerate
symmetric bilinear form on $E$. Via the canonic isomorphism of complex lines
$\det H^*(M;\bar E)=\overline{\det H^*(M;E)}$ we obviously have
\begin{equation}\label{E:anconj}
\tau^\an_{\bar E,\bar\e^*,[\bar b]}=\overline{\tau^\an_{E,\e^*,[b]}}
\end{equation}
where $\e^*\mapsto\bar\e^*$ denotes the complex conjugation of coEuler structures 
introduced in Section~\ref{S:eul}. For this note that $\Delta_{\bar E,g,\bar b}=
\Delta_{E,g,b}$ but the spectrum of $\Delta_{\bar E,g,\bar b}$ is complex
conjugate to the spectrum of $\Delta_{E,g,b}$ and thus ${\det}'(\Delta_{\bar E,g,\bar b,q})
=\overline{{\det}'(\Delta_{E,g,b,q})}$. Also recall \eqref{E:KTconj}.

Suppose $F$ is a flat real vector bundle over $M$.
Let $\e^*\in\Eul(M;\R)$ be a coEuler structure with real coefficients.
Let $h$ be a fiber wise non-degenerate symmetric bilinear form on $F$.
Proceeding exactly as in the complex case we obtain a non-degenerate
bilinear form $\tau^\an_{F,\e,[h]}$ on the real line $\det H^*(M;F)$.
Note that although the Laplacians $\Delta_{F,g,h}$ need not be selfadjoint
their spectra are invariant under complex conjugation and hence
${\det}'(\Delta_{F,g,h,q})$ will be real.
Let $F^\C$ denote the complexification of the flat bundle $F$, and let
$h^\C$ denote the complexification of $h$, a non-degenerate symmetric
bilinear form on $F^\C$. Via the canonic isomorphism
of complex lines $\det H^*(M;F^\C)=(\det H^*(M;F))^\C$ we have:
\begin{equation}\label{E:cxan}
\tau^\an_{F^\C,\e^*,[h^\C]}=\bigl(\tau^\an_{F,\e^*,[h]}\bigr)^\C
\end{equation}
For this note that via $\Omega^*(M;F^\C)=\Omega^*(M;F)^\C$ we have
$\Delta_{F^\C,g,h^\C}=(\Delta_{F,g,h})^\C$ and thus
${\det}'(\Delta_{F^\C,g,h^\C,q})={\det}'(\Delta_{F,g,h,q})$, and also
recall \eqref{E:KTcx}. If $n$ is odd, $H^*(M;F)=0$, and if $h$ is positive
definite, then $\tau^\an_{F,\e^*_\can,[h]}$ is the square of the analytic torsion
considered in \cite{RS71}, see Remark~\ref{R:cancoeul}.

\begin{remark}\label{R:nob}
Not every flat complex vector bundle $E$ admits a fiber wise non-degenerate
symmetric bilinear form $b$. However, since $E$ is flat  
all rational Chern classes of $E$ must vanish. Since $M$ is compact,
the Chern character induces an isomorphism on rational $K$-theory, and hence
$E$ is trivial in rational
$K$-theory. Thus there exists $N\in\N$ so that $E^N=E\oplus\cdots\oplus E$
is a trivial vector bundle. Particularly, there exists a fiber wise
non-degenerate bilinear form $b$ on $E^N$. In view of \eqref{E:578}
the non-degenerate bilinear form 
$\bigl(\tau^\an_{E^N,\e^*,[b]}\bigr)^{1/N}$ on $\det H^*(M;E)$ is a
reasonable candidate for the analytic torsion of $E$.
Note however, that this is only defined up to a root of unity.
\end{remark}

\subsection*{Rewriting the analytic torsion}

Instead of just treating the $0$-eigen space by means of finite 
dimensional linear algebra one can equally well do this with finitely many
eigen spaces of $\Delta_{E,g,b}$. Proposition~\ref{P:A} below makes this precise.
We will make use of this formula when computing the variation of the analytic
torsion through a variation of $g$ and $b$. This is necessary since the
dimension of the $0$-eigen space need not be locally constant through such a 
variation. Note that this kind of problem does not occur in the
selfadjoint situation, i.e.\ when instead of a non-degenerate symmetric 
bilinear form we have a hermitian structure.

Suppose $\gamma$ is a simple closed curve around $0$,
avoiding the spectrum of $\Delta_{E,g,b}$. Let
$\Omega^*_{g,b}(M;E)(\gamma)$
denote the sum of eigen spaces corresponding to eigen values in the interior
of $\gamma$. Using Proposition~\ref{P:spec} we see that the inclusion 
$\Omega^*_{g,b}(M;E)(\gamma)\to\Omega^*(M;E)$
is a quasi isomorphism. We obtain a canonic isomorphism of determinant
lines
\begin{equation}\label{E:987}
\det H(\Omega^*_{g,b}(M;E)(\gamma))=\det H^*(M;E).
\end{equation}
Moreover, the restriction of $\beta_{g,b}$ to $\Omega^*_{g,b}(M;E)(\gamma)$
is non-degenerate. Hence the torsion provides us with a
non-degenerate bilinear form on $\det H(\Omega^*_{g,b}(M;E)(\gamma))$ and via
\eqref{E:987} we get a non-degenerate bilinear form
$\tau^\an_{E,g,b}(\gamma)$ on $\det H^*(M;E)$.
Moreover, introduce
$$
{\det}^\gamma(\Delta_{E,g,b,q})
:=\exp\left(-\frac\partial{\partial s}\Big|_{s=0}
\tr\Bigl(\bigl(\Delta_{E,g,b,q}|{\Omega^q_{g,b}(M;E)(\gamma)^{\perp_{\beta_{g,b}}}}\bigr)^{-s}\Bigr)\right),
$$
the zeta regularized product of eigen values in the exterior of $\gamma$.

\begin{proposition}\label{P:A}
In this situation, as bilinear forms on $\det H^*(M;E)$, we have:
$$
\tau^\an_{E,g,b}(0)\cdot\prod_q\bigl({\det}'(\Delta_{E,g,b,q})\bigr)^{(-1)^qq}
=\tau^\an_{E,g,b}(\gamma)\cdot\prod_q\bigl({\det}^\gamma(\Delta_{E,g,b,q})\bigr)^{(-1)^qq}
$$
\end{proposition}

\begin{proof}
Let $C^*\subseteq\Omega^*_{g,b}(M;E)(\gamma)$ denote the sum of the eigen spaces of
$\Delta_{E,g,b}$ corresponding to non-zero eigen values in the interior of $\gamma$.
Clearly, for every $q$ we have
$$
{\det}'(\Delta_{E,g,b,q})=\det(\Delta_{E,g,b,q}|_{C^q})\cdot{\det}^\gamma(\Delta_{E,g,b,q}).
$$
Particularly,
\begin{equation}\label{E:96}
\prod_q\bigl({\det}'(\Delta_{E,g,b,q})\bigr)^{(-1)^qq}
=\prod_q\bigl(\det(\Delta_{E,g,b,q}|_{C^q})\bigr)^{(-1)^qq}
\cdot\prod_q\bigl({\det}^\gamma(\Delta_{E,g,b,q})\bigr)^{(-1)^qq}.
\end{equation}
Applying Lemma~\ref{L:fintor} to the finite dimensional complex
$\Omega^*_{g,b}(M;E)(\gamma)$ we obtain
\begin{equation}\label{E:97}
\tau^\an_{E,g,b}(\gamma)=\tau^\an_{E,g,b}(0)
\cdot\prod_q\bigl(\det(\Delta_{E,g,b,q}|_{C^q})\bigr)^{(-1)^qq}
\end{equation}
Multiplying \eqref{E:96} with $\tau^\an_{E,g,b}(0)$ and using \eqref{E:97} 
we obtain the statement.
\end{proof}

\section{A Bismut--Zhang, Cheeger, M\"uller type formula}\label{S:BZ}

The conjecture below asserts that the complex valued analytical torsion
defined in Section~\ref{S:anator} coincides with the combinatorial torsion
from Section~\ref{S:comb}. It should be considered as a complex valued
version of a theorem of Cheeger \cite{C77, C79}, M\"uller \cite{M78}
and Bismut--Zhang \cite{BZ92}.

\begin{conjecture}\label{C:main}
Let $M$ be a closed connected smooth manifold with vanishing
Euler--Poincar\'e characteristics. Let $E$ be a flat complex vector bundle
over $M$, and suppose $b$ is a fiber wise non-degenerate symmetric bilinear form on
$E$. Let $\e\in\Eul(M;\Z)$ be an Euler structure. Then, as bilinear forms
on the complex line $\det H^*(M;E)$, we have:
$$
\tau^\comb_{E,\e}=\tau^\an_{E,P(\e),[b]}
$$
Here we slightly abuse notation and let $P$ also denote the composition
$\Eul(M;\Z)\to\Eul(M;\C)\xrightarrow{P}\Eul^*(M;\C)$, see Section~\ref{S:eul}.
\end{conjecture}

We will establish this conjecture in several special cases, see
Remark~\ref{R:evendim}, Theorem~\ref{T:abs1}, Corollary~\ref{C:locconst}, Corollary~\ref{C:complexi} 
and the discussion for the circle below. Some of these results have been established by Braverman--Kappeler
\cite{BK06} and were not contained in the first version of this manuscript.
The proofs we provide below have been inspired by a trick used in \cite{BK06}
but do not rely on the results therein.

\begin{remark}
If Conjecture~\ref{C:main} holds for one Euler structure $\e\in\Eul(M;\Z)$
then it will hold for all Euler structures. This follows immediately from
Remark~\ref{R:aneul}, Remark~\ref{R:combeul} and Lemma~\ref{L:A}.
\end{remark}

\begin{remark}\label{R:bindep}
If Conjecture~\ref{C:main} holds, and if $\e^*\in\Eul^*(M;\C)$ is integral, then
$\tau^\an_{E,\e^*,[b]}$ is independent of $[b]$. This is not obvious from
the definition of the analytic torsion.
\end{remark}

\begin{remark}\label{R:X}
If Conjecture~\ref{C:main} holds, $\e\in\Eul(M;\Z)$ and $\e^*\in\Eul^*(M;\C)$
then:
$$
\tau^\comb_{E,\e}=\tau^\an_{E,\e^*,[b]}\cdot
\Bigl(e^{\langle[\omega_{E,b}]\cup(P(\e)-\e^*),[M]\rangle}\Bigr)^2
$$
This follows from Remark~\ref{R:aneul}.
\end{remark}

\begin{remark}
If Conjecture~\ref{C:main} holds, and $\e^*\in\Eul^*(M;\C)$, then
$\tau^\an_{E,\e^*,[b]}$ does only depend on $E$, $\e^*$ and the induced
homotopy class $[\det b]$ of non-degenerate bilinear forms on $\det E$.
This follows from Remark~\ref{R:X} and the fact that the cohomology class
$[\omega_{E,b}]$ does depend on $\det E$ and the homotopy class $[\det b]$ on $\det E$
only, see Section~\ref{S:eul}.
\end{remark}

\subsection*{Relative torsion}

In the situation above, consider the non-vanishing complex number
$$
\mathcal S_{E,\e,[b]}:=\frac{\tau^\an_{E,P(\e),[b]}}{\tau^\comb_{E,\e}}\in\C^\times.
$$
It follows from Remark~\ref{R:aneul}, Remark~\ref{R:combeul} and Lemma~\ref{L:A}
that this does not depend on $\e\in\Eul(M;\Z)$. We will thus denote it by $\mathcal S_{E,[b]}$.
The number $\mathcal S_{E,[b]}$ will be referred to as the \emph{relative torsion}
associated with the flat complex vector bundle $E$ and the homotopy class $[b]$.
Conjecture~\ref{C:main} asserts that $\mathcal S_{E,[b]}=1$.

Similarly, if $F$ is a real flat vector bundle
over $M$ equipped with a fiber wise non-degenerate symmetric bilinear form $h$, we set
$$
\mathcal S_{F,[h]}:=\frac{\tau^\an_{F,P(\e),[h]}}{\tau^\comb_{F,\e}}\in\R^\times:=\R\setminus\{0\}
$$
where $\e\in\Eul(M;\Z)$ is any Euler structure. The combinatorial torsion
$\tau^\comb_{F,\e}$ and the analytic torsion $\tau^\an_{F,P(\e),[h]}$ on $\det H^*(M;F)$
have been introduced in Sections~\ref{S:comb} and \ref{S:anator}, respectively.
It follows via complexification from the corresponding statements for complex vector bundles that
this does indeed only depend on $F$ and the homotopy class of $h$, see \eqref{E:cxcomb} and \eqref{E:cxan}.

\begin{remark}\label{R:BZ}
If $F$ is a flat real vector bundle equipped with a positive definite symmetric bilinear form $h$,
then the Bismut--Zhang theorem \cite[Theorem~0.2]{BZ92} asserts that $\mathcal S_{F,[h]}=1$.
This follows from the formula in Proposition~\ref{P:Sexp} below (applied to a simple closed curve whose interior
contains the eigen value $0$ only) which, via complexification, provides
an analogous formula for flat real vector bundles. For the relation of the first factor in this formula
with the statement in \cite[Theorem~0.2]{BZ92} see \eqref{E:ZZ}.
\end{remark}

\begin{proposition}\label{P:S}
The following properties hold:
\begin{enumerate}
\item\label{P:S:i}
$\mathcal S_{E_1\oplus E_2,[b_1\oplus b_2]}=
\mathcal S_{E_1,[b_1]}\cdot\mathcal S_{E_2,[b_2]}$
\item\label{P:S:ii}
$\mathcal S_{E'\otimes\mathcal O_M,[b']}=(\mathcal S_{E,[b]})^{(-1)^{n+1}}$
\item\label{P:S:iii}
$\mathcal S_{\bar E,[\bar b]}=\overline{\mathcal S_{E,[b]}}$
\item\label{P:S:iv}
$\mathcal S_{F^\C,[h^\C]}=\mathcal S_{F,[h]}$
\end{enumerate}
\end{proposition}

\begin{proof}
This follows immediately from the basic properties of analytic and 
combinatorial torsion discussed in Sections~\ref{S:comb} and 
\ref{S:anator}. For \itemref{P:S:ii} and \itemref{P:S:iii} one also has to use
$P(\nu(\e))=\nu(P(\e))$ and $\overline{P(\e)}=P(\bar\e)$, see Section~\ref{S:eul}.
\end{proof}

\begin{remark}\label{R:evendim}
Proposition~\ref{P:S}\itemref{P:S:ii} permits to verify Conjecture~\ref{C:main},
up to sign, for even dimensional orientable manifolds and parallel bilinear forms.
More precisely, let $M$ be an even dimensional closed connected orientable smooth manifold with
vanishing Euler--Poincar\'e characteristics. Let $E$ be a flat complex vector
bundle over $M$ and suppose $b$ is a parallel fiber wise non-degenerate
symmetric bilinear form on $E$. Let $\e\in\Eul(M;\Z)$ be an Euler structure.
Then 
\begin{equation}\label{E:291}
\tau^\comb_{E,\e}=\pm\tau^\an_{E,P(\e),[b]}
\end{equation} 
i.e.\ in this situation Conjecture~\ref{C:main} holds up to sign.
To see this, note that the parallel bilinear form $b$ and the choice of an orientation
provides an isomorphism of flat vector bundles $b:E\to E'\otimes\mathcal O_M$ which maps
$b$ to $b'$. Thus $\mathcal S_{E'\otimes\mathcal O_M,[b']}=\mathcal S_{E,[b]}$. 
Combining this with Proposition~\ref{P:S}\itemref{P:S:ii} we obtain
$(\mathcal S_{E,[b]})^2=1$, and hence \eqref{E:291}. Note, however, that in this situation the 
arguments used to establish \eqref{E:dualan} immediately yield
$$
\prod_q\bigl({\det}'(\Delta_{E,g,b,q})\bigr)^{(-1)^qq}=1.
$$
\end{remark}

Corollary~\ref{C:bbb} below has been established by Braverman and Kappeler
see \cite[Theorem~5.3]{BK06} by comparing $\tau^\an_{E,P(\e),[b]}$ with their refined analytic torsion,
see \cite[Theorem~1.4]{BK06}. We will give an elementary proof relying on Proposition~\ref{P:S}
and a trick similar to the one used in the proof of Theorem~1.4 in \cite{BK06}.

\begin{corollary}\label{C:bbb}
Let $M$ be a closed connected smooth orientable manifold of odd dimension. Suppose
$E$ is a flat complex vector bundle over $M$ equipped with a non-degenerate
symmetric bilinear form $b$. Let $\e^*\in\Eul^*(M;\mathbb C)$ be an integral 
coEuler structure. Then, up to sign, $\tau^\an_{E,\e^*,[b]}$ is independent 
of $[b]$, cf.\ Remark~\ref{R:bindep}.
\end{corollary}

\begin{proof}
It suffices to show $(\mathcal S_{E,[b]})^2$ is independent of $[b]$.
The choice of an orientation provides an isomorphism of flat vector bundles
$E'\cong E'\otimes\mathcal O_M$ from which we obtain
$$
\mathcal S_{E',[b']}=\mathcal S_{E'\otimes\mathcal O_M,[b']}
=\mathcal S_{E,[b]}
$$
where the latter equality follows from Proposition~\ref{P:S}\itemref{P:S:ii}.
Together with Proposition~\ref{P:S}\itemref{P:S:i} we thus obtain
\begin{equation}\label{E:1}
(\mathcal S_{E,[b]})^2
=\mathcal S_{E,[b]}\cdot\mathcal S_{E',[b']}
=\mathcal S_{E\oplus E',[b\oplus b']}.
\end{equation}
Observe that on $E\oplus E'$ there exists a canonic (independent of $b$)
symmetric non-degenerate bilinear form $b_\can$ defined by
$$
b_\can\bigl((x_1,\alpha_1),(x_2,\alpha_2)\bigr)
:=\alpha_1(x_2)+\alpha_2(x_1),\qquad x_1,x_2\in E,\ \alpha_1,\alpha_2\in E'.
$$
This bilinear form $b_\can$ is homotopic to $b\oplus b'$, and thus
$$
\mathcal S_{E\oplus E',[b\oplus b']}=
\mathcal S_{E\oplus E',[b_\can]}.
$$
Hence $\mathcal S_{E\oplus E',[b\oplus b']}$ does not depend on $[b]$. 
In view of \eqref{E:1} the same holds for $(\mathcal S_{E,[b]})^2$, and the 
proof is complete.

To see that $b\oplus b'$ is indeed homotopic to $b_\can$ let us consider $b$
as an isomorphism $b:E\to E'$. For $t\in\mathbb R$ consider the endomorphisms
$$
\Phi_t:\eend(E\oplus E'),\qquad
\Phi_t:=\left(\begin{smallmatrix}\id_E\cos t&-b^{-1}\sin t\\b\sin t&\id_{E'}\cos t\end{smallmatrix}\right)
$$
From $\Phi_{t+s}=\Phi_t\Phi_s$ we conclude that every $\Phi_t$ is invertible.
Consider the curve of non-degenerate symmetric bilinear forms
$$
b_t:=\Phi_t^*b_\can,\qquad b_t(X_1,X_2)=b_\can(\Phi_tX_1,\Phi_tX_2),
\qquad X_1,X_2\in E\oplus E'.
$$
Then clearly $b_0=b_\can$. An easy calculation shows $b_{\pi/4}=b\oplus(-b')$.
Clearly, $b\oplus(-b')$ is homotopic to $b\oplus b'$. So we see that $b_\can$
is homotopic to $b\oplus b'$.
\end{proof}

Using a result of Cheeger \cite{C77, C79}, M\"uller \cite{M78} and Bismut--Zhang \cite{BZ92} we 
will next show that the absolute value of the relative torsion is always one.
In odd dimensions this has been established by Braverman and Kappeler, see Theorem~1.10 in \cite{BK06}.
We will again use a trick similar to the one in \cite{BK06}.

\begin{theorem}\label{T:abs1}
Suppose $M$ is a closed connected smooth manifold with vanishing  Euler--Poincar\'e characteristics. 
Let $E$ be a flat complex vector bundle over $M$
equipped with a non-degenerate symmetric bilinear form $b$.
Then $|\mathcal S_{E,[b]}|=1$.
\end{theorem}

\begin{proof}
Note first that in view of Proposition~\ref{P:S}\itemref{P:S:iii} and \itemref{P:S:i} we have
\begin{equation}\label{E:2}
|\mathcal S_{E,[b]}|^2
=\mathcal S_{E,[b]}\cdot\overline{\mathcal S_{E,[b]}}
=\mathcal S_{E\oplus\bar E,[b\oplus\bar b]}.
\end{equation}
Set $k:=\rank E$, and 
observe that $b$ provides a reduction of the structure group of $E$ to 
$O_k(\mathbb C)$. Since the inclusion $O_k(\mathbb R)\subseteq O_k(\mathbb C)$ is a homotopy
equivalence, the structure group can thus be further reduced to $O_k(\mathbb R)$.
In other words, there exists a complex anti-linear involution $\nu:E\to E$
such that
$$ 
\nu^2=\id_E,\qquad
b(\nu x,y)=\overline{b(x,\nu y)},\qquad b(x,\nu x)\geq0,\qquad x,y\in E.
$$ 
Then 
$$
\mu:E\otimes E\to\mathbb C,\qquad \mu(x,y):=b(x,\nu y)
$$
is a fiber wise positive definite Hermitian structure on $E$, anti-linear in the second variable.
Define a non-degenerate symmetric bilinear form $b^\mu$ on $E\oplus\bar E$ by
$$
b^\mu\bigl((x_1,y_1),(x_2,y_2)\bigr)
:=\mu(x_1,y_2)+\mu(x_2,y_1).
$$
We claim that the symmetric bilinear form $b^\mu$ is homotopic to $b\oplus\bar b$.
To see this, consider $\nu:E\to\bar E$ as a complex linear isomorphism. 
For $t\in\mathbb R$, define
$$
\Phi_t\in\eend(E\oplus\bar E),\qquad
\Phi_t:=\left(\begin{smallmatrix}\id_E\cos t&-\nu^{-1}\sin t\\\nu\sin t&\id_{\bar E}\cos t\end{smallmatrix}\right)
$$
From $\Phi_{t+s}=\Phi_t\Phi_s$ we conclude that every $\Phi_t$ is invertible.
Consider the curve of non-degenerate symmetric bilinear forms
$$
b_t:=\Phi_t^*b^\mu,\qquad b_t(X_1,X_2)=b^\mu(\Phi_tX_1,\Phi_tX_2),\qquad
X_1,X_2\in E\oplus\bar E.
$$
Clearly, $b_0=b^\mu$. An easy computation shows $b_{\pi/4}=b\oplus(-\bar b)$.
Since $b\oplus(-\bar b)$ is homotopic to $b\oplus\bar b$ we see that $b^\mu$ is indeed
homotopic to $b\oplus\bar b$. Together with \eqref{E:2} we conclude
\begin{equation}\label{E:3}
|\mathcal S_{E,[b]}|^2=\mathcal S_{E\oplus\bar E,b^\mu}.
\end{equation}
Next, recall that there is a canonic isomorphism of flat vector bundles 
$$
\psi:E^{\mathbb C}\cong E\oplus\bar E,\qquad
\psi(x+\ii y):=(x+\ii y,x-\ii y),\qquad x,y\in E.
$$
Consider the fiber wise positive definite symmetric real bilinear form $h:=\Re\mu$ on 
$E^{\mathbb R}$, the underlying real vector bundle. Its complexification 
$h^{\mathbb C}$ is a non-degenerate
symmetric bilinear form on $E^{\mathbb C}$.
A simple computations shows $\psi^*b^\mu=2h^{\mathbb C}$.
Together with \eqref{E:3} we obtain
$$
|\mathcal S_{E,[b]}|^2=\mathcal S_{E^{\mathbb C},2h^{\mathbb C}}
=\mathcal S_{E^{\mathbb R},2h}
$$
where the last equation follows from Proposition~\ref{P:S}\itemref{P:S:iv}.
The Bismut--Zhang theorem \cite[Theorem~0.2]{BZ92} 
asserts that $\mathcal S_{E^{\mathbb R},2h}=1$, see Remark~\ref{R:BZ}, and the proof is complete.
\end{proof}

\subsection*{Analyticity of the relative torsion}

In this section we will show that the relative torsion $\mathcal S_{E,[b]}$ depends holomorphically on
the flat connection, see Proposition~\ref{P:holoS} below. Combined with Theorem~\ref{T:abs1} this 
implies that $\mathcal S_{E,[b]}$ is locally constant on the space of flat connections on a fixed vector bundle, see
Corollary~\ref{C:locconst} below. We start by establishing an explicit formula for the relative torsion, see Proposition~\ref{P:Sexp}.

Suppose $f:C_1\to C_2$ is a homomorphism of finite dimensional complexes.
Consider the mapping cone $C_2^{*-1}\oplus C_1^*$ with differential
$\left(\begin{smallmatrix}d&f\\0&-d\end{smallmatrix}\right)$. If $C_1^*$ and
$C_2^*$ are equipped with graded non-degenerate symmetric bilinear forms $b_1$ and
$b_2$ we equip the mapping cylinder with the bilinear form $b_2\oplus b_1$.
The resulting torsion $\tau(f,b_1,b_2):=\tau_{C_2^{*-1}\oplus
C_1^*,b_2\oplus b_1}$ is called the relative torsion of $f$. It is a
non-degenerate bilinear form on the determinant line 
$\det H(C_2^{*-1}\oplus C_1^*)$. Recall that if $f$ is a quasi isomorphism
then $C_2^{*-1}\oplus C_1^*$ is acyclic and
\begin{equation}\label{E:AR}
\tau(f,b_1,b_2)=\frac{(\det H(f))(\tau_{C_1^*,b_1})}{\tau_{C_2^*,b_2}}
\end{equation}
where $\det H(f):\det H(C_1^*)\to\det H(C_2^*)$ denotes the isomorphism of
complex lines induces from the isomorphism in cohomology $H(f):H(C_1^*)\to
H(C_2^*)$.

Let us apply this to the integration homomorphism
\begin{equation}\label{E:AX}
\Int:\Omega^*_{g,b}(M;E)(\gamma)\to C^*(X;E)
\end{equation} 
where the notation is as in 
Proposition~\ref{P:A}. Equip $\Omega^*_{g,b}(M;E)(\gamma)$ with the
restriction of $\beta_{g,b}$, and equip $C^*(X;E)$ with the bilinear form
$b|_{\mathcal X}$ obtained by restricting $b$ to the fibers over $\mathcal X$.
Since \eqref{E:AX} is a quasi isomorphism the mapping cylinder is acyclic and 
the corresponding relative torsion is a non-vanishing complex number we will denote
by
$$
\tau\Bigl(\Omega_{g,b}^*(M;E)(\gamma)\xrightarrow{\Int}C^*_b(X;E)\Bigr)\in\C^\times.
$$

\begin{proposition}\label{P:Sexp}
Let $M$ be a closed connected smooth manifold with vanishing
Euler--Poincar\'e characteristics. Let $E$ be a flat complex vector bundle
over $M$. Let $g$ be a Riemannian metric, and let
$X$ be a Morse--Smale vector field on $M$. Suppose $b$ is a fiber wise
non-degenerate symmetric bilinear form on $E$ which is parallel
in a neighborhood of the critical points $\mathcal X$. 
Moreover, let $\gamma$ be a simple closed curve around $0$ which avoids the
spectrum of $\Delta_{E,g,b}$. Then:
\begin{multline*}
\mathcal S_{E,[b]}=
\tau\Bigl(\Omega^*_{g,b}(M;E)(\gamma)\xrightarrow{\Int}C_b^*(X;E)\Bigr)
\\\cdot\prod_q\bigl({\det}^\gamma(\Delta_{E,g,b,q})\bigr)^{(-1)^qq}
\cdot\exp\Bigl(-2\int_{M\setminus\mathcal X}\omega_{E,b}\wedge(-X)^*\Psi(g)\Bigr)
\end{multline*}
The integral is absolutely convergent since $\omega_{E,b}$ vanishes
in a neighborhood of $\mathcal X$.
\end{proposition}

\begin{proof}
Let $x_0\in M$ be a base point. For every critical point $x\in\mathcal X$
choose a path $\sigma_x$ with $\sigma_x(0)=x_0$ and $\sigma_x(1)=x$.
Set $c:=\sum_{x\in\mathcal X}(-1)^{\ind(x)}\sigma_x$ and consider the
Euler structure $\e:=[-X,c]\in\Eul(M;\Z)$.
For the dual coEuler structure $P(\e)=[g,\alpha]$ we have, see \eqref{E:665},
\begin{equation}\label{E:444}
\int_{M\setminus\mathcal X}\omega_{E,b}\wedge\bigl((-X)^*\Psi(g)-\alpha\bigr)=\int_c\omega_{E,b}.
\end{equation}

Let $b_{x_0}$ denote the bilinear form on the fiber $E_{x_0}$ obtained by restricting $b$.
For $x\in\mathcal X$ let $\tilde b_x$ denote the bilinear form obtained from
$b_{x_0}$ by parallel transport along $\sigma_x$.
Let $\tilde b_{\det C^*(X;E)}$ denote the induced bilinear form on
$\det C^*(X;E)$. This is the bilinear form used in the definition of the
combinatorial torsion. We want to compare it with the bilinear form
$b_{\det C^*(X;E)}$ on $\det C^*(X;E)$ induced by the restriction
$b|_{\mathcal X}$ of $b$ to the
fibers over $\mathcal X$. A simple computation similar to the proof of
Lemma~\ref{L:A} yields
\begin{equation}\label{E:333}
\tilde b_{\det C^*(X;E)}=\exp\Bigl(2\int_c\omega_{E,b}\Bigr)
\cdot b_{\det C^*(X;E)}.
\end{equation}

Let $\tau_{C^*(X;E),b|_{\mathcal X}}$ denote the non-degenerate
bilinear form on $\det H^*(M;E)$ obtained from the torsion of the complex
$C^*(X;E)$ equipped with the bilinear form $b|_{\mathcal X}$ via the
isomorphism $\det H^*(M;E)=\det H(C^*(X;E))$, see
\eqref{E:int} and \eqref{E:rm}. 
Then, using \eqref{E:AR},
\begin{equation}\label{E:ZZ}
\frac{\tau^\an_{E,g,b}(\gamma)}{\tau_{C^*(X;E),b|_{\mathcal X}}}
=\tau\Bigl(\Omega^*_{g,b}(M;E)(\gamma)\xrightarrow{\Int}C^*_b(X;E)\Bigr).
\end{equation}
Moreover, \eqref{E:333} implies
\begin{equation}\label{E:AB}
\tau^\comb_{E,\e}
=\tau_{C^*(X;E),b|_{\mathcal X}}\cdot\exp\Bigl(2\int_c\omega_{E,b}\Bigr).
\end{equation}
From Proposition~\ref{P:A} we obtain
\begin{equation}\label{E:AC}
\tau^\an_{E,P(\e),[b]}=\tau^\an_{E,g,b}(\gamma)
\cdot\prod_q\bigl({\det}^\gamma(\Delta_{E,g,b,q})\bigr)^{(-1)^qq}\cdot
\exp\Bigl(-2\int_M\omega_{E,b}\wedge\alpha\Bigr).
\end{equation}
Combining \eqref{E:444}, \eqref{E:ZZ}, \eqref{E:AB} and \eqref{E:AC} we obtain the
statement of the proposition.
\end{proof}

Consider an open subset $U\subseteq\C$ and a family
of flat complex vector bundles $\{E^z\}_{z\in U}$. Such a family is called 
\emph{holomorphic} if the underlying vector bundles are the same for all $z\in U$ 
and the mapping $z\mapsto\nabla^{E^z}$ is holomorphic into the affine
Fr\'echet space of linear connections equipped with the $C^\infty$-topology.

\begin{proposition}\label{P:holoS}
Let $M$ be a closed connected smooth manifold with vanishing
Euler--Poincar\'e characteristics. Let $\{E^z\}_{z\in U}$ be a holomorphic
family of flat complex vector bundles over $M$, and let $b^z$ be a holomorphic 
family of fiber wise non-degenerate symmetric bilinear forms on $E^z$. 
Then $\mathcal S_{E^z,[b^z]}$ depends holomorphically on $z$.
\end{proposition}

\begin{proof}
Let $X$ be a Morse--Smale vector field on $M$. Let $g$ be a Riemannian
metric on $M$.
In view of Theorem~\ref{T:anom} we may w.l.o.g.\ assume $\nabla^{E^z}b^z=0$
in a neighborhood of $\mathcal X$.
W.l.o.g.\ we may assume that there exists a simple closed curve $\gamma$ around $0$
so that the spectrum of $\Delta_{E^z,g,b^z}$ avoids $\gamma$ for all
$z\in U$. From Proposition~\ref{P:Sexp} we know:
\begin{multline*}
\mathcal S_{E^z,[b^z]}=
\tau\Bigl(\Omega^*_{g,b^z}(M;E^z)(\gamma)\xrightarrow{\Int} C^*_{b^z}(X;E^z)\Bigr)
\\\cdot\prod_q\bigl({\det}^\gamma(\Delta_{E^z,g,b^z,q})\bigr)^{(-1)^qq}
\cdot\exp\Bigl(-2\int_{M\setminus\mathcal X}\omega_{E^z,b^z}\wedge(-X)^*\Psi(g)\Bigr)
\end{multline*}
Since $\Delta_{E^z,g,b^z}$ depends holomorphically on $z$, each of the three
factors in this expression for $S_{E^z,[b^z]}$ will depend holomorphically
on $z$ too.
\end{proof}

In odd dimensions the following result has been established by Braverman and Kappeler, see Theorem~1.10 in \cite{BK06}.

\begin{corollary}\label{C:locconst}
Let $M$ be a closed connected smooth manifold with vanishing Euler--Poincar\'e characteristics.
Let $E$ be a complex vector bundle over $M$, and let $b$ be a fiber wise non-degenerate
symmetric bilinear form on $E$. Then the assignment $\nabla\mapsto\mathcal S_{(E,\nabla),[b]}$ is locally constant,
and of absolute value one, on the space of flat connections on $E$.
\end{corollary}

\begin{proof}
Note that in view of Theorem~\ref{T:abs1} and Proposition~\ref{P:holoS} the relative torsion
$\mathcal S_{(E,\nabla^z),[b]}$ is constant along every holomorphic path of flat connections $z\mapsto\nabla^z$ on $E$.
Moreover, note that two flat connections, contained in the same connected component, can always
be joined by a piecewise holomorphic path of flat connections. 
\end{proof}

Using the Bismut--Zhang, Cheeger, M\"uller theorem again, we are able to verify Conjecture~\ref{C:main}
for flat connections contained in particular connected components of the space of flat connections on a fixed
complex vector bundle. More precisely, we have\footnote{In a recent preprint \cite{H06} R.-T.~Huang verified a similar 
statement for flat connections whose connected component contains a flat connection which admits a parallel Hermitian structure.}

\begin{corollary}\label{C:complexi}
Let $M$ be a closed connected smooth manifold with vanishing Euler--Poincar\'e characteristics.
Let $(F,\nabla^F)$ be a flat real vector bundle over $M$ equipped with a fiber wise Hermitian structure $h$.
Let $(E,\nabla^E)$ denote the flat complex vector bundle obtained by complexifying $(F,\nabla^F)$, and let
$b$ denote the fiber wise non-degenerate symmetric bilinear form on $E$ obtained by complexifying $h$.
Then, for every flat connection $\nabla$ on $E$ which is contained in the connected component of $\nabla^E$, we have
$\mathcal S_{(E,\nabla),[b]}=1$.
\end{corollary}

\begin{proof}
In view of Corollary~\ref{C:locconst} it suffices to show $\mathcal S_{(E,\nabla^E),[b]}=1$. From Proposition~\ref{P:S}\itemref{P:S:iv}
we have $\mathcal S_{(E,\nabla^E),[b]}=\mathcal S_{(F,\nabla^F),[h]}$. In view of \cite[Theorem~0.2]{BZ92}, see Remark~\ref{R:BZ}, 
we indeed have $\mathcal S_{(F,\nabla^F),[h]}=1$, and the statement follows.
\end{proof}

\subsection*{The circle, a simple explicit example}

Consider $M:=S^1$. In this case it is possible to explicitly compute
the combinatorial and analytic torsion, see below. It turns out that
Conjecture~\ref{C:main} holds true for every flat vector bundle over the
circle.

We think of $S^1$ as $\{z\in\C\mid|z|=1\}$. Equip
$S^1$ with the standard Riemannian metric $g$ of circumference $2\pi$. Orient
$S^1$ in the standard way. Let $\theta$ denote the angular
\lq coordinate\rq. Let $\frac\partial{\partial\theta}$ denote the
corresponding vector field which is of length $1$ and induces the
orientation. For the dual $1$-form we write $d\theta$.

Let $k\in\N$ and suppose $a\in C^\infty(S^1,\gl_k(\C))$. Let $E^a$ denote the
trivial vector bundle $S^1\times\C^k$ equipped with the flat connection
$\nabla=\frac\partial{\partial\theta}+a$. Here and in what follows we use
the identifications $\Omega^0(M;E^a)=C^\infty(S^1;\C^k)=\Omega^1(M;E^a)$
where the latter stems from the global coframe $d\theta$.

Let $b\in C^\infty(S^1,\Sym^\times_k(\C))$ where $\Sym^\times_k(\C)$
denotes the space of complex non-de\-ge\-ner\-ate symmetric $k\times k$-matrices.
We consider $b$ as a fiber wise non-degenerate symmetric bilinear form on
$E^a$. For the induced bilinear form on $\Omega^*(S^1;E^a)$ we have:
\begin{align*}
\beta_{g,b}(v,w)&=\int_{S^1}v^tbw\ d\theta, &v,w\in\Omega^0(S^1;E^a)=C^\infty(S^1,\C^k)
\\
\beta_{g,b}(v,w)&=\int_{S^1}v^tbw\ d\theta, &v,w\in\Omega^1(S^1;E^a)=C^\infty(S^1,\C^k)
\end{align*}
A straight forward computations yields:
\begin{align*}
d_{E^a}&=\tfrac\partial{\partial\theta}+a
\\
d^\sharp_{E^a,g,b}&=-\tfrac\partial{\partial\theta}-b^{-1}b'+b^{-1}a^tb
\\
\Delta_{E^a,g,b,0}&=-\bigl(\tfrac\partial{\partial\theta}\bigr)^2
+\bigl(b^{-1}a^tb-b^{-1}b'-a\bigr)\tfrac\partial{\partial\theta}
+\bigl(b^{-1}a^tba-b^{-1}b'a-a'\bigr)
\\
\Delta_{E^a,g,b,1}&=-\bigl(\tfrac\partial{\partial\theta}\bigr)^2
+\bigl(b^{-1}a^tb-b^{-1}b'-a\bigr)\tfrac\partial{\partial\theta}
\\&\qquad+\bigl((b^{-1}a^tb)'-(b^{-1}b')'-ab^{-1}b'+ab^{-1}a^tb\bigr)
\\
b^{-1}\nabla_{\frac\partial{\partial\theta}}b
&=b^{-1}b'-b^{-1}a^tb-a
\\
\omega_{E^a,b}
&=-\tfrac12\tr\bigl(b^{-1}b'-b^{-1}a^tb-a\bigr)d\theta
=-\tfrac12\bigl(\tr(b^{-1}b')-2\tr(a)\bigr)d\theta
\end{align*}
Here $b':=\frac\partial{\partial\theta}b$ and $a':=\frac\partial{\partial\theta}a$.

Let us write $A\in\GL_k(\C)$ for the holonomy in $E^a$ along the standard
generator of $\pi_1(S^1)$. Recall that $\det
A=\exp\bigl(\int_{S^1}\tr(a)d\theta\bigr)$.
Using the explicit formula in \cite[Theorem~1]{BFK91} we get:
\begin{align*}
\det(\Delta_{E^a,g,b,1})
&=\ii^{2k}\exp\Bigl(\frac\ii2\int_{S^1}\tr\bigl(\ii(b^{-1}a^tb-b^{-1}b'-a)\bigr)d\theta\Bigr)
\det\left(1-\left(\begin{smallmatrix}A^{-1}&*\\0&A^t\end{smallmatrix}\right)\right)
\\
&=\exp\Bigl(\frac12\int_{S^1}\tr(b^{-1}b')d\theta\Bigr)\det(A-1)^2\det A^{-1}
\\
&=\exp\Bigl(\frac12\int_{S^1}\bigl(\tr(b^{-1}b')-2\tr(a)\bigr)d\theta\Bigr)\det(A-1)^2
\end{align*}

Consider the Euler structure
$\e:=[-\frac\partial{\partial\theta},0]\in\Eul(S^1;\Z)$, and the coEuler
structure $\e^*:=[g,\frac12]\in\Eul^*(S^1;\C)$. Then $P(\e)=\e^*$, see
\eqref{E:665}. Assuming acyclicity, i.e.\ $1$ is not an eigen value of $A$,
we conclude:
\begin{equation}\label{E:ans1}
\tau^\an_{E^a,\e^*,[b]}=\det(A-1)^{-2}.
\end{equation}
Observe that this is independent of $[b]$, cf.\ Remark~\ref{R:bindep}.

Considering a Morse--Smale vector field $X$ with two critical points 
and the Euler structure $\e$ we obtain a 
Morse complex $C^*(X;E^a)$ isomorphic to
$$
\C^k\xrightarrow{A-1}\C^k
$$
equipped with the standard bilinear form. 
From Example~\ref{Ex:simtor} we obtain
$$
\tau_{E^a,\e}^\comb=\det(A-1)^{-2}
$$
which coincides with \eqref{E:ans1}. So we see that $\tau_{E^a,\e}^\comb=\tau^\an_{E^a,\e^*,[b]}$, i.e.\
$\mathcal S_{E^a,[b]}=1$,
whenever $E^a$ is acyclic. From Proposition~\ref{P:holoS} we conclude $\mathcal S_{E^a,[b]}=1$
for all, not necessarily acyclic, $E^a$. Thus Conjecture~\ref{C:main} holds for $M=S^1$.

\begin{remark}
Recall the canonic coEuler structure $\e^*_\can=[g,0]$ defined as the unique fixed
point of the involution on $\Eul^*(S^1;\C)$, see Section~\ref{S:eul}.
Note that $\e^*_\can$ is not integral. The computations above show that for 
the analytic torsion we have
$$
\tau^\an_{E^a,\e^*_\can,[b]}=s_{[b]}\det A\det(A-1)^{-2}
$$
where 
$$
s_{[b]}=\exp\Bigl(-\frac12\int_{S^1}\tr(b^{-1}b')\Bigr)\in\{\pm1\}
$$
does depend on $b$. Note that this sign $s_{[b]}$ appears, although we consider
the torsion as a bilinear form, i.e.\ we essentially consider the
square of what is traditionally called the torsion.

On odd dimensional manifolds one often
considers the analytic torsion without a correction term, i.e.\
one considers $\tau^\an_{E,\e^*_\can,[b]}$. Let us give two reasons
why this is not such a natural choice as it might seem.
First, the celebrated fact that the Ray--Singer torsion on odd dimensional manifolds
does only depend on the flat connection, is no longer true in the complex
setting as the appearance of the sign $s_{[b]}$ shows. Of course a different
definition of complex valued analytic torsion might circumvent this problem.
More serious is the second point. One would
expect that the analytic torsion as considered above is the square of a
rational function on the space of acyclic representations of the fundamental
group. As the computation for the circle shows, this cannot
be true for $\tau^\an_{E,\e^*_\can,[b]}$, simply because $\sqrt{\det A}$ cannot be
rational in $A\in\GL_k(\C)$. Any reasonable definition of complex valued analytic torsion
will have to face this problem.

If one is willing to consider
$\tau^\an_{E,\e^*,[b]}$ where $\e^*$ is an integral coEuler structure
both problems disappear, assuming $E$ admits a non-degenerate symmetric bilinear form and
Conjecture~\ref{C:main} is true.
Then $\tau^\an_{E,\e^*,[b]}$ is indeed independent of 
$[b]$, see Remark~\ref{R:bindep}, and the dependence on $\e^*$ is very
simple, see Remark~\ref{R:aneul}. More importantly, $\tau^\an_{E,\e^*,[b]}$
is the square of a rational function on the space of acyclic
representations of the fundamental group. This follows from the
fact that $\tau^\comb_{E,\e}$ with $P(\e)=\e^*$ is the square of such a rational function, see
\cite{BH03}.
\end{remark}

\section{Proof of the anomaly formula}\label{S:anom}

We continue to use the notation of Section~\ref{S:anator}.
The proof of Theorem~\ref{T:anom} is based on the following two results whose
proof we postpone till Section~\ref{S:anomprop}.

\begin{proposition}\label{P:strbD}
Suppose $\phi\in\Gamma(\eend(E))$. Then
$$
\LIM_{t\to0}\str\bigl(\phi e^{-t\Delta_{E,g,b}}\bigr)
=\int_M\tr(\phi)\ec(g).
$$
Here $\LIM$ denotes the renormalized limit, see \cite[Section~9.6]{BGV92},
which in this case is actually an ordinary limit.
\end{proposition}

\begin{proposition}\label{P:strgD}
Suppose $\xi\in\Gamma(\eend(TM))$ is symmetric with respect to $g$, 
and let $\Lambda^*\xi\in\eend(\Lambda^*T^*M)$ denote
its extension to a derivation on $\Lambda^*T^*M$. Then
$$
\LIM_{t\to0}\str\Bigl(\bigl(\Lambda^*\xi-\tfrac12\tr(\xi)\bigr)
e^{-t\Delta_{E,g,b}}\Bigr)
=\int_M\tr(b^{-1}\nabla^Eb)\wedge(\partial_2\cs)(g,g\xi).
$$
Again $\LIM$ denotes the renormalized limit,
which in this case is just the constant term of the asymptotic expansion
for $t\to0$. Moreover, we use the notation
$(\partial_2\cs)(g,g\xi):=\frac\partial{\partial t}|_0\cs(g,g+tg\xi)$.
\end{proposition}

Let us now give a proof of Theorem~\ref{T:anom}.
Suppose $g_u$ and $b_u$ depend smoothly on a real parameter
$u$. Let $\gamma$ be a simple closed curve around $0$ and assume w.l.o.g.\
that $u$ varies in an open set $U$ so that
the spectrum of $\Delta_u:=\Delta_{E,g_u,b_u}$ avoids the curve $\gamma$ for
all $u\in U$. Let $Q_u$ denote the spectral projection onto the eigen spaces
corresponding to eigen values in the exterior of $\gamma$, and $Q_{u,q}$ the part acting in
degree $q$. Let us write $\dot\Delta_u:=\frac\partial{\partial u}\Delta_u$, and
$\dot\Delta_{u,q}$ for the part acting in degree $q$.
From the variation formula for the determinant of generalized Laplacians,
see for instance \cite[Proposition~9.38]{BGV92}, we obtain
\begin{align}
\frac\partial{\partial u}
\log\prod_q\bigl({\det}^\gamma(\Delta_{u,q})\bigr)^{(-1)^qq}
&=\sum_q(-1)^qq\Bigl(\frac\partial{\partial u}\log{\det}^\gamma(\Delta_{u,q})\Bigr)
\notag\\
&=\sum_q(-1)^qq\Bigl(
\LIM_{t\to0}\tr\bigl(\dot\Delta_{u,q}(\Delta_{u,q})^{-1}Q_{u,q}e^{-t\Delta_{u,q}}\bigr)
\Bigr)
\notag\\&=\label{E:zeta}
\LIM_{t\to0}\str\bigl(N\dot\Delta_u\Delta^{-1}_uQ_ue^{-t\Delta_u}\bigr)
\end{align}
where $N$ denotes the grading operator which
acts by multiplication with $q$ on $\Omega^q(M;E)$.

Choose $u_0\in U$ and define $G_u\in\Gamma(\Aut(TM))$ by
$$
g_u(a,b)=g_{u_0}(G_ua,b)=g_{u_0}(a,G_ub)
$$ 
and similarly $B_u\in\Gamma(\Aut(E))$ by
$$
b_u(e,f)=b_{u_0}(B_ue,f)=b_{u_0}(e,B_uf).
$$
Let $\Lambda^*G_u^{-1}$ denote the natural extension of $G_u^{-1}$ to
$\Gamma(\Aut(\Lambda^*T^*M))$ and define 
$$
A_u=\det(G_u)^{1/2}\cdot\Lambda^*G_u^{-1}\otimes B_u\in\Gamma(\Aut(\Lambda^*T^*M\otimes E)).
$$
Then
\begin{equation}\label{E:betauu}
\beta_{g_u,b_u}(v,w)=\beta_{g_{u_0},b_{u_0}}(A_uv,w)=\beta_{g_{u_0},b_{u_0}}(v,A_uw),\qquad v,w\in\Omega(M;E).
\end{equation}
Abbreviating $d_u^\sharp:=d^\sharp_{E,g_u,b_u}$ we immediately get
$$
d^\sharp_u:=A_u^{-1}d_{u_0}^\sharp A_u.
$$ 
Writing $\dot A_u:=\frac\partial{\partial u}A_u$, $\dot
g_u:=\frac\partial{\partial u}g_u$ and $\dot b_u:=\frac\partial{\partial
u}b_u$ we have 
\begin{equation}\label{E:ada}
A_u^{-1}\dot A_u=
\bigl(-\Lambda^*(g_u^{-1}\dot g_u)+\tfrac12\tr(g_u^{-1}\dot g_u)\bigr)\otimes 1
+1\otimes(b_u^{-1}\dot b_u)
\in\Gamma(\eend(\Lambda^*T^*M\otimes E))
\end{equation}
where $\Lambda^*(g_u^{-1}\dot g_u)$ denotes the extension of 
$g_u^{-1}\dot g_u\in\Gamma(\eend(TM))$ to a derivation on $\Lambda^*T^*M$.

Let us write $d:=d_E$ and 
$\dot d_u^\sharp:=\frac\partial{\partial u}d_u^\sharp$.
Using the obvious relations $\dot\Delta_u=[d,\dot d_u^\sharp]$,
$[N,d]=d$, $[d,\Delta_u]=0$, $[d,Q_u]=0$, $\dot
d_u^\sharp=[d_u^\sharp,A_u^{-1}\dot A_u]$ and the fact that the super trace
vanishes on super commutators we get:
\begin{align*}
\str\bigl(N\dot\Delta_u\Delta^{-1}_uQ_ue^{-t\Delta_u}\bigr)
&
=\str\bigl(Nd\dot d^\sharp_u\Delta^{-1}_uQ_ue^{-t\Delta_u}\bigr)
+\str\bigl(N\dot d^\sharp_ud\Delta^{-1}_uQ_ue^{-t\Delta_u}\bigr)
\\&
=\str\bigl(d\dot d_u^\sharp\Delta_u^{-1}Q_ue^{-t\Delta_u}\bigr)
\\&
=\str\bigl(dd_u^\sharp A^{-1}_u\dot A_u\Delta_u^{-1}Q_ue^{-t\Delta_u}\bigr)
\\&\qquad-\str\bigl(dA_u^{-1}\dot A_ud_u^\sharp\Delta_u^{-1}Q_ue^{-t\Delta_u}\bigr)
\\&=\str\bigl(A_u^{-1}\dot A_u(dd_u^\sharp+d_u^\sharp d)\Delta_u^{-1}Q_ue^{-t\Delta_u}\bigr)
\\&=\str\bigl(A_u^{-1}\dot A_uQ_ue^{-t\Delta_u}\bigr)
\end{align*}
Together with \eqref{E:zeta} this gives
\begin{equation}\label{E:offK}
\frac\partial{\partial u}
\log\prod_q\bigl({\det}^\gamma(\Delta_{u,q})\bigr)^{(-1)^qq}
=\LIM_{t\to0}\str\bigl(A_u^{-1}\dot A_uQ_ue^{-t\Delta_u}\bigr)
\end{equation}

Let us write $\Omega^*_u:=\Omega^*_{E,g_u,b_u}(M;E)(\gamma)$. Note that this is a family
of finite dimensional complexes smoothly parametrized by $u\in U$. Let $P_u=1-Q_u$ denote the spectral
projection of $\Delta_u$ onto $\Omega^*_u$.
Note that since $\str P_uP_u=\const$ we have $\str P_u\dot P_u=0$.
For sufficiently small $w-u$ the restriction of the spectral
projection $P_w|_{\Omega^*_u}:\Omega^*_u\to\Omega^*_w$ is an isomorphism of
complexes. We get a commutative diagram of determinant lines:
$$
\xymatrix{
\det \Omega^*_u \ar[r] \ar[d]_-{\det(P_w|_{\Omega^*_u})} 
& \det H(\Omega^*_u) \ar[d]^-{\det H(P_w|_{\Omega^*_u})} \ar[r]
& \det H^*(M;E) \ar[d]^-{\det H(P_w)=1}
\\
\det \Omega_w^* \ar[r] & \det H(\Omega_w^*) \ar[r] & \det H^*(M;E)
}
$$
Writing $\beta_u:=\beta_{E,g_u,b_u}$ and $\tau^\an_u(\gamma):=\tau^\an_{E,g_u,b_u}(\gamma)$,
we obtain, for sufficiently small $w-u$,
\begin{equation}\label{E:L1}
\frac{\tau^\an_w(\gamma)}{\tau^\an_u(\gamma)}
=\sdet\Bigl((\beta_u|_{\Omega_u^*})^{-1}(P_w|_{\Omega_u^*})^*\beta_w\Bigr).
\end{equation}
Here the two non-degenerate bilinear forms $\beta_u|_{\Omega_u^*}$ and 
$(P_w|_{\Omega_u^*})^*\beta_w$ on $\Omega^*_u$ are considered as 
isomorphisms from $\Omega_u^*$ to its dual, hence
$(\beta_u|_{\Omega_u^*})^{-1}(P_w|_{\Omega_u^*})^*\beta_w$ is an automorphism of $\Omega_u^*$.
Using \eqref{E:betauu} we find
$$
(\beta_u|_{\Omega_u^*})^{-1}(P_w|_{\Omega_u^*})^*\beta_w=P_uA_u^{-1}A_wP_w|_{\Omega_u^*}.
$$
Using \eqref{E:L1} we thus obtain
$$
\frac{\tau^\an_w(\gamma)}{\tau^\an_u(\gamma)}
=\sdet\bigl(P_uA_u^{-1}A_wP_w|_{\Omega^*_u}\bigr).
$$
In view of $\str(P_u\dot P_u)=0$ we get
$$
\frac\partial{\partial w}\Big|_u\Biggl(\frac{\tau^\an_w(\gamma)}{\tau^\an_u(\gamma)}\Biggr)
=\str\bigl(P_uA_u^{-1}\dot A_uP_u+P_uA_u^{-1}A_u\dot P_u\bigr)
=\str\bigl(A_u^{-1}\dot A_uP_u\bigr).
$$
Combining this with \eqref{E:offK} and Proposition~\ref{P:A} we obtain
\begin{equation}\label{E:main}
\frac\partial{\partial w}\Big|_u
\left(\frac{\tau^\an_w(0)\cdot\prod_q({\det}'\Delta_{w,q})^{(-1)^qq}}
{\tau^\an_u(0)\cdot\prod_q({\det}'\Delta_{u,q})^{(-1)^qq}}\right)
=\LIM_{t\to0}\str\bigl(A_u^{-1}\dot A_ue^{-t\Delta_u}\bigr).
\end{equation}

Applying Proposition~\ref{P:strbD} to $\phi=b_u^{-1}\dot b_u$ we obtain 
$$
\LIM_{t\to0}\str\bigl(b_u^{-1}\dot b_ue^{-t\Delta_u}\bigr)
=\int_M\tr(b_u^{-1}\dot b_u)\ec(g_u).
$$
Using Proposition~\ref{P:strgD} with $\xi=g_u^{-1}\dot g_u$ we get
\begin{multline*}
\LIM_{t\to0}\str\Bigl(\bigl(\Lambda^*(g_u^{-1}\dot g_u)
-\tfrac12\tr(g_u^{-1}\dot g_u)\bigr)e^{-t\Delta_u}\Bigr)
\\=\int_M\tr(b^{-1}_u\nabla^Eb_u)\wedge(\partial_2\cs)(g_u,\dot g_u).
\end{multline*}
Using \eqref{E:ada} we conclude
\begin{multline}\label{E:varD}
\LIM_{t\to0}\str\bigl(A_u^{-1}\dot A_ue^{-t\Delta_u}\bigr)
\\=\int_M\tr(b_u^{-1}\dot b_u)\ec(g_u)
-\int_M\tr(b^{-1}_u\nabla^Eb_u)\wedge(\partial_2\cs)(g_u,\dot g_u).
\end{multline}

Let us finally turn to the correction term. If
$[g_u,\alpha_u]\in\Eul^*(M;\mathbb C)$ represent the
same coEuler structure then $\alpha_w-\alpha_u=\cs(g_u,g_w)$ and thus
$$
\frac\partial{\partial u}\alpha_u=
\frac\partial{\partial w}\Big|_u\cs(g_u,g_w)
=(\partial_2\cs)(g_u,\dot g_u).
$$
Moreover, we have
\begin{align*}
\frac\partial{\partial u}
\tr(b_u^{-1}\nabla^Eb_u)
&=\tr\bigl(-b_u^{-1}\dot b_ub_u^{-1}\nabla^Eb_u\bigr)+\tr\bigl(b_u^{-1}\nabla^E\dot b_u\bigr)
\\&
=\tr\bigl(-b_u^{-1}(\nabla^Eb_u)b_u^{-1}\dot b_u\bigr)+\tr\bigl(b_u^{-1}\nabla^E\dot b_u\bigr)
\\&
=\tr\bigl((\nabla^Eb_u^{-1})\dot b_u\bigr)+\tr\bigl(b_u^{-1}\nabla^E\dot b_u\bigr)
\\&
=\tr\bigl(\nabla^E(b_u^{-1}\dot b_u)\bigr)
\\&
=d\tr(b_u^{-1}\dot b_u).
\end{align*}
Using $-2\omega_{E,b_u}=\tr(b^{-1}_u\nabla^Eb_u)$, $d\alpha_u=\ec(g_u)$ and Stokes'
theorem we get
\begin{multline}\label{E:corr}
\frac\partial{\partial u}\int_M-2\omega_{E,b_u}\wedge\alpha_u
=\int_Md\tr(b_u^{-1}\dot b_u)\wedge\alpha_u+\int_M\tr(b_u^{-1}\nabla^Eb_u)\wedge(\partial_2\cs)(g_u,\dot g_u)
\\=
-\int_M\tr(b_u^{-1}\dot b_u)\ec(g_u)
+\int_M\tr(b_u^{-1}\nabla^Eb_u)\wedge(\partial_2\cs)(g_u,\dot g_u).
\end{multline}
Combining \eqref{E:main}, \eqref{E:varD} and \eqref{E:corr} 
we obtain
$$
\frac\partial{\partial w}\Bigm|_u
\frac{\tau^\an_{E,g_w,b_w,\alpha_w}}{\tau^\an_{E,g_u,b_u,\alpha_u}}
=0.
$$
This completes the proof of Theorem~\ref{T:anom}.

\section{Asymptotic expansion of the heat kernel}\label{S:ass_ex}

In this section we will consider Dirac operators associated to a 
class of Clifford super connections. The main result Theorem~\ref{T:ass_exp}
below computes the leading and subleading terms of the asymptotic expansion
of the corresponding heat kernels. In Section~\ref{S:anomprop} we will apply these
results to the Laplacians introduced in Section~\ref{S:anator} which are 
squares of such Dirac operators. We refer to
\cite{BGV92} for background on the Clifford super connection formalism.

Let $(M,g)$ be a closed Riemannian manifold of dimension $n$.
Let $\Cl=\Cl(T^*M,g)$ denote the corresponding Clifford bundle.
Recall that $\Cl=\Cl^+\oplus\Cl^-$ is a bundle of $\mathbb Z_2$-graded filtered algebras, and let us
write $\Cl_k$ for the subbundle of filtration degree $k$.
Recall that we have the symbol map
$$
\sigma:\Cl\to\Lambda^*T^*M,\qquad\sigma(a):=c(a)\cdot 1
$$
where $c$ denotes the usual Clifford action on $\Lambda^*T^*M$. Explicitly,
for $a\in T^*_xM\subseteq\Cl_x$ and $\alpha\in\Lambda^*T^*_xM$ we have
$c(a)\cdot\alpha=a\wedge\alpha-i_{\sharp a}\alpha$, where $\sharp a=g^{-1}a\in T_xM$
and $i_{\sharp a}$ denotes contraction with $\sharp a$. Here the metric is
considered as an isomorphism $g:TM\to T^*M$ and $g^{-1}$ denotes its
inverse. Recall that $\sigma$ is an isomorphism of filtered 
$\mathbb Z_2$-graded vector bundles inducing an
isomorphism on the associated graded bundles of algebras.

Let $\mathcal E=\mathcal E^+\oplus\mathcal E^-$ 
be a $\mathbb Z_2$-graded complex Clifford module over $M$.
The forms with values in $\mathcal E$ inherit a $\mathbb Z_2$-grading
which will be denoted by:
$$
\Omega(M;\mathcal E)=\Omega(M;\mathcal E)^+\oplus\Omega(M;\mathcal E)^-
$$
We have $\Omega(M;\mathcal E)^+=\Omega^\even(M;\mathcal E^+)\oplus
\Omega^\odd(M;\mathcal E^-)$ and similarly for $\Omega(M;\mathcal E)^-$.
Let us write $\eend_{\Cl}(\mathcal E)$ for the bundle
of algebras of endomorphisms of $\mathcal E$ which (super) commute with the Clifford action,
and let us indicate its $\mathbb Z_2$-grading by:
$$
\eend_{\Cl}(\mathcal E)
=\eend^+_{\Cl}(\mathcal E)\oplus\eend^-_{\Cl}(\mathcal E)
$$
Recall that we have a canonic isomorphism of bundles of $\mathbb Z_2$-graded algebras
\begin{equation}\label{E:endE}
\eend(\mathcal E)=\Cl\otimes\eend_{\Cl}(\mathcal E).
\end{equation}

Suppose $\mathbb A:\Omega(M;\mathcal E)^\pm\to\Omega(M;\mathcal E)^\mp$
is a Clifford super connection, see \cite[Definition~3.39]{BGV92}.
Recall that with respect to \eqref{E:endE} 
its curvature $\mathbb A^2\in\Omega(M;\eend(\mathcal E))^+$
decomposes as
\begin{equation}\label{E:Q4}
\mathbb A^2=R^{\Cl}\otimes1+1\otimes F^{\mathcal E/S}_{\mathbb A}
\end{equation}
where $R^{\Cl}\in\Omega^2(M;\Cl^2)$ with $\Cl^2:=\sigma^{-1}(\Lambda^2T^*M)\subseteq\Cl^+$ 
is a variant of the Riemannian curvature
\begin{equation}\label{E:Q6}
R^{\Cl}(X,Y)=\tfrac14\sum_{i,j}g(R_{X,Y}e_i,e_j)c^ic^j
\end{equation}
and $F^{\mathcal E/S}_{\mathbb A}\in\Omega(M;\eend_{\Cl}(\mathcal E))^+$ is called
the twisting curvature, see \cite[Proposition~3.43]{BGV92}. Here $e_i$ is a local
orthonormal frame of $TM$, $e^i:=ge_i$ denotes its dual local coframe
and $c^i=c(e^i)$ denotes Clifford multiplication with $e^i$.

Recall that the Dirac operator $D_{\mathbb A}$ associated to the Clifford
super connection $\mathbb A$ is given by the composition
$$
\Gamma(\mathcal E)\xrightarrow{\mathbb A}
\Omega(M;\mathcal E)=\Gamma(\Lambda^*T^*M\otimes\mathcal E)
\xrightarrow{\sigma^{-1}\otimes 1_{\mathcal E}}
\Gamma(\Cl\otimes\mathcal E)\xrightarrow{c}
\Gamma(\mathcal E)
$$
where $c:\Cl\otimes\mathcal E\to\mathcal E$ denotes Clifford
multiplication.

We will from now on restrict to very special Clifford super connections on $\mathcal E$
which are of the form
$$
\mathbb A=\nabla+A
$$
where $\nabla:\Omega^*(M;\mathcal E^\pm)\to\Omega^{*+1}(M;\mathcal E^\pm)$ 
is a Clifford connection on $\mathcal E$, and
$A\in\Omega^0(M;\eend_{\Cl}^-(\mathcal E))$.
For the associated Dirac operator acting on $\Gamma(\mathcal E)$ we have
$$
D_{\mathbb A}=D_\nabla+A.
$$
Consider the induced connection
$\nabla:\Omega^*(M;\eend^\pm(\mathcal E))\to\Omega^{*+1}(M;\eend^\pm(\mathcal E))$.
Since $\nabla$ is a Clifford connection this induced connection preserves the subbundle
$\eend_{\Cl}(\mathcal E)$. Moreover, we have
$[D_\nabla,A]=c(\nabla A)$ and thus
\begin{equation}\label{E:D2}
D_{\mathbb A}^2=D_\nabla^2+c(\nabla A)+A^2.
\end{equation}
Here $\nabla A\in\Omega^1(M;\eend_{\Cl}^-(\mathcal E))$,
$A^2\in\Omega^0(M;\eend_{\Cl}^+(\mathcal E))$, and 
the Clifford action $c(B)$ of $B\in\Omega(M;\eend(\mathcal E))$ on 
$\Gamma(\mathcal E)$ is given by the composition:
$$
\Gamma(\mathcal E)\xrightarrow{B}
\Omega(M;\mathcal E)=\Gamma(\Lambda^*T^*M\otimes\mathcal E)
\xrightarrow{\sigma^{-1}\otimes 1_{\mathcal E}}\Gamma(\Cl\otimes\mathcal E)
\xrightarrow{c}\Gamma(\mathcal E)
$$
Note that for $B\in\Omega^0(M;\eend(\mathcal E))$ the Clifford action 
coincides with the usual action $c(B)=B$.

\begin{theorem}\label{T:ass_exp}
Let $\mathcal E$ be a $\mathbb Z_2$-graded complex Clifford bundle
over a closed Riemannian manifold $(M,g)$ of dimension $n$. Suppose $\nabla$ is a
Clifford connection on $\mathcal E$ and
$A\in\Omega^0(M;\eend_{\Cl}^-(\mathcal E))$. Consider the Clifford super
connection $\mathbb A=\nabla+A$ and the associated
Dirac operator $D_{\mathbb A}$ acting on $\Gamma(\mathcal E)$.
Let $\Omega_g\in\Omega^n(M;\mathcal O_M^\C)$ denote the volume density
associated with the Riemannian metric $g$.
Let $k_t\in\Gamma(\eend(\mathcal E))$ so that $k_t\Omega_g$ is
the restriction of the kernel of 
$e^{-tD^2_{\mathbb A}}$ to the diagonal in $M\times M$. Consider its asymptotic expansion
\begin{equation}\label{E:ass_exp}
k_t\sim(4\pi t)^{-n/2}\sum_{i\geq0}t^i\tilde k_i
\qquad\text{as $t\to0$}
\end{equation}
with $\tilde k_i\in\Gamma(\eend(\mathcal E))$, see \cite[Theorem~2.30]{BGV92}. Then
\begin{equation}\label{E:thmi}
\tilde k_i\in\Gamma(\Cl_{2i}\otimes\eend_{\Cl}(\mathcal E))
\subseteq\Gamma(\Cl\otimes\eend_{\Cl}(\mathcal E))
=\Gamma(\eend(\mathcal E)).
\end{equation}
Moreover, with the help of the symbol map
$$
\sigma:\Gamma(\eend(\mathcal E))=\Gamma(\Cl\otimes\eend_{\Cl}(\mathcal E))
\xrightarrow{\sigma\otimes1}
\Omega^*(M;\eend_{\Cl}(\mathcal E))
$$
and writing $\alpha_{[j]}$ for the $j$-form piece of $\alpha$ we have
\begin{equation}\label{E:thmii}
\sum_{i\geq0}\sigma(\tilde k_i)_{[2i]}=
\hat A_g\wedge\exp(-F^{\mathcal E/S}_\nabla).
\end{equation}
Here $\hat A_g\in\Omega^{4*}(M;\mathbb R)$ denotes the $\hat A$-genus
$$
\hat A_g={\det}^{1/2}\biggl(\frac{R/2}{\sinh(R/2)}\biggr)
$$
and $R\in\Omega^2(M;\eend(TM))$ the Riemannian curvature.
Moreover, we have
\begin{equation}\label{E:thmiii}
\sum_{i\geq0}\sigma(\tilde k_i)_{[2i-1]}=
-\nabla\Biggl(\hat A_g\wedge
\biggl(\frac{e^{\ad(-F^{\mathcal E/S}_\nabla)}-1}{\ad(-F^{\mathcal E/S}_\nabla)}A\biggr)
\wedge\exp\bigl(-F^{\mathcal E/S}_\nabla\bigr)\Biggr)
\end{equation}
where $\ad(F^{\mathcal E/S}_\nabla):\Omega^*(M;\eend_{\Cl}^\pm(\mathcal E))
\to\Omega^{*+2}(M;\eend_{\Cl}^\pm(\mathcal E))$, is given by
$$
\ad(F^{\mathcal E/S}_\nabla)\phi
:=F^{\mathcal E/S}_\nabla\wedge\phi-\phi\wedge F^{\mathcal E/S}_\nabla.
$$
\end{theorem}

\begin{remark}
Note that \eqref{E:thmi} and \eqref{E:thmii} tell that on this level
the asymptotic expansions for
$e^{-tD_{\mathbb A}^2}$ and $e^{-tD_\nabla^2}$ are the same.
\end{remark}

\begin{proof}
The proof below parallels the one of Theorem~4.1 in \cite{BGV92} where the
case $A=0$ is treated. It too is based on Getzler's scaling techniques, see \cite{G86}.
In order to prove Theorem~\ref{T:ass_exp} we need to compute one more term
in the asymptotic expansion of the rescaled operator.

The calculation is local.
Let $x_0\in M$. Use normal coordinates, i.e.\ the
exponential mapping of $g$, to identify a convex neighborhood $U$
of $0\in T_{x_0}M$ with a neighborhood of $x_0$.
Choose an orthonormal basis $\{\partial_i\}$ of $T_{x_0}M$ and linear
coordinates $\mathbf x=(x^1,\dotsc,x^n)$ on $T_{x_0}M$ such that
$\{dx^i\}$ is dual to $\{\partial_i\}$.
Let $\mathcal R:=\sum_ix^i\partial_i$ denote the radial vector field.
Note that every affinely parametrized line through the origin in $T_{x_0}M$ is a geodesic.
Let $\{e_i\}$ denote the local orthonormal frame of $TM$
obtained from $\{\partial_i\}$ by parallel transport along $\mathcal R$,
i.e.\ $\nabla^g_{\mathcal R}e_i=0$ and $e_i(x_0)=\partial_i$.
Let $\{e^i\}$ denote the dual local coframe.

Trivialize $\mathcal E$
with the help of radial parallel transport by $\nabla$. Use this trivialization to identify
$\Gamma(\mathcal E|_U)$ with $C^\infty(U,\mathcal E_0)$, where 
$\mathcal E_0:=\mathcal E_{x_0}$. Let $\omega\in\Omega^1(U;\eend(\mathcal E_0))$
denote the connection one form of this trivialization, i.e.\ $\nabla_{\partial_i}
=\partial_i+\omega(\partial_i)$. For the curvature $F$ of $\nabla$ we then have
$F=d\omega+\omega\wedge\omega\in\Omega^2(M;\eend^+(\mathcal E_0))$.
By the choice of trivialization of $\mathcal E|_U$ we have $i_{\mathcal R}\omega=0$
and thus $i_{\mathcal R}F=i_{\mathcal R}(d\omega+\omega\wedge\omega)
=i_{\mathcal R}d\omega$. Contracting this with $\partial_i$ and using 
$[\partial_i,\mathcal R]=\partial_i$ we obtain
\begin{equation}\label{E:abc}
-F(\partial_i,\mathcal R)
=F(\mathcal R,\partial_i)
=(d\omega)(\mathcal R,\partial_i)
=(L_{\mathcal R}+1)(\omega(\partial_i))
\end{equation}
where $L_{\mathcal R}$ denotes Lie derivative with respect to the vector field $\mathcal R$.
Let $\omega(\partial_i)\sim\sum_\alpha\frac{\partial_\alpha\omega(\partial_i)_{x_0}}{\alpha!}x^\alpha$
denote the Taylor expansion of $\omega(\partial_i)$ at $x_0$, written with the help of multi index 
notation. Using $L_{\mathcal R}x^\alpha=|\alpha|x^\alpha$ we obtain the following Taylor expansion
$(L_{\mathcal R}+1)(\omega(\partial_i))\sim\sum_\alpha\bigl(|\alpha|+1\bigr)
\frac{\partial_\alpha\omega(\partial_i)_{x_0}}{\alpha!}x^\alpha$. If 
$F(\partial_i,\partial_j)\sim\sum_\alpha\frac{\partial_\alpha F(\partial_i,\partial_j)_{x_0}}{\alpha!}x^\alpha$
denotes the Taylor expansion of $F(\partial_i,\partial_j)$ at $x_0$ then we obtain the Taylor expansion
$F(\partial_i,\mathcal R)\sim\sum_{j,\alpha}\frac{\partial_\alpha F(\partial_i,\partial_j)_{x_0}}{\alpha!}x^jx^\alpha$.
Comparing the Taylor expansions of both sides of \eqref{E:abc} we obtain the Taylor expansion,
cf.\ Proposition~1.18 in \cite{BGV92},
$$
\nabla_{\partial_i}-\partial_i
=\omega(\partial_i)
\sim-\sum_{j,\alpha}\frac{\partial_\alpha F(\partial_i,\partial_j)_{x_0}}{(|\alpha|+2)\alpha!}
x^jx^\alpha.
$$
For the first few terms this gives:
\begin{equation}\label{E:bcd}
\nabla_{\partial_i}=\partial_i
-\tfrac12\sum_jF(\partial_i,\partial_j)_{x_0}x^j
-\tfrac13\sum_{j,k}\partial_kF(\partial_i,\partial_j)_{x_0}x^jx^k
+O(|\mathbf x|^3)
\end{equation}

Let $c^i:=c(e^i)\in \Gamma(\mathcal E|_U)=C^\infty(U,\eend(\mathcal E_0))$ denote Clifford multiplication with 
$e^i$. Since $\nabla^g_{\mathcal R}e^i=0$ and since $\nabla$ is a Clifford connection we have
$\nabla_{\mathcal R}c^i=c(\nabla^g_{\mathcal R}e^i)=0$. So we see that $c^i$ is actually
a constant in $\eend(\mathcal E_0)$, cf.\ \cite[Lemma~4.14]{BGV92}.
Particularly, our trivialization of $\mathcal E|_U$ identifies
$\Gamma(\eend_{\Cl}(\mathcal E|_U))$ with $C^\infty(U,\eend_{\Cl}(\mathcal E_0))$.
Recall that 
$$
F(\partial_i,\partial_j)=\tfrac14\sum_{l,m}g(R_{\partial_i,\partial_j}e_l,e_m)c^lc^m
+F^{\mathcal E/S}_\nabla(\partial_i,\partial_j)
$$
with $F^{\mathcal E/S}_\nabla\in\Omega^2(U;\eend^+_{\Cl}(\mathcal E_0))$.
From \eqref{E:bcd} we thus obtain, cf.\ \cite[Lemma~4.15]{BGV92},
\begin{align}\label{E:565}
\notag
\nabla_{\partial_i}
=&\ \partial_i-\tfrac18\sum_{j,l,m}g(R_{\partial_i,\partial_j}e_l,e_m)_{x_0}x^jc^lc^m
\\&\notag
-\tfrac1{12}\sum_{j,k,l,m}\partial_kg(R_{\partial_i,\partial_j}e_l,e_m)_{x_0}
x^jx^kc^lc^m
\\&
+\sum_{l,m}u_{ilm}(\mathbf x)c^lc^m
+v_i(\mathbf x)
\end{align}
with
$u_{ilm}(\mathbf x)=O(|\mathbf x|^3)\in C^\infty(U)$ 
and
$v_i(\mathbf x)=O(|\mathbf x|)\in C^\infty(U,\eend_{\Cl}(\mathcal E_0))$.

Let $\Delta$ denote the connection Laplacian given by the composition
$$
\Gamma(\mathcal E)\xrightarrow{\nabla}\Gamma(T^*M\otimes\mathcal E)
\xrightarrow{\nabla^g\otimes 1+1\otimes\nabla}
\Gamma(T^*M\otimes T^*M\otimes\mathcal E)
\xrightarrow{-\tr_g}\Gamma(\mathcal E)
$$
Let $r$ denote the scalar curvature of $g$ and recall Lichnerowicz' formula
\cite[Theorem~3.52]{BGV92}
$$
D^2_\nabla=\Delta+c(F^{\mathcal E/S}_\nabla)+\tfrac r4.
$$
Recall our Clifford super connection $\mathbb A=\nabla+A$ with
$A\in\Omega^0(M;\eend_{\Cl}^-(\mathcal E))$. Since
$D^2_{\mathbb A}=D^2_\nabla+c(\nabla A)+A^2$ we obtain
\begin{equation}\label{E:extL}
D^2_{\mathbb A}=\Delta+c\bigl(F^{\mathcal E/S}_\nabla+\nabla A\bigr)+A^2+\tfrac r4.
\end{equation}

Use the symbol map to identify
$\eend(\mathcal E_0)=\Lambda^*T^*_{x_0}M\otimes\eend_{\Cl}(\mathcal E_0)$.
For $0<u\leq1$ and $\alpha\in C^\infty\bigl(\mathbb R^+\times
U,\Lambda^*T_{x_0}^*M\otimes\eend_{\Cl}(\mathcal E_0)\bigr)$ define
Getzler's rescaling
$$
(\delta_u\alpha)(t,\mathbf x)
:=\sum_iu^{-i/2}\alpha(ut,u^{1/2}\mathbf x)_{[i]}.
$$
Consider the kernel 
$p\in C^\infty\bigl(\mathbb R^+\times U,\Lambda^*T^*_{x_0}M\otimes\eend_{\Cl}(\mathcal E_0)\bigr)$
of $e^{-tD^2_{\mathbb A}}$, $p(t,\mathbf x)=p_t(\mathbf x,x_0)$. Note that
$p(t,0)=k_t(x_0)$.
Define the rescaled kernel $r_u:=u^{n/2}\delta_up$ and 
the rescaled operator $L_u:=u\delta_uD^2_{\mathbb A}\delta_u^{-1}$.
Since $(\partial_t+D_{\mathbb A}^2)p=0$ and 
$\delta_u\partial_t\delta_u^{-1}=u^{-1}\partial_t$ we have
\begin{equation}\label{E:sc_heat}
(\partial_t+L_u)r_u=0.
\end{equation}
Note that setting $t=1$ and $x=0$ and using \eqref{E:ass_exp} 
we get an asymptotic expansion
\begin{equation}\label{E:cde}
r_u(1,0)\sim(4\pi)^{-n/2}\sum_{j\geq-n}u^{j/2}\sum_{i\geq0}\sigma\bigl(\tilde k_i(x_0)\bigr)_{[2i-j]}
\qquad\text{as $u\to0$.}
\end{equation}
The claim \eqref{E:thmi} just states that the terms for $-n\leq j<0$ vanish, i.e.\
there are no Laurent terms in \eqref{E:cde}. Statements \eqref{E:thmii} and \eqref{E:thmiii}
are explicit expressions for the term $j=0$ and $j=1$ in \eqref{E:cde}.

Let us compute the first terms in the asymptotic expansion of $L_u$ in powers of
$u^{1/2}$. Let us write $\varepsilon^j$ for the exterior multiplication with
$e^j$, and $\iota^j$ for the contraction with $e_j$. Note that
$$
\delta_u\varepsilon^j\delta_u^{-1}=u^{-1/2}\varepsilon^j,
\qquad
\delta_u\iota^j\delta_u^{-1}=u^{1/2}\iota^j,
\qquad
\delta_u\partial_i\delta_u^{-1}=u^{-1/2}\partial_i
$$
and recall that $c^j=\varepsilon^j-\iota^j$. 
Let us look at the simplest part first. Clearly,
\begin{equation}\label{E:sc1}
u\delta_u(A^2+\tfrac r4)\delta_u^{-1}=O(u)
\qquad\text{as $u\to0$.}
\end{equation}
Next we have $\nabla A=\sum_i(\nabla_{e_i}A)e^i$, hence
$c(\nabla A)=\sum_i(\nabla_{e_i}A)c^i$ and therefore
\begin{equation}\label{E:sc2}
u\delta_uc(\nabla A)\delta_u^{-1}=u^{1/2}\mathsf A'+O(u^{3/2})
\qquad\text{as $u\to0$,}
\end{equation}
where $\mathsf A':=\sum_i(\nabla_{\partial_i}A)_{x_0}\varepsilon^i$.
Moreover, 
$F^{\mathcal E/S}_\nabla=\tfrac12\sum_{i,j}F^{\mathcal 
E/S}_\nabla(e_i,e_j)e^i\wedge e^j$, hence
$c(F^{\mathcal E/S}_\nabla)=\tfrac12\sum_{i,j}F^{\mathcal 
E/S}_\nabla(e_i,e_j)c^ic^j$, and thus
\begin{equation}\label{E:sc3}
u\delta_uc(F^{\mathcal E/S}_\nabla)\delta_u^{-1}=\mathsf F+O(u)
\qquad\text{as $u\to0$,}
\end{equation}
where $\mathsf F:=\tfrac12\sum_{i,j}F^{\mathcal
E/S}_\nabla(\partial_i,\partial_j)_{x_0}\varepsilon^i\varepsilon^j$.
From \eqref{E:565} we easily get
$$
u^{1/2}\delta_u\nabla_{\partial_i}\delta_u^{-1}
=\partial_i-\tfrac14\sum_j\mathsf R_{ij}x^j
+u^{1/2}R'_i+O(u)
\qquad\text{as $u\to0$,}
$$
where
$\mathsf
R_{ij}:=\tfrac12\sum_{l,m}g(R_{\partial_i,\partial_j}e_l,e_m)_{x_0}\varepsilon^l\varepsilon^m$
and $\mathsf R'_i$ is an operator which acts on
$C^\infty\bigl(U,\Lambda^{\even/\odd}T^*_{x_0}M\otimes\eend_{\Cl}(\mathcal E_0)\bigr)$
in a way which preserves the parity of the form degree.
Using the formula
$\Delta=-\sum_i\bigl((\nabla_{e_i})^2-\nabla_{\nabla^g_{e_i}e_i}\bigr)$
and the fact that $\nabla^g_{e_i}e_i$ vanishes at $x_0$ we obtain
\begin{equation}\label{E:sc4}
u\delta_u\Delta\delta_u^{-1}=
-\sum_i\Bigl(\partial_i-\tfrac14\sum_j\mathsf R_{ij}x^j\Bigr)^2
+u^{1/2}K^\even+O(u)
\qquad\text{as $u\to0$,}
\end{equation}
where $K^\even$ acts on
$C^\infty\bigl(U,\Lambda^{\even/\odd}T^*_{x_0}M\otimes\eend_{\Cl}(\mathcal
E_0)\bigr)$ in a parity preserving way.
Let us write
$$
\mathsf K:=-\sum_i\Bigl(\partial_i-\tfrac14\sum_j\mathsf
R_{ij}x^j\Bigr)^2+\mathsf F.
$$
Then \eqref{E:sc1}, \eqref{E:sc2}, \eqref{E:sc3}, \eqref{E:sc4} together
with \eqref{E:extL} finally give
\begin{equation}\label{E:sc_asy}
L_u=\mathsf K+u^{1/2}\bigl(\mathsf A'+K^\even\bigr)+O(u)
\qquad\text{as $u\to0$.}
\end{equation}

Recall, see \cite[Lemma~4.19]{BGV92}, that there exist 
$\Lambda^*T^*M\otimes\eend_{\Cl}(\mathcal E_0)$ valued polynomials $\tilde r_i$
on $\R\times U$ so that we have an asymptotic expansion
\begin{equation}\label{E:sc_ker}
r_u(t,\mathbf x)\sim q_t(\mathbf x)\sum_{j\geq-n}u^{j/2}\tilde r_j(t,\mathbf x)
\qquad\text{as $u\to0$,}
\end{equation}
where $q_t(\mathbf x)=(4\pi t)^{-n/2}e^{-|\mathbf x|^2/4t}$.
Moreover, the initial condition for the heat kernel translates to
\begin{equation}\label{E:init}
\tilde r_j(0,0)=\delta_{j,0}.
\end{equation}
Setting $t=1$, $\mathbf x=0$ in \eqref{E:sc_ker} we get
\begin{equation}\label{E:X}
r_u(1,0)\sim(4\pi)^{-n/2}\sum_{j\geq-n}u^{j/2}\tilde r_j(1,0)
\qquad\text{as $u\to0$.}
\end{equation}
Comparing this with \eqref{E:cde} we obtain
\begin{equation}\label{E:X2}
\tilde r_j(1,0)=\sum_{i\geq0}\sigma(\tilde k_i)_{[2i-j]}(x_0).
\end{equation}

Expanding the equation $(\partial_t+L_u)r_u=0$ in a power series in
$u^{1/2}$ with the help of \eqref{E:sc_ker} and \eqref{E:sc_asy}
the leading term $q\tilde r_l$ satisfies $(\partial_t+\mathsf K)(q\tilde r_l)=0$.
Because of the initial condition \eqref{E:init} and the uniqueness of formal
solutions \cite[Theorem~4.13]{BGV92} we must have $l\geq0$ and thus
$\tilde r_j=0$ for $j<0$. In view of \eqref{E:X2} this proves \eqref{E:thmi}.

So $q\tilde r_0$ satisfies $(\partial_t+\mathsf K)(q\tilde r_0)=0$
with initial condition $\tilde r_0(0,0)=1$, see \eqref{E:init}. Mehler's formula
\cite[Theorem~4.13]{BGV92} provides an explicit solution:
\begin{multline*}
q_t(\mathbf x)\tilde r_0(t,\mathbf x)
\\=(4\pi t)^{-n/2}{\det}^{1/2}
\biggl(\frac{t\mathsf R/2}{\sinh(t\mathsf R/2)}\biggr)\wedge
\exp\Bigl(-\tfrac1{4t}\bigl\langle\mathbf x\bigm|\tfrac{t\mathsf
R}2\coth(\tfrac{t\mathsf R}2)\bigm|\mathbf x\bigr\rangle\Bigr)\wedge\exp(-t\mathsf F)
\end{multline*}
Setting $t=1$, $\mathbf x=0$ we obtain
\begin{equation}\label{E:fgh}
\tilde r_0(1,0)
\\={\det}^{1/2}\biggl(\frac{\mathsf R/2}{\sinh(\mathsf R/2)}\biggr)\wedge\exp(-\mathsf F).
\end{equation}
In view of \eqref{E:X2} we thus have established \eqref{E:thmii}.

The term $q\tilde r_1$ satisfies $(\partial_t+\mathsf K)(q\tilde r_1)
=-(\mathsf A'+K^\even)(\tilde qr_0)$. Let us write 
$$
\tilde r_1(t,\mathbf x)
=\tilde r_1^\even(t,\mathbf x)+\tilde r_1^\odd(t,\mathbf x)
$$
with $\tilde r_1^\even(t,\mathbf x)\in\Lambda^\even T^*_{x_0}M\otimes\eend_{\Cl}(\mathcal E_0)$
and $\tilde r_1^\odd(t,\mathbf x)\in\Lambda^\odd T^*_{x_0}M\otimes\eend_{\Cl}(\mathcal E_0)$.
Note that in view of \eqref{E:X2} we have
\begin{equation}\label{E:X3}
\tilde r_1(1,0)=\sum_{i\geq0}\sigma(\tilde k_i)_{[2i-1]}(x_0)
=\tilde r_1^\odd(1,0).
\end{equation}
It thus suffices to determine $\tilde r_1^\odd$.
Since
\begin{align*}
(K^\even(q\tilde r_0))(t,\mathbf x)&\in\Lambda^\even T^*_{x_0}M\otimes\eend_{\Cl}(\mathcal E_0)
\\
(\mathsf A'(q\tilde r_0))(t,\mathbf x)&\in\Lambda^\odd T^*_{x_0}M\otimes\eend_{\Cl}(\mathcal E_0)
\end{align*}
we must have $(\partial_t+\mathsf K)(q\tilde r_1^\odd)=-\mathsf A'(q\tilde r_0)$.
We make the following ansatz, we suppose that $\tilde r_1^\odd=B\tilde r_0$
with $B\in C^\infty(\R,\Lambda^\odd T^*_{x_0}M\otimes\eend_{\Cl}(\mathcal E_0))$. Then
\begin{align*}
(\partial_t+\mathsf K)(q\tilde r_1^\odd)
&=(\partial_tB)q\tilde r_0+B\partial_t(q\tilde r_0)+\mathsf K(Bq\tilde r_0)
\\&=(\partial_tB)q\tilde r_0-B\mathsf K(q\tilde r_0)+\mathsf K(Bq\tilde r_0)
\\&=
(\partial_tB)q\tilde r_0-B\mathsf Fq\tilde r_0+\mathsf FBq\tilde r_0
\\&=
\bigl(\partial_tB+\ad(\mathsf F)B\bigr)q\tilde r_0
\end{align*}
Hence we have to solve $\partial_tB=\ad(-\mathsf F)B-\mathsf A'$ with
initial condition $B(0)=0$. This is easily carried out and we
find the solution:
$$
B(t)=-\frac{e^{\ad(-t\mathsf F)}-1}{\ad(-\mathsf F)}\mathsf A'
$$
Thus $qB\tilde r_0$ satisfies $(\partial_t+\mathsf K)(qB\tilde r_0)=-\mathsf A'(q\tilde r_0)$
with initial condition $(B\tilde r_0)(0,0)=0$.
Again, the uniqueness of formal solutions of the heat equation guarantees
that we actually have $\tilde r_1^\odd=B\tilde r_0$. Setting $t=1$, $\mathbf x=0$
and using \eqref{E:fgh} we get
$$
\tilde r_1^\odd(1,0)=B(1)\tilde r_0(1,0)
=-{\det}^{1/2}\biggl(\frac{\mathsf R/2}{\sinh(\mathsf R/2)}\biggr)\wedge
\biggl(\frac{e^{\ad(-\mathsf F)}-1}{\ad(-\mathsf F)}\mathsf A'\biggr)\wedge\exp(-\mathsf F).
$$
Using \eqref{E:X3} we conclude
$$
\sum_{i\geq0}\sigma(\tilde k_i)_{[2i-1]}=
-\hat A_g\wedge
\biggl(\frac{e^{\ad(-F^{\mathcal E/S}_\nabla)}-1}{\ad(-F^{\mathcal E/S}_\nabla)}\nabla A\biggr)
\wedge\exp\bigl(-F^{\mathcal E/S}_\nabla\bigr).
$$
The Bianchi identity $\nabla F^{\mathcal E/S}_\nabla=0$ implies
$\nabla\exp(-F^{\mathcal E/S}_\nabla)=0$, $\nabla\ad(-F^{\mathcal
E/S}_\nabla)=0$, and similarly $d\hat A_g=0$, from which we finally
obtain \eqref{E:thmiii}.
\end{proof}

\subsection*{Certain heat traces}

Since $\mathcal E$ is $\Z_2$-graded we have a \emph{super trace}
$\str_{\mathcal E}:\Gamma(\eend(\mathcal E))\to\Omega^0(M;\mathbb C)$.
If $n$ is even we will also make use of the so called
\emph{relative super trace,} see \cite[Definition~3.28]{BGV92},
$\str_{\mathcal E/S}:\Gamma(\eend_{\Cl}(\mathcal E))\to\Omega^0(M;\mathcal O_M^\C)$
$$
\str_{\mathcal E/S}(b):=2^{-n/2}\str_{\mathcal E}(c(\Gamma)b).
$$ 
Here $\Gamma\in\Gamma(\Cl\otimes\mathcal O_M^\C)$ denotes the chirality element, see
\cite[Lemma~3.17]{BGV92}. With respect to a local orthonormal frame $\{e_i\}$ of $TM$
and its dual local coframe $\{e^i\}$
the chirality element $\Gamma$ is given as $\ii^{n/2}e^1\cdots e^n$ times the orientation
of $(e_1,\dotsc,e_n)$.
 This relative super trace gives rise to 
$$
\str_{\mathcal E/S}:\Omega^*(M;\eend_{\Cl}(\mathcal E))\to\Omega^*(M;\mathcal O_M^\C)
$$
which will be denoted by the same symbol. For every $\phi\in\Gamma(\eend(\mathcal E))$ 
we have
\begin{equation}\label{E:Q1}
(\str_{\mathcal E}(\phi))\cdot\Omega_g
=(\ii/2)^{-n/2}\str_{\mathcal E/S}\bigl(\sigma(\phi)_{[n]}\bigr),
\end{equation}
where $\Omega_g\in\Omega^n(M;\mathcal O_M^\C)$ denotes the volume density
associated with $g$. To see \eqref{E:Q1} note first that
\begin{equation}\label{E:Q2}
\Cl_{n-1}=[\Cl,\Cl]
\end{equation}
where $\Cl_k$ denotes the filtration on $\Cl$, see \cite[Proof of Proposition~3.21]{BGV92}. 
Hence both sides of \eqref{E:Q1} vanish on $\Gamma(\Cl_{n-1}\otimes\eend_{\Cl}(\mathcal E))$.
It remains to check \eqref{E:Q1} on sections of $\Cl/\Cl_{n-1}\otimes\eend_{\Cl}(\mathcal E)$,
but for these the desired equality follows immediately from the definition of the relative
super trace.

\begin{lemma}\label{L:str}
Let $D_{\mathbb A}$ be a Dirac operator 
and $\tilde k_i\in\Gamma(\eend(\mathcal E))$ as in Theorem~\ref{T:ass_exp}.
Moreover, let $\Phi\in\Gamma(\eend(\mathcal E))$. Then, for even $n$, we have
$$
\LIM_{t\to0}\str\bigl(\Phi e^{-tD^2_{\mathbb A}}\bigr)
=(2\pi\ii)^{-n/2}\int_M\str_{\mathcal E/S}\bigl(\sigma(\Phi\tilde k_{n/2})_{[n]}\bigr),
$$
whereas $\LIM_{t\to0}\str\bigl(\Phi e^{-tD^2_{\mathbb A}}\bigr)=0$ if $n$ is odd.
Here $\LIM$ denotes the renormalized limit \cite[Section~9.6]{BGV92} which in this 
case is just the constant term in the asymptotic expansion for $t\to0$.
\end{lemma}

\begin{proof}
For odd $n$ this follows immediately from \eqref{E:ass_exp}. So assume $n$ is even.
Recall from \cite[Proposition~2.32]{BGV92} that
\begin{equation}\label{E:007}
\str\bigl(\Phi e^{-tD^2_{\mathbb A}}\bigr)
=\int_M\str_{\mathcal E}(\Phi k_t)\cdot\Omega_g.
\end{equation}
Combining this with \eqref{E:Q1} we obtain
$$
\str\bigl(\Phi e^{-tD^2_{\mathbb A}}\bigr)
=(\ii/2)^{-n/2}\int_M\str_{\mathcal E/S}\bigl(\sigma(\Phi k_t)_{[n]}\bigr)
$$
We thus get an asymptotic expansion, see \eqref{E:ass_exp},
$$
\str\bigl(\Phi e^{-tD^2_{\mathbb A}}\bigr)
\sim(2\pi\ii t)^{-n/2}\sum_{i\geq0}t^i\int_M\str_{\mathcal E/S}\bigl(\sigma(\Phi\tilde k_i)_{[n]}\bigr)
\qquad\text{as $t\to0$,}
$$
from which the desired formula follows at once.
\end{proof}

\begin{corollary}\label{C:tr_ass_exp}
Let $D_{\mathbb A}$ be a Dirac operator as in Theorem~\ref{T:ass_exp}.
Moreover, let $U\in\Gamma(\eend_{\Cl}(\mathcal E))$. Then, for even $n$, we have
\begin{equation}\label{E:aaa}
\LIM_{t\to0}\str(Ue^{-tD_{\mathbb A}^2})
=(2\pi\ii)^{-n/2}\int_M\hat A_g\wedge\str_{\mathcal E/S}
\bigl(U\exp(-F^{\mathcal E/S}_\nabla)\bigr),
\end{equation}
whereas $\LIM_{t\to0}\str(Ue^{-tD_{\mathbb A}^2})=0$ if $n$ is odd.
\end{corollary}

\begin{proof}
For odd $n$ this follows immediately from Lemma~\ref{L:str}. So assume $n$ is even.
Since $\sigma(U\tilde k_i)_{[n]}=U\sigma(\tilde k_i)_{[n]}$ Theorem~\ref{T:ass_exp}
yields
$$
\str_{\mathcal E/S}\bigl(\sigma(U\tilde k_{n/2})_{[n]}\bigr)
=\Bigl(\hat A_g\wedge\str_{\mathcal E/S}\bigl(U\exp(-F^{\mathcal E/S}_\nabla)\bigr)\Bigr)_{[n]}.
$$
Equation~\eqref{E:aaa} then follows from Lemma~\ref{L:str}.
\end{proof}

\begin{corollary}\label{C:tr_ass_exp_2}
Let $D_{\mathbb A}$ be a Dirac operator as in Theorem~\ref{T:ass_exp}.
Moreover, suppose $V\in\Omega^1(M;\eend_{\Cl}(\mathcal E))$,
let $c(V)\in\Gamma(\eend(\mathcal E))$ denote Clifford multiplication
with $V$, and consider $\nabla V\in\Omega^2(M;\eend_{\Cl}(\mathcal E))$. Then,
for even $n$, we have
\begin{multline}\label{E:bbb}
\LIM_{t\to0}\str\Bigl(c(V)e^{-tD_{\mathbb A}^2}\Bigr)=
\\
-(2\pi\ii)^{-n/2}\int_M
\hat A_g\wedge
\str_{\mathcal E/S}\biggl(
\biggl(\frac{e^{\ad(-F^{\mathcal E/S}_\nabla)}-1}{\ad(-F^{\mathcal
E/S}_\nabla)}A\biggr)\wedge
\exp\bigl(-F^{\mathcal E/S}_\nabla\bigr)\wedge\nabla V\biggr),
\end{multline}
whereas $\LIM_{t\to0}\str\Bigl(c(V)e^{-tD_{\mathbb A}^2}\Bigr)=0$ if $n$ is odd.
\end{corollary}

\begin{proof}
If $n$ is odd the statement follows immediately from Lemma~\ref{L:str}. So assume $n$ is even.
Since $\sigma\bigl(c(V)\tilde k_i\bigr)_{[n]}=V\wedge\sigma(\tilde k_i)_{[n-1]}$
Theorem~\ref{T:ass_exp} yields
\begin{align*}
\str_{\mathcal E/S}\Bigl(\sigma\bigl(c(&V)\tilde k_{n/2}\bigr)_{[n]}\Bigr)
\\&=-\str_{\mathcal E/S}\biggl(V\wedge\nabla\biggl(\hat A_g\wedge
\biggl(\frac{e^{\ad(-F^{\mathcal E/S}_\nabla)}-1}{\ad(-F^{\mathcal E/S}_\nabla)}
A\biggr)
\wedge\exp\bigl(-F^{\mathcal E/S}_\nabla\bigr)\biggr)\biggr)_{[n]}.
\\&=-\str_{\mathcal E/S}\biggl(\nabla\biggl(\hat A_g\wedge
\biggl(\frac{e^{\ad(-F^{\mathcal E/S}_\nabla)}-1}{\ad(-F^{\mathcal E/S}_\nabla)}
A\biggr)
\wedge\exp\bigl(-F^{\mathcal E/S}_\nabla\bigr)\biggr)\wedge V\biggr)_{[n]}.
\end{align*}
Applying Lemma~\ref{L:str} and using Stokes' theorem we obtain \eqref{E:bbb}.
\end{proof}

\section{Application to Laplacians}\label{S:anomprop}

Below we will see that the Laplacians $\Delta_{E,g,b}$ introduced in 
Section~\ref{S:anator} are the squares of 
Dirac operators of the kind considered in Section~\ref{S:ass_ex}.
Applying Corollaries~\ref{C:tr_ass_exp} and \ref{C:tr_ass_exp_2} will lead to
a proofs of Propositions~\ref{P:strbD} and \ref{P:strgD}, respectively.

\subsection*{The exterior algebra as Clifford module}

Let $(M,g)$ be a Riemannian manifold of dimension $n$. In order to understand the
Clifford module structure of $\Lambda:=\Lambda^*T^*M$ we first note that
$\Lambda$ is a Clifford module for $\hat\Cl:=\Cl(T^*M,-g)$ too.
Let us write $\hat c$ for the Clifford multiplication of $\hat\Cl$
on $\Lambda$. Explicitly, for $a\in T_x^*M\subseteq\hat\Cl$
and $\alpha\in\Lambda^*T^*_xM$ we have $\hat
c(a)\alpha
=a\wedge\alpha+i_{\sharp a}\alpha$, where $\sharp a:=g^{-1}a\in T_xM$ and
$i_{\sharp a}$ denotes contraction with $\sharp a$. 
It follows from this formula that every $\hat c(a)$ 
commutes with the Clifford action of $\Cl$. We thus obtain an isomorphism
of $\Z_2$-graded filtered algebras
$$
\hat c:\hat\Cl\to\eend_{\Cl}(\Lambda).
$$
Let us write
$$
\hat\sigma:\hat\Cl\to\Lambda,\qquad\hat\sigma(a):=\hat c(a)\cdot1
$$
for the symbol map of $\hat\Cl$.

As in \eqref{E:Q6} define $R^{\hat\Cl}\in\Omega^2(M;\hat\Cl)$ by
$$
R^{\hat\Cl}(X,Y):=-\tfrac14\sum_{i,j}g\bigl(R(X,Y)e_i,e_j\bigr)\hat c^i\hat c^j
$$
where $X$ and $Y$ are two vector fields, $\{e_i\}$ is a local orthonormal frame, 
$\{e^i\}$ denotes its dual local coframe, and $\hat c^i:=\hat c(e^i)$.
For the twisting curvature 
$F_{\nabla^g}^{\Lambda/S}\in\Omega^2(M;\eend_{\Cl}(\Lambda))$ 
we then have, see \cite[Page~145]{BGV92},
\begin{equation}\label{E:Q5}
F^{\Lambda/S}_{\nabla^g}
=(1\otimes\hat c)(R^{\hat\Cl})
\in\Omega^2(M;\eend_{\Cl}(\Lambda)),
\end{equation}
where $(1\otimes\hat c):\Omega(M;\hat\Cl)\to\Omega(M;\eend_{\Cl}(\Lambda))$.
Indeed, the curvature of $\Lambda$,
$R^\Lambda\in\Omega^2(M;\eend(\Lambda))$, can be written as
$$
R^\Lambda(X,Y)=\sum_{i,j}g\bigl(R(X,Y)e_i,e_j\bigr)\tfrac12(\varepsilon^j\iota^i-\varepsilon^i\iota^j)
\in\Gamma(\eend(\Lambda))
$$
where $\varepsilon^j\in\Gamma(\eend(\Lambda))$ denotes exterior multiplication with $e^j$, 
and $\iota^i\in\Gamma(\eend(\Lambda))$ denotes contraction with $e_i$.
Using $\varepsilon^i=\frac12(c^i+\hat c^i)$ and $\iota^i=-\frac12(c^i-\hat c^i)$ 
one easily deduces
$$
\tfrac12(\varepsilon^j\iota^i-\varepsilon^i\iota^j)
=\tfrac14\bigl(\tfrac12(c^ic^j-c^jc^i)\bigr)
-\tfrac14\bigl(\tfrac12(\hat c^i\hat c^j-\hat c^j\hat c^i)\bigr)
$$
from which we read off \eqref{E:Q5}, see \eqref{E:Q4}. Also note that
we have
\begin{equation}\label{E:Q7}
(1\otimes\hat\sigma)(R^{\hat\Cl})=-\tfrac12R\in\Omega^2(M;\Lambda^2T^*M),
\end{equation}
where $(1\otimes\hat\sigma):\Omega(M;\hat\Cl)\to\Omega(M;\Lambda)$.

If $n$ is even then the relative super trace
$$
\str_{\Lambda/S}:\eend_{\Cl}(\Lambda)\to\mathcal O_M^\C
$$ 
is given by
\begin{equation}\label{E:jkl}
\str_{\Lambda/S}(\hat c(a))=(\ii/2)^{-n/2}T(\hat\sigma(a))
\qquad a\in\hat\Cl,
\end{equation}
where $T:\Lambda\to\mathcal O_M^\C$ denotes the Berezin integration
associated with $g$. 
Indeed, since $[\hat\Cl,\hat\Cl]=\hat\Cl_{n-1}$, see \eqref{E:Q2},
both sides of \eqref{E:jkl} vanish for $a\in\hat\Cl_{n-1}$. 
Checking \eqref{E:jkl} on $\hat\Cl/\hat\Cl_{n-1}$ is straight forward.
We will also make use of the formula
\begin{equation}\label{E:Q9}
\str_{\Lambda/S}\bigl(\exp((1\otimes\hat c)a)\bigr)
=(\ii/2)^{-n/2}T\bigl(\exp_\Lambda((1\otimes\hat\sigma)a)\bigr)\qquad
a\in\Omega^2(M;\hat\Cl_2),
\end{equation}
where $1\otimes\hat c:\Omega(M;\hat\Cl)\to\Omega(M;\eend_{\Cl}(\Lambda))$,
$1\otimes\hat\sigma:\Omega(M;\hat\Cl)\to\Omega(M;\Lambda)$ and
$T:\Omega(M;\Lambda)\to\Omega(M;\mathcal O_M^\C)$ denotes Berezin integration.
To check this equation note that the assumption on the form degree and the
filtration degree of $a$ implies:
\begin{align*}
\str_{\Lambda/S}\bigl(\exp((1\otimes\hat c)a)\bigr)
&=\str_{\Lambda/S}\bigl(\tfrac1{n!}((1\otimes\hat c)a)^{n/2}\bigr)
\\
T\bigl(\exp_\Lambda((1\otimes\hat\sigma)a)\bigr)
&=T\bigl(\tfrac1{n!}((1\otimes\hat\sigma)a)^{n/2}\bigr)
\end{align*}
Using the fact that $1\otimes\hat c$ is an algebra isomorphism and
\eqref{E:jkl} we obtain
\begin{multline*}
\str_{\Lambda/S}\bigl(((1\otimes\hat c)a)^{n/2}\bigr)
=\str_{\Lambda/S}\bigl((1\otimes\hat c)(a^{n/2})\bigr)
\\=(\ii/2)^{-n/2}T\bigl((1\otimes\hat\sigma)(a^{n/2})\bigr)
=(\ii/2)^{-n/2}T\bigl(((1\otimes\hat\sigma)a)^{n/2}\bigr)
\end{multline*}
where we made use of the fact that $1\otimes\hat\sigma$ induces
an isomorphism on the level of associated graded algebras, for the last
equality. Combined with the previous two equations this proves \eqref{E:Q9}.

\begin{lemma}\label{L:char}
Let $(M,g)$ be a Riemannian manifold of even dimension $n$. Then\footnote{Since the degree $0$ part of $\hat A_g$ is $1$,
this formula is easily seen to be equivalent to 
$\ec(g)=(2\pi\ii)^{-n/2}\hat A_g\wedge\str_{\Lambda/S}\bigl(\exp\bigl(-F^{\Lambda/S}_{\nabla^g}\bigr)\bigr)$ which
can be found in \cite[Proposition~4.6]{BGV92}.}
$$
\ec(g)=(2\pi\ii)^{-n/2}\str_{\Lambda/S}\bigl(\exp\bigl(-F^{\Lambda/S}_{\nabla^g}\bigr)\bigr).
$$
\end{lemma}

\begin{proof}
Consider the negative of the Riemannian
curvature $-R\in\Omega^2(M;\Lambda^2T^*M)$ and its exponential
$\exp_\Lambda(-R)\in\Omega(M;\Lambda)$. Recall that
$$
\ec(g):=(2\pi)^{-n/2}T(\exp_\Lambda(-R))\in\Omega^n(M;\mathcal O_M^\C).
$$
Using \eqref{E:Q7}, \eqref{E:Q9} and \eqref{E:Q5} we conclude:
\begin{align*}
\ec(g)
&=(2\pi)^{-n/2}T\Bigl(\exp_\Lambda\bigl((1\otimes\hat\sigma)(2R^{\hat\Cl})\bigr)\Bigr)
\\&=(-\pi)^{-n/2}T\Bigl(\exp_\Lambda\bigl((1\otimes\hat\sigma)(-R^{\hat\Cl})\bigr)\Bigr)
\\&=(2\pi\ii)^{-n/2}\str_{\Lambda/S}\Bigl(\exp\bigl(-(1\otimes\hat c)(R^{\hat\Cl})\bigr)\Bigr)
\\&=(2\pi\ii)^{-n/2}\str_{\Lambda/S}\bigl(\exp\bigl(-F^{\Lambda/S}_{\nabla^g}\bigr)\bigr)
\qedhere
\end{align*}
\end{proof}

\begin{lemma}\label{L:char2}
Let $(M,g)$ be a Riemannian manifold of even dimension $n$. Suppose
$\tilde\xi\in\Gamma(T^*M\otimes T^*M)$ is symmetric, use 
$1\otimes\hat c:T^*M\otimes T^*M\to T^*M\otimes\eend_{\Cl}(\Lambda)$
to define $V:=\frac12(1\otimes\hat c)(\tilde\xi)\in\Omega^1(M;\eend_{\Cl}(\Lambda))$,
and consider $\nabla^g V\in\Omega^2(M;\eend_{\Cl}(\Lambda))$.
Then, for every closed one form $\omega\in\Omega^1(M;\C)$, we have
$$
\omega\wedge(\partial_2\cs)(g,\tilde\xi)
=\tfrac12(2\pi\ii)^{-n/2}\str_{\Lambda/S}\Bigl(\hat c(\omega)\wedge
\exp\bigl(-F^{\Lambda/S}_{\nabla^g}\bigr)\wedge\nabla^gV\Bigr)
$$
in $\Omega^n(M;\mathcal O_M^\C)/d\Omega^{n-1}(M;\mathcal O_M^\C)$.
\end{lemma}

\begin{proof}
Set $\tilde M:=M\times\R$ and consider the two natural projections $p:\tilde M\to M$
and $t:\tilde M\to\R$. Consider the bundle $\tilde TM:=p^*TM$ over $\tilde M$, and equip
it with the fiber metric $\tilde g:=p^*(g+t\tilde\xi)$.
For sufficiently small $t$, this will indeed be non-degenerate.
For $t\in\R$ let $\inc_t:M\to\tilde M$ denote the inclusion $x\mapsto(x,t)$.
Define a connection $\tilde\nabla$ on $\tilde TM$ so that
$\inc_t^*\tilde\nabla=\nabla^{g+t\tilde\xi}$ for sufficiently small $t$,
where $\nabla^{g+t\tilde\xi}$ denotes the Levi--Civita connection
of $g+t\tilde\xi$, and so that $\tilde\nabla_{\partial t}=\partial_t+\frac12\tilde g^{-1}(p^*\tilde\xi)$.
It is not hard to check that $\tilde g$
is parallel with respect to $\tilde\nabla$, i.e.\ $\tilde\nabla\tilde g=0$.
Let $\ec(\tilde TM,\tilde g,\tilde\nabla)\in\Omega^n(\tilde M;\mathcal O_{\tilde TM})$
denote the Euler form of this Euclidean bundle. Recall that
$$
\cs(g,g+\tau\tilde\xi)=\int_0^\tau\inc_t^*i_{\partial_t}\ec(\tilde TM,\tilde g,\tilde\nabla)dt
$$
and thus
\begin{align*}
(\partial_2\cs)(g,\tilde\xi)
&=\inc_0^*i_{\partial_t}\ec(\tilde TM,\tilde g,\tilde\nabla)
\\&=(2\pi)^{-n/2}\cdot\inc_0^*i_{\partial_t}T\bigl(\exp_\Lambda(-R^{\tilde\nabla})\bigr)
\\&=(-2\pi)^{-n/2}\cdot T\bigl(\inc_0^*i_{\partial_t}\exp_\Lambda(R^{\tilde\nabla})\bigr)
\\&=(-2\pi)^{-n/2}\cdot T\bigl(\inc_0^*\bigl(\exp_\Lambda(R^{\tilde\nabla})\wedge i_{\partial_t}R^{\tilde\nabla}\bigr)\bigr)
\\&=(-2\pi)^{-n/2}\cdot T\bigl(\exp_\Lambda(\inc_0^*R^{\tilde\nabla})\wedge\inc_0^*i_{\partial_t}R^{\tilde\nabla}\bigr)
\end{align*}
where $R^{\tilde\nabla}\in\Omega^2(\tilde M;\Lambda^2\tilde TM)$
denotes the curvature of $\tilde\nabla$. Let 
$$
S:\Lambda\otimes\Lambda\to\Lambda\otimes\Lambda,\qquad
S(\alpha\otimes\beta):=(-1)^{|\alpha||\beta|}\beta\otimes\alpha
$$
denote the isomorphism of graded algebras obtained by interchanging variables.
Consider $\tilde\xi\in\Omega^1(M;T^*M)$, $\nabla^g\tilde\xi\in\Omega^2(M;T^*M)$
and $S(\nabla^g\tilde\xi)\in\Omega^1(M;\Lambda^2T^*M)$. With this notation we have
\begin{align*}
\inc_0^*R^{\tilde\nabla}&=R\in\Omega^2(M;\Lambda^2T^*M)
\\
\inc_0^*i_{\partial _t}R^{\tilde\nabla}&=S(\tfrac12\nabla^g\tilde\xi)
\in\Omega^1(M;\Lambda^2T^*M)
\end{align*}
where $R$ denotes the Riemannian curvature of $g$.
We obtain
$$
(\partial_2\cs)(g,\tilde\xi)
=(-2\pi)^{-n/2}T\bigl(\exp_\Lambda(R)\wedge S(\tfrac12\nabla^g\tilde\xi)\bigr)
$$
and wedging with $\omega$ we get
$$
\omega\wedge(\partial_2\cs)(g,\tilde\xi)
=(-2\pi)^{-n/2}T\bigl(\exp_\Lambda(R)\wedge\omega\wedge S(\tfrac12\nabla^g\tilde\xi)\bigr).
$$
Next, note that for $a\in\Lambda^nT^*M\otimes\Lambda^nT^*M$ we have $T(S(a))=T(a)$, for $n$ is supposed to be even.
Together with the symmetries of the Riemann curvature, $S(R)=R$, we obtain
\begin{align*}
\omega\wedge(\partial_2\cs)(g,\tilde\xi)
&=(-2\pi)^{-n/2}T\Bigl(S\bigl(\exp_\Lambda(R)\wedge\omega\wedge S(\tfrac12\nabla^g\tilde\xi)\bigr)\Bigr)
\\&=(-2\pi)^{-n/2}T\Bigl(S(\exp_\Lambda(R))\wedge S(\omega)\wedge\tfrac12\nabla^g\tilde\xi\Bigr)
\\&=(-2\pi)^{-n/2}T\Bigl(\exp_\Lambda(R)\wedge(1\otimes\omega)\wedge\tfrac12\nabla^g\tilde\xi\Bigr)
\\&=(-2\pi)^{-n/2}\frac\partial{\partial s}\Bigm|_{s=0}T\Bigl(\exp_\Lambda\Bigl(R+s(1\otimes\omega)\wedge\tfrac12\nabla^g\tilde\xi\Bigr)\Bigr)
\end{align*}
In view of \eqref{E:Q7} we have:
$$
R+s(1\otimes\omega)\wedge\tfrac12(\nabla^g\tilde\xi)
=(1\otimes\hat\sigma)\Bigl(-2R^{\hat\Cl}+s(1\otimes\omega)\wedge\tfrac12\nabla^g\tilde\xi\Bigr)
$$
Moreover, using \eqref{E:Q5} and $(1\otimes\hat c)(\frac12\nabla^g\tilde\xi)=\nabla^gV$ we also have:
$$
(1\otimes\hat c)\Bigl(-2R^{\hat\Cl}+s(1\otimes\omega)\wedge\tfrac12\nabla^g\tilde\xi\Bigr)
=-2F^{\Lambda/S}_{\nabla^g}+s\hat c(\omega)\wedge\nabla^gV
$$
Using these two equations and applying \eqref{E:Q9} we obtain
\begin{align*}
\omega\wedge(\partial_2\cs)(g,\tilde\xi)
&=(4\pi\ii)^{-n/2}\frac\partial{\partial s}\Bigm|_{s=0}
\str_{\Lambda/S}\Bigl(\exp\Bigl(-2F^{\Lambda/S}_{\nabla^g}+s\hat c(\omega)\wedge\nabla^gV
\Bigr)\Bigr)
\\&=(4\pi\ii)^{-n/2}
\str_{\Lambda/S}\Bigl(\exp\bigl(-2F^{\Lambda/S}_{\nabla^g}\bigr)\wedge\hat c(\omega)\wedge\nabla^gV
\Bigr)
\\&=\tfrac12(2\pi\ii)^{-n/2}
\str_{\Lambda/S}\Bigl(\exp\bigl(-F^{\Lambda/S}_{\nabla^g}\bigr)\wedge\hat c(\omega)\wedge\nabla^gV
\Bigr)
\\&=\tfrac12(2\pi\ii)^{-n/2}
\str_{\Lambda/S}\Bigl(\hat c(\omega)\wedge\exp\bigl(-F^{\Lambda/S}_{\nabla^g}\bigr)\wedge\nabla^gV
\Bigr)
\qedhere
\end{align*}
\end{proof}

\subsection*{The Laplacians as squares of Dirac operators}

Let $E$ be a flat complex
vector bundle equipped with a fiber wise non-degenerate symmetric bilinear form $b$. 
Let $\nabla^E$ denote the flat connection on $E$.
Consider $b^{-1}\nabla^Eb\in\Omega^1(M;\eend(E))$ and introduce
the connection, cf.\ \cite[Section~4]{BZ92},
$$
\nabla^{E,b}:=\nabla^E+\tfrac12b^{-1}\nabla^Eb
$$ 
on $E$. 
Consider the Clifford bundle $\mathcal E:=\Lambda\otimes E$
with Clifford connection
$$
\nabla^{E,g,b}:=\nabla^g\otimes1_E+1_\Lambda\otimes\nabla^{E,b}.
$$
Since $(\nabla^{E,g,b})^2=(\nabla^g)^2+(\nabla^{E,b})^2$ the 
twisting curvature is
\begin{equation}\label{E:FF}
F^{\mathcal E/S}_{\nabla^{E,g,b}}
=F^{\Lambda/S}_{\nabla^g}+(\nabla^{E,b})^2.
\end{equation}
Since the two summands commute we obtain
\begin{equation}\label{E:expexp}
\exp\bigl(-F^{\mathcal E/S}_{\nabla^{E,g,b}}\bigr)
=\exp\bigl(-F^{\Lambda/S}_{\nabla^g}\bigr)
\wedge\exp\bigl(-(\nabla^{E,b})^2\bigr).
\end{equation}
An easy computation shows that the
Dirac operator associated to the Clifford connection $\nabla^{E,g,b}$ is
\begin{equation}\label{E:dir}
D_{\nabla^{E,g,b}}
=d_E+d_{E,g,b}^\sharp
+\hat c\bigl(\tfrac12b^{-1}\nabla^Eb\bigr).
\end{equation}
Setting
\begin{equation}\label{E:Q11}
A_{E,g,b}:=-\hat c\bigl(\tfrac12b^{-1}\nabla^Eb\bigr)
\in\Omega^0(M;\eend^-_{\Cl}(\mathcal E))
\end{equation}
we obtain a Clifford super connection
\begin{equation}\label{E:agb}
\mathbb A_{E,g,b}:=\nabla^{E,g,b}+A_{E,g,b}.
\end{equation}
For the associated Dirac
operator $D_{\mathbb A_{E,g,b}}=d_E+d^\sharp_{E,g,b}$ we find
\begin{equation}\label{E:DD2}
(D_{\mathbb A_{E,g,b}})^2
=(d_E+d_{E,g,b}^\sharp)^2
=\Delta_{E,g,b}.
\end{equation}
So we see that the Laplacians introduced in Section~\ref{S:anator} are
indeed squares of Dirac operators of the type considered in
Theorem~\ref{T:ass_exp}.

\subsection*{Proof of Proposition~\ref{P:strbD}}

For odd $n$ the statement follows immediately from Lem\-ma~\ref{L:str}.
So let us assume that $n$ is even.
We will apply Corollary~\ref{C:tr_ass_exp} to the Clifford super connection
\eqref{E:agb} and $U:=\phi$. 
From \eqref{E:expexp} and Lemma~\ref{L:char} we get:
\begin{multline*}
\str_{\mathcal E/S}\Bigl(\phi\exp\bigl(-F^{\mathcal E/S}_{\nabla^{E,g,b}}\bigr)\Bigr)
=\str_{\mathcal E/S}\Bigl(\phi\exp\bigl(-F^{\Lambda/S}_{\nabla^g}\bigr)
\wedge\exp\bigl(-(\nabla^{E,b})^2\bigr)
\Bigr)
\\
=\str_{\Lambda/S}\Bigl(
\exp\bigl(-F^{\Lambda/S}_{\nabla^g}\bigr)
\Bigr)
\wedge\tr_E\Bigl(\phi\exp\bigl(-(\nabla^{E,b})^2\bigr)\Bigr)
=(2\pi\ii)^{n/2}\ec(g)\tr(\phi)
\end{multline*}
Here we also used the fact that the form
$\str_{\Lambda/S}\bigl(\exp\bigl(-F^{\Lambda/S}_{\nabla^g}\bigr)\bigr)
=(2\pi\ii)^{n/2}\ec(g)$ has degree $n$, and thus the only contributing part of
$\tr_E\bigl(\phi\exp\bigl(-(\nabla^{E,b})^2\bigr)\bigr)$ is the one of form degree $0$, which is just $\tr(\phi)$.
Using again the fact that $\ec(g)$ has maximal form degree, we conclude
$$
(2\pi\ii)^{-n/2}\hat A_g\wedge
\str_{\mathcal E/S}\Bigl(\phi\exp\bigl(-F^{\mathcal E/S}_{\nabla^{E,g,b}}\bigr)\Bigr)
=\tr(\phi)\ec(g),
$$
since the degree $0$ part of $\hat A_g$ is just $1$.
Proposition~\ref{P:strbD} now follows from 
Corollary~\ref{C:tr_ass_exp} and \eqref{E:DD2}.

\subsection*{Proof of Proposition~\ref{P:strgD}}

For odd $n$ the statement follows immediately from Lem\-ma~\ref{L:str}.
So let us assume $n$ is even.
Consider $\tilde\xi:=g\xi\in\Gamma(T^*M\otimes T^*M)$, and use the bundle map
$1\otimes\hat c:T^*M\otimes T^*M\to T^*M\otimes\eend_{\Cl}^-(\Lambda)$
to define 
$$
V:=\tfrac12(1\otimes\hat c)(\tilde\xi)\in\Omega^1(M;\eend_{\Cl}^-(\Lambda)).
$$
We claim
\begin{equation}\label{E:ghj}
c(V)=\Lambda^*\xi-\tfrac12\tr(\xi).
\end{equation}
To check this let $\{e_i\}$ be a local orthonormal frame and let $\{e^i\}$ be
its dual local coframe. Then
$$
\Lambda^*\xi=\sum_{i,j}g(\xi e_i,e_j)\tfrac12\bigl(\varepsilon^i\iota^j+\varepsilon^j\iota^i\bigr),
$$
where $\varepsilon^i\in\Gamma(\eend(\Lambda))$ denotes exterior multiplication with $e^i$,
and $\iota^i\in\Gamma(\eend(\Lambda))$ denotes contraction with $e_i$. Writing
$c^i:=c(e^i)$, $\hat c^i:=\hat c(e^i)$ and using $\varepsilon^i=\frac12(c^i+\hat c^i)$
as well as $\iota^i=-\frac12(c^i-\hat c^i)$ one easily checks
$$
\tfrac12\bigl(\varepsilon^i\iota^j+\varepsilon^j\iota^i\bigr)
=\tfrac14\bigl(c^i\hat c^j+c^j\hat c^i\bigr)+\tfrac12\delta^{ij}.$$
We conclude
$$
\Lambda^*\xi
=\sum_{i,j}g(\xi e_i,e_j)\Bigl(\tfrac14\bigl(c^i\hat c^j+c^j\hat c^i\bigr)+\tfrac12\delta^{ij}\Bigr)
=\sum_{i,j}g(\xi e_i,e_j)\tfrac12 c^i\hat c^j+\tfrac12\tr(\xi).
$$
On the other hand we clearly have
$c(V)=\sum_{i,j}g(\xi e_i,e_j)\tfrac12 c^i\hat c^j$
and thus \eqref{E:ghj} is established.
We will apply Corollary~\ref{C:tr_ass_exp_2} to the Clifford super connection
\eqref{E:agb} and this $V$.

Next we claim that for all integers $k\geq1$ and $l\geq0$ we have
\begin{equation}\label{E:fin}
\str_{\mathcal E/S}\Bigl(
\bigl(\bigl(\ad(-F^{\mathcal E/S}_{\nabla^{E,g,b}})\bigr)^kA_{E,g,b}\bigr)
\wedge\bigl(-F^{\mathcal E/S}_{\nabla^{E,g,b}}\bigr)^l\wedge\nabla V\Bigr)
=0.
\end{equation}
To see this let us write $\eend_{\Cl}(\mathcal E)_i$ for the subspace of 
$\eend_{\Cl}(\mathcal E)$ which via the isomorphism $\hat c\otimes1:\hat\Cl\otimes\eend(E)\to
\eend_{\Cl}(\Lambda)\otimes\eend(E)=\eend_{\Cl}(\mathcal E)$
corresponds to the filtration subspace $\hat\Cl_i\otimes\eend(E)$. Then
$-F^{\mathcal E/S}_{\nabla^{E,g,b}}\in\Omega^2(M;\eend_{\Cl}(\mathcal E)_2)$,
$\nabla V\in\Omega^2(M;\eend_{\Cl}(\mathcal E)_1)$ and $A_{E,g,b}\in\Omega^0(M;\eend_{\Cl}(\mathcal E)_1)$. 
Looking at the form degree, we see that \eqref{E:fin} holds whenever $2k+2l+2>n$. Moreover, since $k\geq1$
we have $\bigl(\ad(-F^{\mathcal E/S}_{\nabla^{E,g,b}})\bigr)^kA_{E,g,b}\in\Omega^{2k}(M;\eend_{\Cl}(\mathcal E)_{2k})$, for
$[\hat\Cl_2,\hat\Cl_1]\subseteq\hat\Cl_2$.
Thus, considering the filtration degree, we see that \eqref{E:fin} holds whenever $2k+2l+1<n$, for $\str_{\mathcal E/S}$ vanishes
on $\Omega(M;\eend_{\Cl}(\mathcal E)_{n-1})$. This establishes \eqref{E:fin}. We conclude
\begin{multline}\label{E:AA}
\str_{\mathcal E/S}\Biggl(
\biggl(\frac{e^{\ad(-F^{\mathcal E/S}_{\nabla^{E,g,b}})}-1}
{\ad(-F^{\mathcal E/S}_{\nabla^{E,g,b}})}A_{E,g,b}
\biggr)
\wedge\exp\bigl(-F^{\mathcal E/S}_{\nabla^{E,g,b}}\bigr)\wedge\nabla V\Biggr)
\\=\str_{\mathcal E/S}\Bigl(A_{E,g,b}
\wedge\exp\bigl(-F^{\mathcal E/S}_{\nabla^{E,g,b}}\bigr)\wedge\nabla^gV\Bigr)
\end{multline}
Here we wrote $\nabla V=\nabla^gV$ to emphasize that this form does not depend on the flat connection on $E$,
but only on the Levi--Civita connection. Using \eqref{E:expexp} and $(\nabla^{E,b})^2\in\Omega^2(M;\eend_{\Cl}(\mathcal E)_0)$
and considering form and filtration degree we easily obtain:
\begin{align}
\notag\str_{\mathcal E/S}\Bigl(A_{E,g,b}&
\wedge\exp\bigl(-F^{\mathcal E/S}_{\nabla^{E,g,b}}\bigr)\wedge\nabla^gV\Bigr)
\\\notag&=\str_{\mathcal E/S}\Bigl(A_{E,g,b}
\wedge\exp\bigl(-F^{\Lambda/S}_{\nabla^g}\bigr)\wedge\nabla^gV\Bigr)
\\\label{E:BB}&=\str_{\Lambda/S}\Bigl(\tr_E(A_{E,g,b})
\wedge\exp\bigl(-F^{\Lambda/S}_{\nabla^g}\bigr)\wedge\nabla^gV\Bigr)
\end{align}
Using \eqref{E:Q11} and applying Lemma~\ref{L:char2} to the closed one form
$\tr(b^{-1}\nabla^Eb)$ we find
\begin{align}
\notag\str_{\Lambda/S}\Bigl(\tr_E(&A_{E,g,b})
\wedge\exp\bigl(-F^{\Lambda/S}_{\nabla^g}\bigr)\wedge\nabla^gV\Bigr)
\\\notag&=-\tfrac12\str_{\Lambda/S}\Bigl(\hat c\bigl(\tr(b^{-1}\nabla^Eb)\bigr)
\wedge\exp\bigl(-F^{\Lambda/S}_{\nabla^g}\bigr)\wedge\nabla^gV\Bigr)
\\\label{E:CC}&=-(2\pi\ii)^{n/2}\tr(b^{-1}\nabla^Eb)\wedge(\partial_2\cs)(g,\tilde\xi)
\end{align}
Combining \eqref{E:AA}, \eqref{E:BB} and \eqref{E:CC} we conclude:
\begin{multline*}
-(2\pi\ii)^{-n/2}\hat A_g\wedge\str_{\mathcal E/S}\Biggl(
\biggl(\frac{e^{\ad(-F^{\mathcal E/S}_{\nabla^{E,g,b}})}-1}
{\ad(-F^{\mathcal E/S}_{\nabla^{E,g,b}})}A_{E,g,b}
\biggr)
\wedge\exp\bigl(-F^{\mathcal E/S}_{\nabla^{E,g,b}}\bigr)\wedge\nabla V\Biggr)
\\
=\tr(b^{-1}\nabla^Eb)\wedge(\partial_2\cs)(g,g\xi)
\end{multline*}
Now apply Corollary~\ref{C:tr_ass_exp_2} and use \eqref{E:ghj} as well as \eqref{E:DD2}
to complete the proof of Proposition~\ref{P:strgD}.

\section{The case of non-vanishing Euler--Poincar\'e characteristics}\label{S:gen}

It is not necessary to restrict to manifolds with vanishing Euler
characteristics. In the general situation \cite{BH03, BH06} Euler
structures, coEuler structures, the combinatorial torsion and the analytic
torsion depend on the choice of a base point. Given a 
path connecting two such base points everything associated with the
first base point identifies in an equivariant way with the everything associated to the
other base point. However, these identifications do depend on the homotopy class
of such a path. Below we sketch a natural way to conveniently
deal with this situation.

In general the the set of Euler structures $\Eul_{x_0}(M;\Z)$ depends on a base point $x_0\in M$. 
One defines the set of Euler structures based at $x_0$ as equivalence classes $[X,c]$ where
$X$ is a vector field with non-degenerate zeros and $c\in C^\sing_1(M;\Z)$ is such that
$\partial c=\ec(X)-\chi(M)x_0$. Two such pairs $(X_1,c_1)$ and $(X_2,c_2)$ are equivalent iff
$c_2-c_1=\cs(X_1,X_2)$ mod boundaries. Again this is an affine version of $H_1(M;\Z)$, the action is
defined as in Section~\ref{S:eul}.
Given a path $\sigma$ from $x_0$ to $x_1$, the assignment $[X,c]\mapsto[X,c-\chi(M)\sigma]$ 
defines an $H_1(M;\Z)$-equivariant isomorphism from
$\Eul_{x_0}(M;\Z)$ to $\Eul_{x_1}(M;\Z)$. Since this isomorphism depends on the homotopy class of $\sigma$
only, we can consider the set of Euler structures as a flat principal bundle $\Eul(M;\Z)$ over $M$
with structure group $H_1(M;\Z)$. Its fiber over $x_0$ is just $\Eul_{x_0}(M;\Z)$, and its 
holonomy is given by the composition
$$
\pi_1(M)\to H_1(M;\Z)\xrightarrow{-\chi(M)}H_1(M;\Z).
$$

Similarly, the set of Euler structures with complex coefficients can be
considered as a flat principal bundle $\Eul(M;\C)$ over $M$ with structure
group $H_1(M;\C)$ and holonomy given by the composition
$$
\pi_1(M)\to H_1(M;\Z)\xrightarrow{-\chi(M)} H_1(M;\Z)\to H_1(M;\C).
$$
There is an obvious parallel homomorphism of flat principal bundles over $M$
\begin{equation}\label{E:D}
\iota:\Eul(M;\Z)\to\Eul(M;\C)
\end{equation}
which is equivariant over the homomorphism of structure groups 
$H_1(M;\Z)\to H_1(M;\C)$.

The set of coEuler structures $\Eul^*_{x_0}(M;\C)$ depends
on the choice of a base point $x_0\in M$. It can be defined as the set of equivalence classes $[g,\alpha]$,
where $g$ is a Riemannian metric and $\alpha\in\Omega^{n-1}(M\setminus\{x_0\};\mathcal O_M^\C)$ is such that
$\ec(g)=d\alpha$ on $M\setminus\{x_0\}$. Two such pairs $[g_1,\alpha_1]$ and $[g_2,\alpha_2]$ are equivalent
iff $\alpha_2-\alpha_1=\cs(g_1,g_2)$ mod coboundaries, see \cite[Section~3.2]{BH03}.
Every homotopy class of paths connecting
$x_0$ and $x_1$ provides an identification between
$\Eul^*_{x_0}(M;\C)$ and $\Eul^*_{x_1}(M;\C)$. Again, one can consider the
set of coEuler structures as a flat principal bundle $\Eul^*(M;\C)$ over $M$
with structure group $H^{n-1}(M;\mathcal O_M^\C)$. Its fiber over $x_0$ is $\Eul_{x_0}^*(M;\C)$, and its holonomy is given by
the composition
$$
\pi_1(M)\to H_1(M;\Z)\xrightarrow{-\chi(M)}H_1(M;\Z)\to H_1(M;\C)\to H^{n-1}(M;\mathcal O_M^\C)
$$
where the last arrow indicates Poincar\'e duality.

The affine version of Poincar\'e duality introduced in Section~\ref{S:eul}
can be consider as a parallel isomorphism of flat principal bundles over $M$
\begin{equation}\label{E:A}
P:\Eul(M;\C)\to\Eul^*(M;\C)
\end{equation}
which is equivariant over the homomorphism of structure groups
$H_1(M;\C)\to H^{n-1}(M;\mathcal O_M^\C)$ provided by Poincar\'e duality. We have
$P([X,c])=[g,\alpha]$ iff
$$
\int_{M\setminus(\mathcal X\cup\{x_0\})}\omega\wedge(X^*\Psi(g)-\alpha)=\int_c\omega
$$
for all closed one forms $\omega$ which vanish in a neighborhood of $\mathcal X\cup\{x_0\}$.

If $E$ is a flat complex vector bundle over $M$ we consider the flat line
bundle
$$
\Det(M;E):=\det H^*(M;E)\otimes(\det E)^{-\chi(M)}.
$$
Let $\Det^\times(M;E)$ denote its frame bundle, a flat principal bundle
over $M$ with structure group $\C^\times$ and holonomy given by
$$
\pi_1(M)\to H_1(M;\Z)\xrightarrow{(\theta_E)^{\chi(M)}}\C^\times.
$$
We will also consider the flat principal bundle $\Det^\times(M;E)^{-2}$ over
$M$ with structure group $\C^\times$ and holonomy given by the composition
$$
\pi_1(M)\to H_1(M;\Z)\xrightarrow{(\theta_E)^{-2\chi(M)}}\C^\times.
$$
Note that elements in $\Det^\times(M;E)^{-2}$ can be considered as 
non-degenerate bilinear forms on the corresponding fiber of $\Det(M;E)$.

The combinatorial torsion defines a parallel homomorphism of flat principal bundles
\begin{equation}\label{E:B}
\tau_E^\comb:\Eul(M;\Z)\to\Det^\times(M;E)^{-2}
\end{equation}
which is equivariant over the homomorphism of structure groups
$$
(\theta_E)^2:H_1(M;\Z)\to\C^\times.
$$
This formulation encodes in a rather natural way the combinatorial torsion's
dependence on the Euler structure and its base point.
Concerning the definition of \eqref{E:B}, recall that the corresponding construction in Section~\ref{S:comb}
assigns to an Euler structure $\e_{x_0}\in\Eul_{x_0}(M;\Z)$ and a bilinear form $b_{x_0}$ on $E_{x_0}$ a 
bilinear form on $\det H^*(M;E)$. Tensorizing this with the 
bilinear form on $(\det E_{x_0})^{-\chi(M)}$ induced by $b_{x_0}$, we obtain an element of 
$\Det_{x_0}^\times(M;E)^{-2}$ which does not depend on the choice of $b_{x_0}$. By definition this is the 
combinatorial torsion $\tau^\comb_E(\e_{x_0})$ in \eqref{E:B}.

If $b$ is a fiber wise non-degenerate symmetric bilinear form on $E$, its
analytic torsion provides a parallel homomorphism of flat principal bundles
\begin{equation}\label{E:C}
\tau^\an_{E,[b]}:\Eul^*(M;\C)\to\Det^\times(M;E)^{-2}
\end{equation}
which is equivariant over the homomorphism of structure groups
$$
H^{n-1}(M;\mathcal O_M^\C)\to\C^\times,
\quad\beta\mapsto\Bigl(e^{\langle[\omega_{E,b}]\cup\beta,[M]\rangle}\Bigr)^2.
$$
The definition of \eqref{E:C} is essentially the same as in Section~\ref{S:anator}. To be more precise,
we represent the coEuler structure $\e_{x_0}^*\in\Eul^*_{x_0}(M;\C)$ as $\e^*_{x_0}=[g,\alpha]$, where
$\alpha\in\Omega^{n-1}(M\setminus\{x_0\};\mathcal O_M^\C)$ is such that $\ec(g)=d\alpha$. We write 
$b_{(\det E_{x_0})^{-\chi(M)}}$ for the induced bilinear form on $(\det E_{x_0})^{-\chi(M)}$, and set
$$
\tau_{E,g,b,\alpha}^\an:=\tau^\an_{E,g,b}(0)\cdot\prod_q\bigl({\det}'(\Delta_{E,g,b,q})\bigr)^{(-1)^qq}\cdot
\exp\Bigl(-2\int_M\omega_{E,b}\wedge\alpha\Bigr)\otimes b_{(\det E_{x_0})^{-\chi(M)}}.
$$
If $\chi(M)\neq0$, then $\alpha$ will be singular at $x_0$ and the integral
$\int_M\omega_{E,b}\wedge\alpha$ has to be regularized, see \cite{BH03, BH06}. Due to this regularization
the additional term $\chi(M)\tr(b_u^{-1}\dot b_u)(x_0)$ will appear on the right hand side of 
\eqref{E:corr} and cancel the variation of $b_{(\det E_{x_0})^{-\chi(M)}}$.
Other than that the proof of Theorem~\ref{T:anom} remains the same. Thus $\tau_{E,g,b,\alpha}$
depends on $E$, $\e_{x_0}^*$ and $[b]$ only. By definition this is the analytic torsion
$\tau^\an_{E,[b]}(\e^*_{x_0})$ in \eqref{E:C}.

In this language the extension of Conjecture~\ref{C:main} to non-vanishing
Euler--Poincar\'e characteristics asserts that
for all $b$ we have
$$
\tau^\an_{E,[b]}\circ P\circ\iota=\tau^\comb_E
$$
as an equality of homomorphism of principal bundles over $M$, see \eqref{E:D}, \eqref{E:A},
\eqref{E:B} and \eqref{E:C}.

As in Section~\ref{S:BZ} one defines the relative torsion as the quotient of analytic
and combinatorial torsion. This is a non-vanishing complex number independent of the
Euler structure and its base point. Its properties in Proposition~\ref{P:S} remain true as stated.
With little more effort one shows that the relative torsion in general is given by the
formula in Proposition~\ref{P:Sexp}. Proving the generalization of
Conjecture~\ref{C:main} thus amounts to show that the right hand side of the equation
in Proposition~\ref{P:Sexp} equals $1$, even if $\chi(M)\neq0$.
In view of the anomaly formula it suffices to check this for a single Riemannian metric and 
any representative of the homotopy class $[b]$.

\end{document}